\RequirePackage{amsthm}

\documentclass[sn-mathphys-num]{sn-jnl}

\usepackage{graphicx}%
\usepackage{multirow}%
\usepackage{amsmath,amssymb,amsfonts}%
\usepackage{mathrsfs}%
\usepackage[title]{appendix}%
\usepackage{xcolor}%
\usepackage{textcomp}%
\usepackage{manyfoot}%
\usepackage{booktabs}%
\usepackage{algorithm}%
\usepackage{algorithmicx}%
\usepackage{algpseudocode}%
\usepackage{listings}%
\usepackage{tikz-cd}
\usepackage{adjustbox}
\usepackage{todonotes}
\usepackage{mathtools}

\theoremstyle{thmstyleone}%
\theoremstyle{thmstyletwo}%
\newtheorem{remark}{Remark}%

\theoremstyle{thmstylethree}%
\newcommand{\RR}{\mathbb{R}} 

\DeclareMathOperator{\im}{im} %

\begin{document}

\title{Persistent de Rham-Hodge Laplacians in Eulerian representation for manifold topological learning}

\author{Zhe Su$^{1}$, Yiying Tong$^{2}$,
\footnote{Corresponding author.		Email: ytong@msu.edu}
and Guo-Wei Wei$^{1,3,4}$ 
\footnote{Corresponding author.		Email: weig@msu.edu}
\\
$^{1}$Department of Mathematics,\\ Michigan State University, East Lansing, MI 48824, USA\\
$^{2}${Department of Computer Science and Engineering,\\ Michigan State University,
East Lansing, MI 48824, USA}\\
$^{3}$Department of Biochemistry and Molecular Biology,\\ 
Michigan State University, East Lansing, MI 48824, USA\\
$^{4}$Department of Electrical and Computer Engineering, \\
Michigan State University, East Lansing, MI 48824, USA}

\abstract{
Recently, topological data analysis has become a trending topic in data science and engineering. However, the key technique of topological data analysis, i.e., persistent homology, is defined on point cloud data, which does not work directly for data on manifolds. Although earlier evolutionary de Rham-Hodge theory deals with data on manifolds, it is inconvenient for machine learning applications because of the numerical inconsistency caused by remeshing the involving manifolds in the Lagrangian representation. In this work, we introduce persistent de Rham-Hodge Laplacian, or persistent Hodge Laplacian (PHL) as an abbreviation, for manifold topological learning. Our PHLs are constructed in the Eulerian representation via structure-persevering Cartesian grids, avoiding the numerical inconsistency over the multiscale manifolds. To facilitate the manifold topological learning, we propose a persistent Hodge Laplacian learning algorithm for data on manifolds or volumetric data. As a proof-of-principle application of the proposed manifold topological learning model, we consider the prediction of protein-ligand binding affinities with two benchmark datasets. Our numerical experiments highlight the power and promise of the proposed method. 

}



\maketitle

\section{Introduction}\label{sec.introduction}

Recent years have witnessed a fast growth of topological data analysis (TDA) in data science and engineering \cite{wasserman2018topological}. The growth is driven by the great promise of topological approaches to real-world data that are distinguished from any other statistical, mathematical, physical, and engineering methods~\cite{mischaikow2013morse,carlsson2009topology}. Typically, TDA offers a multi-scale topological characterization of data, which is the case with persistent homology \cite{edelsbrunner2008persistent, zomorodian2005computing}, a key method employed in TDA. A major feature of persistent homology is its multi-scale analysis, which creates a family of topological spaces from the original data to track the topological persistence, i.e., the lifespan of topological invariants across scales \cite{ghrist2008barcodes,bubenik2015statistical}. The other major feature of persistent homology is its topological description of a space (like connected components, loops, and voids) in terms of topological invariants, such as Betti numbers.  As such, persistent homology-based TDA  leads to much topological simplification of the geometric information in the data \cite{dey2014computing,adams2017persistence}. Consequently, TDA typically works extremely well for data with intricate complexity \cite{xia2014persistent,townsend2020representation}. Unfortunately, for data without geometric complexes, TDA may give rise to an oversimplification of key geometric characteristics, leading to a less competitive approach. 

For many years, persistent homology has been used in qualitative analysis, which is somewhat counterintuitive and unproductive for nonexperts. The power of persistent homology was not demonstrated until it was utilized in quantitative and predictive analysis via machine learning algorithms \cite{macpherson2012measuring,cang2015topological}. 
Topological deep learning (TDL), coined in 2017 \cite{cang2017topologynet},  was introduced to deal with large and intrinsically complex datasets using both persistent homology and deep neural networks. More recently, simplicial neural networks and other topological neural techniques have been applied in TDL to the design of neural network architecture.  TDL has become an emerging paradigm in data science and machine learning \cite{papamarkou2024position}. However, an increasing concern associated with this rising popularity is whether TDL brings any practical benefit beyond its mathematical elegance. There are many applications where TDL has demonstrated superiority to other competitive methods \cite{nguyen2020review}. Perhaps some of the most compelling examples are TDL's dominant wining of D3R Grand Challenges, an annual worldwide competition series in computer-aided drug design \cite{nguyen2020mathdl, nguyen2019mathematical}, its discovery of the mechanisms of severe acute respiratory syndrome coronavirus 2 (SARS-CoV-2) evolution \cite{chen2020mutations,wang2021mechanisms}, and its successful forecast of emerging dominant SARS-CoV-2 variants BA.2 \cite{chen2022omicronBA2} and BA.4/BA.5 about two months in advance \cite{chen2022persistent}.     

It is interesting to understand why TDL (or TDA) was so successful in the aforementioned examples, but was not competitive in many other situations in the literature \cite{pun2018persistent}. First, biomolecular data, which is intricately complex in their internal structures\cite{xia2014persistent}, was involved in the above compelling examples. As such, topological simplification was a productive process, whereas TDL leads to the severe loss of crucial geometric information in many other data that is relatively simple in their internal structures. Additionally, it was element-specific persistent homology, rather than the plain persistent homology, that was applied in the above examples. This approach captures physical and biological interactions in the biomolecular data \cite{cang2017topologynet}.   
In fact, in the forecast of emerging dominant SARS-CoV-2 variants BA.4/BA.5, persistent Laplacian, rather than persistent homology,  was utilized.  This happens because persistent homology has many drawbacks or limitations \cite{wei2023persistentb}. 
First, the topological invariant extracted from persistent homology is qualitative, rather than quantitative. For example, the barcode from persistent homology does not distinguish a five-number from a six-number ring.  
Additionally, persistent homology is incapable of dealing with different elements in a point cloud, which is ineffective with the physics and chemistry of (bio)molecular data.  Moreover, persistent homology cannot describe non-topological changes, i.e., homotopic shape evolution during the multi-scale (or filtration) analysis.
Further, persistent homology is incapable of handling directed networks and digraphs, such as polarization, regulation, and control issues in applications. Finally, persistent homology is unable to characterize structured data, e.g., hypergraphs, directed networks, etc. 
These challenges call for innovative new topological methods. 

To address these challenges, the persistent spectral graph, also known as persistent combinatorial Laplacian or persistent Laplacian (PL), was introduced in 2019 \cite{wang2020persistent}.   The harmonic spectra of PLs fully recover the topological invariants of persistent homology. However, the nonharmonic spectra of PLs capture the homotopic shape evolution during the multi-scale analysis that cannot be observed with persistent homology. Computational algorithms \cite{wang2021hermes,dong2024faster} and mathematical analysis \cite{liu2023algebraic, memoli2022persistent} of PLs have been reported. 
In the past few years, much effort has been given to extend persistent Laplacian to further address other limitations of persistent homology \cite{gulen2023generalization}, leading to  
persistent sheaf Laplacians~\cite{wei2024persistent},
persistent path Laplacians, 
persistent hypergraph and hyperdigraph Laplacians~\cite{liu2021persistent},
persistent directed flag Laplacians,
persistent Mayer Laplacians,
and persistent interaction Laplacians
\cite{wei2023persistentb}.
PLs have been shown to outperform persistent homology in many applications \cite{chen2022persistent,meng2021persistent}.  

However, defined on point cloud data, neither persistent homology nor PL can directly deal with two other commonly occurring data formats, namely, data on manifolds \cite{chen2021evolutionary}, such as electron density \cite{yang1984electron}, cryogenic electron microscopy density, and computed tomography images \cite{chen2017low}, and curves embedded in the three-dimensional Euclidean space, such as knots, links, and tangles, and their generalizations \cite{khovanov2000categorification,panagiotou2019topological}. Multi-scale Gauss link integral \cite{shen2024knot} and evolutionary Khovanov homology have been proposed to deal with embedded curve data \cite{shen2024evolutionary}. Evolutionary Khovanov homology integrates algebraic topology, geometric topology, and metric analysis for the first time. However, effective computational algorithms are needed for this approach to be widely used in practical applications.  

To carry out manifold topological analysis of data on manifolds, the evolutionary de Rham-Hodge method was introduced  \cite{chen2021evolutionary}. This approach creates a family of multi-scale manifolds with boundaries from a given data and then builds evolutionary Hodge Laplacian operators on the multi-scale manifolds with appropriate boundary conditions. While originated from sharply different topological spaces,   evolutionary Hodge Laplacian and PLs share the same algebraic structure and capture topological invariants in their harmonic spectra \cite{ribando2024graph}.  Case studies have been given to demonstrate evolutionary de Rham-Hodge theory-based manifold topological analysis of data on manifolds \cite{chen2021evolutionary}.  However, this approach was based on discrete exterior calculus \cite{desbrun2006discrete,dodziuk1976finite} or finite element exterior calculus \cite{arnold2006finite} in the Lagrangian representation, which is not efficient for multi-scale analysis and machine learning studies. Specifically, the regeneration of the evolving manifolds at different scales with different Lagrangian meshes causes numerical inconsistencies and becomes expensive for practical applications in machine learning studies.  This challenge calls for new effective manifold topological analysis approaches for data on manifolds.   

The objective of this work is to develop a persistent de Rham-Hodge theory on the Euler representation for manifold topological learning (MTL). To this end, we solve Hodge Laplacians on a pre-designed structure-persevering Cartesian grid for all scales to avoid numerical inconsistency. We construct a natural mapping of differential forms from a manifold with boundary embedded in $\RR^3$ to a large manifold,  use it to produce persistent cohomology mapping, and construct a persistent Hodge Laplacian with built-in boundary conditions.
Our new approach draws on differential geometry,  algebraic topology, partial differential equations, metric analysis, and numerical analysis. To give a proof-of-principle demonstration, we pair the proposed persistent de Rham-Hodge Laplacians with an effective machine learning algorithm to predict protein-ligand binding affinities. Based on two benchmark datasets in the Protein Data Bank (PDB), PDBbind v2007 and  PDBbind v2016, we show that our MTL model gives rise to cutting-edge performance. 

The rest of this paper is organized as follows: Section 2 offers a primer on the de Rham-Hodge theory on manifolds with boundaries; Section 3  presents our discretization for evolutionary de Rham-Hodge theory based on spectrum calculation of Laplacians associated with sublevel sets on Cartesian grids; Section 4 presents our construction for persistent de Rham-Hodge Laplacians both in the continuous setting and for given level set functions on Cartesian grids; Section 5 showcases preliminary studies on the applications of MTL; and Section 6 concludes the paper.

\section{De Rham-Hodge Theory}\label{sec.deRhamHodge}

Let $M$ be an $m$-dimensional smooth, orientable, compact Riemannian manifold with boundary. Denote by $\Omega^k(M)$ the space of all differential $k$-forms on $M$, i.e., the space of all smooth antisymmetric covariant tensor fields on $M$ of degree $k$. The \emph{differential} $d$, also called exterior derivative, is the unique $\RR$-linear mapping from the space of $k$-forms $\Omega^k(M)$ to the space of $(k\!+\!1)$-forms $\Omega^{k+1}(M)$ satisfying the Leibniz rule with respect to the wedge product $\wedge$ and the nilpotent property $dd = 0$.
A key property of differential forms is that they can be integrated over any orientable $k$-submanifolds of $M$. For any oriented $(k\!+\!1)$-submanifold $S\subset M$ with boundary $\partial S$, Stokes' theorem, as a generalization of the Newton-Leibniz rule, states that the integral of a differential $k$-form $\omega$ over $\partial S$ is equal to the integral of its differential over $S$, i.e., 
\begin{align}\label{eq.stokesthm}
	\int_S d\omega= \int_{\partial S}\omega.
\end{align}
The differential $d$ generalizes and unifies the classical operators in vector calculus, such as gradient $\nabla$, curl $\nabla\times$, and divergence $\nabla\cdot$ in $\RR^2$ and $\RR^3$. For instance, in $\RR^3$, $0$-forms and $3$-forms can be identified with scalar fields, while $1$-forms and $2$-forms can be identified with vector fields. In this case, the differential $d$ corresponds to the gradient operator $\nabla$ when applied to $0$-forms, the curl operator $\nabla\times$ when applied to $1$-forms, or the divergence operator $\nabla\cdot$ when applied $2$-forms. The nilpotent property $dd=0$ directly leads to the vector field analysis identities $\nabla\times\nabla=0$ and $\nabla\cdot\nabla\times=0$.

A differential form $\omega\in\Omega^k(M)$ is called \emph{closed} if $d\omega = 0$, or \emph{exact} if there is a $(k\!-\!1)$-form $\zeta\in\Omega^{k-1}(M)$ such that $\omega = d\zeta$. Due to the property $dd=0$, every exact form is closed. Thus, the differential $d$ links the sequence of the spaces of differential forms on $M$ into a chain complex
\begin{align}
	0\overset{}{\longrightarrow}\Omega^0(M) \overset{d}{\longrightarrow} \Omega^1(M) \overset{d}{\longrightarrow} \cdots
	\overset{d}{\longrightarrow} \Omega^{m-1}(M) \overset{d}{\longrightarrow} \Omega^m(M) \xrightarrow{} 0.
\end{align}
The $k$-th \emph{de Rham cohomology} group, denoted by $H^k_{dR}(M)$, is then defined to be the $k$-th homology of this chain complex, i.e., the quotient space of closed $k$-forms modulo the space of exact $k$-forms, i.e.,
\begin{align}
	H^k_{dR}(M) = \frac{\ker(d:\Omega^k(M)\to\Omega^{k+1}(M))}{\im(d:\Omega^{k-1}(M)\to\Omega^k(M))}.
\end{align}
The de Rham cohomology, by the de Rham theorem, is naturally isomorphic to the singular cohomology, and thus depends only on the manifold topology.

Let $g$ be a Riemannian metric on $M$ and $\langle\cdot,\cdot\rangle_g$ be the pointwise inner product induced by $g$ on $\Omega^k(M)$. The \emph{Hodge star} operator $\star$ provides an isomorphism from the space of differential $k$-forms $\Omega^k(M)$ to the space of $(m\!-\!k)$-forms $\Omega^{m-k}(M)$, defined by the following formula
\begin{equation}\label{eq.innerprod.p}
	\omega\wedge\star\eta = \langle\omega,\eta\rangle_g\;\mu_g,
\end{equation}
where $\mu_g$ is the volume form on $M$ induced by $g$. The Hodge $L^2$-inner product on the space of $k$-forms $\Omega^k(M)$ can then be obtained by taking the integral of the formula \eqref{eq.innerprod.p}
\begin{equation}\label{eq.l2innerprod.forms}
	(\omega, \eta) = \int_M\omega\wedge\star\eta.
\end{equation}
The \emph{codifferential} $\delta: \Omega^k(M)\to\Omega^{k-1}(M)$ is defined by
\begin{equation}\label{eq.codifferential}
	\delta = (-1)^{m(k-1)+1}\star d\star,
\end{equation}
which also has the nilpotent property $\delta\delta = 0$. We call a differential form $\omega\in\Omega^k(M)$ \emph{coclosed} if $\delta\omega = 0$, or \emph{coexact} if there is a $(k+1)$-form $\eta\in\Omega^{k+1}(M)$ such that $\omega = \delta\eta$. 
The codifferential $\delta$, as the differential $d$, also extends the classical gradient, curl and divergence in vector calculus. In $\RR^3$, it corresponds to $-\nabla\cdot$, $\nabla\times$ and $-\nabla$ when applied to $1$-forms, $2$-forms and $3$-forms, respectively.

The \emph{Hodge Laplacian} for differential forms is defined as $\Delta = d\delta + \delta d: \Omega^k(M)\to\Omega^k(M)$. Its kernel, consisting of all differential $k$-forms $\omega$ on $M$ with $\Delta\omega = 0$, is called the space of \emph{harmonic} $k$-forms.
We denote by $\mathcal{H}^k_{\Delta}(M)$ the space of harmonic $k$-forms and by $\mathcal{H}^k(M)$ the space of $k$-forms that are both closed and coclosed, i.e., $\mathcal{H}^k(M) = \ker d\cap \ker\delta$. The latter space $\mathcal{H}^k(M)$, known as the space of harmonic $k$-fields, is in general only a subset of the space of harmonic forms $\mathcal{H}^k(M)\subset\mathcal{H}^k_{\Delta}(M)$, and is infinite-dimensional \cite{schwarz2006hodge}. However, in the case of closed manifolds where $\partial M=\emptyset$, the space of harmonic forms $\mathcal{H}^k_{\Delta}(M)$ reduces to the space $\mathcal{H}^k(M),$ as any harmonic form is both closed and coclosed. The result follows directly from the following formula
\begin{align}\label{eq.identification.Laplacian}
	0 = (\Delta\omega, \omega) = ((d\delta + \delta d)\omega, \omega) =  (d\omega, d\omega) + (\delta\omega, \delta\omega),
\end{align}
due to the $L^2$-adjointness of the codifferential $\delta$ and the differential $d$ on closed manifolds, i.e., $(d\omega, \eta) = (\omega, \delta\eta)$.

The classical Hodge decomposition theorem for closed manifolds states that the space of differential $k$-forms $\Omega^k(M)$ can be decomposed as
\begin{align}\label{eq.hd.closedMfld}
	\Omega^k(M)= d\Omega^{k-1}(M)\oplus\delta\Omega^{k+1}(M) \oplus \mathcal{H}^k_{\Delta}(M).
\end{align}
These three subspaces are mutually orthogonal with respect to the inner product \eqref{eq.l2innerprod.forms}. Moreover, Hodge theorem identifies the harmonic space $\mathcal{H}^k_{\Delta}(M)$ with the $k$-th de Rham cohomology group $H^k_{dR}(M)$, which states that each harmonic form corresponds to exactly one equivalence class in $H^k_{dR}(M)$. Therefore, the harmonic space $\mathcal{H}^k_{\Delta}(M)$ is fully determined by the manifold topology, and is finite-dimensional with its dimension given by the Betti number $\dim\mathcal{H}^k_{\Delta}(M) = \beta_k$.

%


\subsection{Hodge decomposition for manifolds with boundary}
In the presence of a non-empty boundary $\partial M$, the two operators $d$ and $\delta$ are not $L^2$-adjoint, as integration by parts leads to
\begin{align}
	(d\omega, \eta) = (\omega, \delta\eta) + \int_{\partial M}\omega\wedge\star\eta,
\end{align}
which contains a boundary term that may not vanish, and thus the decomposed subspaces in \eqref{eq.hd.closedMfld} are not orthogonal. However, certain boundary conditions can be enforced, ensuring the adjointness of the differential $d$ and the codifferential $\delta$, thereby inducing an orthogonal decomposition of the space of differential forms. 

The most common choices of boundary conditions ensuring the adjointness of $d$ and $\delta$ are the {normal} (Dirichlet) and {tangential} (Neumann) boundary conditions. 
A differential form $\omega\in\Omega^k(M)$ is called \emph{normal} (Dirichlet) if it gives zero when applied to tangent vectors of the boundary, or \emph{tangential} (Neumann) if the same holds for its dual $\star\omega$ instead. Denote by $\Omega^k_n(M)$ the set of normal differential $k$-forms and by $\Omega^k_t(M)$ the set of tangential differential forms, i.e.,
\begin{align}
	\Omega^k_n(M) &= \{\omega\in\Omega^k(M)\, \vert\quad \omega\vert_{\partial M} = 0\}\\
	\Omega^k_t(M) &= \{\omega\in\Omega^k(M)\, \vert\quad \star\omega\vert_{\partial M} = 0\}.
\end{align}
Following their definitions, the spaces $\Omega^k_n(M)$ and $\Omega^{m-k}_t(M)$ are isomorphic under the Hodge star operator $\star$, also known as the Hodge duality. Moreover, the differential $d$ preserves the normal boundary conditions, while the codifferential $\delta$ preserves the tangential boundary conditions.

The Hodge-Morrey decomposition~\cite{morrey1956variational} states that there is a 3-component $L^2$-orthogonal decomposition
\begin{align}\label{eq.morreyDecomp}
	\Omega^k(M) = d\Omega^{k-1}_n(M)\oplus\delta\Omega^{k+1}_t(M) \oplus \mathcal{H}^k(M),
\end{align}
The orthogonality of the decomposition directly comes from the adjointness of $\delta$ and $d$ when enforcing the normal or tangential boundary conditions. For $\omega\in\Omega^k(M),$ there is a unique decomposition of $\omega$ given as follows:
\begin{align}\label{eq.morreyDecomp.omega}
	\omega = d\alpha_n + \delta\beta_t + \eta,
\end{align}
where $\alpha_n\in \Omega^{k-1}_n(M)$, $\beta_t\in\Omega_t^{k+1}(M)$, and $\eta\in\mathcal{H}^k(M)$. Note that the potentials $\alpha_n$ and $\beta_t$ are not uniquely determined as all $\alpha_n+d\eta$ and $\beta_t+\delta\gamma$ with any $\eta \in\Omega_n^{k-2}(M)$ and
$\gamma\in\Omega_t^{k+2}(M)$ serve as potentials for the same components. However, the issue can be addressed by enforcing \emph{gauge} conditions, such as
\begin{align} 
	\label{eq.normal.gauge}\delta\alpha_n = 0,\\
	\label{eq.tangential.gauge}d\beta_t = 0.
\end{align}
The potentials $\alpha_n$ and $\beta_t$ can then be uniquely determined by the following equations
\begin{align}
	\begin{cases}
		\Delta\alpha_n = \delta\omega\\
		\Delta\beta_t = d\omega,
	\end{cases}	
\end{align}
by resolving the (finite) rank deficiencies of $\Delta$ under these boundary conditions. 

\begin{remark}
	In the case that $M$ is a closed manifold, i.e., $\partial M = \emptyset$, both the spaces $\Omega^k_n(M)$ and $\Omega^k_t(M)$ coincide with the space of differential forms $\Omega^k(M)$, and the space of harmonic fields is identical to the space of harmonic forms $\mathcal{H}^k(M)=\mathcal{H}^k_{\Delta}(M)$. The Hodge decomposition \eqref{eq.morreyDecomp} then reduces to the classical Hodge decomposition \eqref{eq.hd.closedMfld} for closed manifolds.
\end{remark}

\begin{remark}
	The Hodge-Morrey decomposition \eqref{eq.morreyDecomp} in the low dimensional Euclidean spaces $\RR^2$ and $\RR^3$, often referred to as the Helmholtz-Hodge decomposition in vector calculus, states that any vector field $\mathbf{v}$ defined on a compact domain can be orthogonality decomposed as
	\begin{align}\label{eq.hhd.3d}
		\mathbf{v} = \nabla f + \nabla\times \mathbf{u} + \mathbf{h},
	\end{align}
	where $f$ is a scalar potential that vanishes on the boundary of the domain, $\mathbf{u}$ is a vector field orthogonal to the boundary, and $\mathbf{h}$ is the harmonic vector field satisfying $\nabla\times \mathbf{h} = 0$ and $\nabla\cdot \mathbf{h} = 0$.
	The first component $\nabla f$ and the second component $\nabla\times\mathbf{u}$ are often called the curl-free and divergence-free parts of the vector field $\mathbf{v}$ respectively. Note that in the presence of a boundary, the resulting scalar potential $f$ is also called satisfying the normal boundary of $0$-forms, and the vector field $\mathbf{u}$ is called satisfying the tangential boundary condition of $2$-forms, which are direct counterparts of the potentials $\alpha_n$ and $\beta_t$ in \eqref{eq.morreyDecomp}. For a complete correspondence between scalar or vector fields, and differential forms under the normal and tangential boundary conditions, see \cite{zhao20193d}.
\end{remark}

The space of harmonic fields $\mathcal{H}^k$, in general, is infinite-dimensional, and thus has no direct correspondence with the cohomology of the manifold. However, as noted by \cite{zhao20193d}, one can restrict to the space of normal harmonic fields, namely $\mathcal{H}^k_n(M) = \mathcal{H}^k(M)\cap\Omega^k_n(M)$, and the space of tangential harmonic fields, $\mathcal{H}^k_t(M)=\mathcal{H}^k(M)\cap\Omega^k_t(M)$. As a consequence of the de Rham map, these two subspaces $\mathcal{H}^k_n(M)$ and $\mathcal{H}^k_t(M)$ are fully determined by the topology of $M$: the space of normal harmonic fields $\mathcal{H}^k_n(M)$ is isomorphic to the relative de Rham cohomology $H^k_{dR}(M,\partial M)$, while the space of tangential harmonic fields $\mathcal{H}^k_t(M)$ is isomorphic to the absolute de Rham cohomology $H^k_{dR}(M)$~\cite{friedrichs1955differential}. The two subspaces
$\mathcal{H}^k_n(M)$ and $\mathcal{H}^k_t(M)$ are thus finite-dimensional, with dimensions given by the Betti numbers: $ \dim\mathcal{H}^{k}_n(M) = \beta_{m-k}$ and $ \dim\mathcal{H}^{k}_t(M) = \beta_k$. Furthermore, the kernels of the Hodge Laplacian $\Delta$, when restricted to the space of normal forms $\Omega^k_n(M)$ and the space of tangential forms $\Omega^k_t(M)$ with gauge conditions on the boundary, can be identified to the space of normal harmonic fields and the space of tangential harmonic fields, respectively. Denote by $\Delta_{n}$ and $\Delta_{t}$ the restrictions of the Hodge Laplacian $\Delta$ on the space of normal fields $\Omega^k_n(M)$ satisfying Eq.~\eqref{eq.normal.gauge} and the space of tangential fields $\Omega^k_t(M)$ satisfying Eq.~\eqref{eq.tangential.gauge}, i.e., $\Delta_n:\Omega^k_n(M)\to\Omega^k(M)$ and $\Delta_t:\Omega^k_t(M)\to\Omega^k(M)$. Then immediately we have $\ker\Delta_n = \mathcal{H}^k(M)\cap\Omega^k_n(M)=\mathcal{H}^k_n(M)$ and $\ker\Delta_t = \mathcal{H}^k(M)\cap\Omega^k_t(M) = \mathcal{H}^k_t(M)$. The result follows directly from Eq.~\eqref{eq.identification.Laplacian}. These identifications, finally, enable us to study the topology of the underlying manifold $M$ through the Hodge Laplacians on normal and tangential forms.

\begin{remark}
	In fact, let $\mathcal{H}^k_{co} = \mathcal{H}^k(M)\cap \delta\Omega^{k+1}(M)$ and $\mathcal{H}^k_{ex} = \mathcal{H}^k(M)\cap d\Omega^{k-1}(M)$. The space of harmonic fields $\mathcal{H}^k(M)$ can be further orthogonally decomposed for smooth manifolds
	\begin{align}
		\mathcal{H}^k(M) = &\mathcal{H}^k_{co}(M)\oplus\mathcal{H}^k_n(M) \\ = &\mathcal{H}^k_{ex}(M)\oplus\mathcal{H}^k_t(M),
	\end{align}
which results in the Hodge-Morrey-Friedrichs decomposition given as follows
\begin{align}
	\Omega^k(M) &= d\Omega^{k-1}_n(M)\oplus\delta\Omega^{k+1}_t(M) \oplus \mathcal{H}^k_{co}(M)\oplus\mathcal{H}^k_n(M)\\
	& = d\Omega^{k-1}_n(M)\oplus\delta\Omega^{k+1}_t(M) \oplus \mathcal{H}^k_{ex}(M)\oplus\mathcal{H}^k_t(M).
\end{align}
In particular, if $M$ is a compact domain in Euclidean spaces, then there is a unique orthogonal $5$-component decomposition 
\begin{align}\label{eq.hd.5subspaces}
	\Omega^k(M) = d\Omega^{k-1}_n(M)\oplus\delta\Omega^{k+1}_t(M) \oplus \mathcal{H}^{k}_n(M) \oplus \mathcal{H}^{k}_t(M)\oplus (d\Omega^{k-1}(M)\cap\delta\Omega^{k+1}(M)),
\end{align}
as the spaces $\mathcal{H}^k_n(M)$ and $\mathcal{H}^k_t(M)$ are  $L^2$-orthogonal, instead of just being transversal for compact manifolds in general~\cite{shonkwiler2009poincare}. Due to the correspondence between differential forms and vector fields in the low-dimensional Euclidean spaces, the implementation of this $5$-component Hodge decomposition has been applied and implemented to the study of vector fields for surface triangle meshes, for tetrahedral meshes \cite{zhao20193d} and for regular Cartesian grids \cite{su2024hodge}.

\end{remark}


As we mainly focus on applications of compact domains in $\RR^3$, to study the geometric and topological information of the underlying manifolds, there are eight Laplacians to be considered, which are defined on the spaces of differential $k$-forms with $k= 0,1,2,3$ satisfying either the normal or the tangential boundary conditions. However, thanks to the duality between the space of normal fields and tangential fields, the study of the spectra of these eight Laplacians reduces to that of four Laplacians on one of the two types of boundary conditions, and finally to the singular spectra of three differential operators, applied to differential forms of degree $k = 0,1,2,3$~\cite{chen2021evolutionary}. Further details will be discussed in the next section for the discretization of Laplacians.

\section{Discretization and construction of Laplacians}

In this section, we elaborate on the discretization of the Hodge Laplacian and introduce the Boundary-Induced Graph (BIG) Laplacian for compact domains in low-dimensional Euclidean spaces. Although the theory works for 2D compact domains, we focus only, for the remainder of the paper, compact domains in $\RR^3$, as we target mainly 3D applications. We use DEC to discretize all differential operators and differential forms on regular Cartesian grids, as it allows for efficient and accurate numerical algorithms relying on just matrix algebra, while keeping the $L^2$ orthogonality between different components in Hodge decomposition. In addition, the constructed discrete differential operators and differential forms in DEC approximate their smooth analogs. For the characterization of the underlying manifold, we choose the Eulerian formulation, where the manifold is given as a sublevel set of a level set function defined on a regular Cartesian grid. Another common way, called the Lagrangian formulation, discretizes the manifold as simplicial meshes, i.e., triangular or tetrahedral meshes in 2D or 3D. The spectrum analysis of the Hodge Laplacians has been discussed in \cite{zhao20193d} for the Lagrangian formulation and in \cite{su2024hodge} for the Eulerian formulation. Compared to the Lagrangian case, the Eulerian representation uses vertices, edges, faces and cells all fixed in a Cartesian grid, which significantly simplifies the data structures and algorithms. The Hodge stars, in the latter case, are close to rescaled identity matrices. This fact simplifies the study of Hodge Laplacians to that of BIG Laplacians with no Hodge stars involved, and thus leads to algorithms with efficient computations.

\subsection{Discretization on entire grid}
\begin{figure}[t]
	\centering
	\includegraphics[scale=0.8]{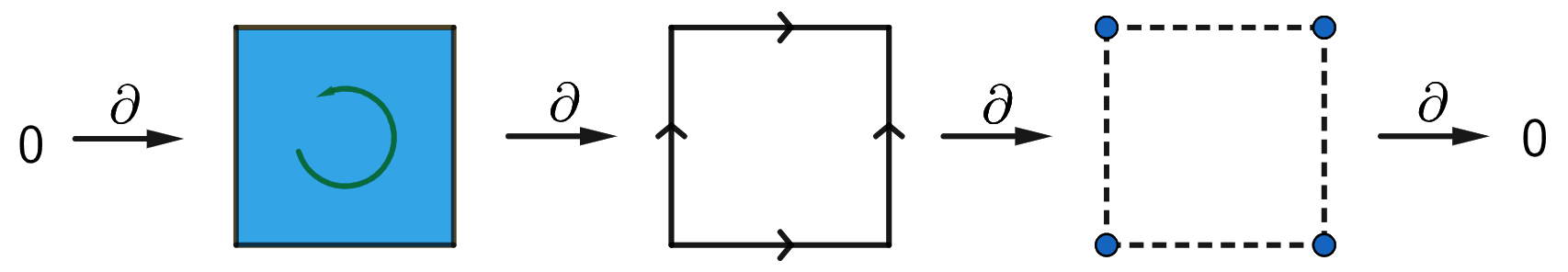}
	\caption{The chain complex of a single-cell grid formed by the boundary operator: from the face, to its edges, and to their vertices.}
	\label{fig.chainboundary}
\end{figure}

Denote by $I_m$ a rectangular $m$-dimensional regular Cartesian grid with $k$-cells oriented according to their alignments with the coordinate axes. The entire grid $I_m$ can be treated as a cell complex tessellating a rectangular domain in $\RR^m$, where each $k$-cell is a $k$-dimensional hypercube with edge length $\ell$. A continuous differential $k$-form $\omega$ on $I_m$, following the de Rham map, can be discretized by its integral value over each oriented $k$-cell $\sigma_i$, given as $W^i=\int_{\sigma_i} \omega$~\cite{desbrun2006discrete}. The discrete differential on discrete $k$-forms of the grid $I_m$ is then encoded by a sparse matrix $D^I_k$, which stores the signed incidence between $(k\!+\!1)$-cells and $k$-cells and is given as the transpose of the cell boundary operator $\partial_{k+1}^T$ on $(k\!+\!1)$-cells following from Stokes' theorem $\int_{\sigma} d\omega= \int_{\partial \sigma} \omega$. An illustration of the chain complex formed by boundary operator $\partial$ for a simple grid complex with a single 2D cell can be seen in Fig.~\ref{fig.chainboundary}, which is a straightforward generalization of the chain complex on simplicial complexes. Note that the boundary of the boundary of a cell always results in a 0 chain, i.e., $\partial\partial = 0$, whose transpose immediately produces $D^I_{k+1}D^I_k = 0$, thus preserving the nilpotent property in the continuous setting.

\begin{figure}[t]
	\centering
	\includegraphics[height=7cm, trim=0 1cm 0 0,clip]{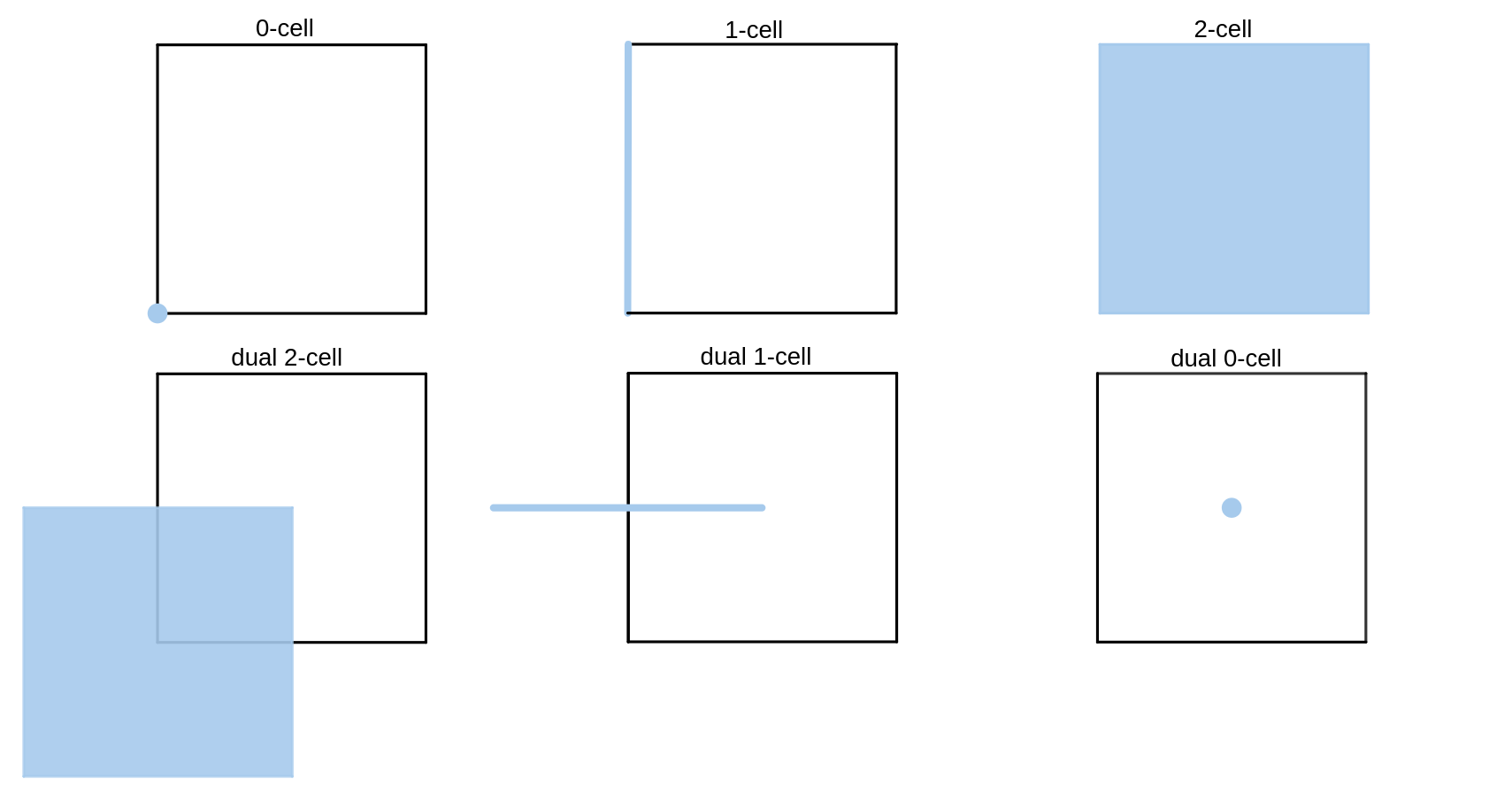}
	\caption{An example of the primal and dual grid cells for the 2D case. The top row highlights the primal cells, and the bottom row presents their corresponding dual cells.}
	\label{fig.primal.dual}
\end{figure}

The discrete Hodge star establishes a one-to-one correspondence between discrete $k$-forms on the primal grid $I_m$ and discrete $(m\!-\!k)$-forms on its dual grid, given as the translated grid with grid points located at the $m$-cell centers of $I_m$, based on the following formula
\begin{align}
	\frac{1}{|\sigma_k|}\int_{\sigma_k}\omega \approx \frac{1}{|\star\sigma_k|}\int_{\star\sigma_k}\star\omega,
\end{align}
where $\star\sigma_k$ is the dual $(m\!-\!k)$-cell formed by the dual grid points located at the centers of the primal $m$-cells incident to $\sigma_k$. See Fig.~\ref{fig.primal.dual} for an illustration for the correspondences between the primal and dual cells in the Cartesian grid case. Following from the discretization of differential forms, this correspondence leads to a diagonal matrix $S^I_k$ with diagonal entries given by the ratio between the volumes of the dual  $(m\!-\!k)$-cells and the primal $k$-cells, $\ell^{m-k}/\ell^k=\ell^{m-2k}$. The associated discrete Hodge $L^2$-inner product \eqref{eq.l2innerprod.forms} of two discrete $k$-forms $V_k$ and $W_k$ on grid $I_m$ is then given by
\begin{align}
	(V_k, W_k)^I = V_k^TS^I_kW_k.
\end{align}
The discrete codifferential, by definition of its smooth counterpart \eqref{eq.codifferential}, can be assembled from the discrete differential and Hodge star operators as $\delta^I_k = (S^I_{k-1})^{-1}D^I_{k-1}S^I_k$. Note that the discrete counterpart of the Hodge Laplacian $\Delta = d\delta + \delta d$ by replacing the differential and codifferential operators results in a nonsymmetric matrix. Instead, we consider the counterpart of $\star\Delta$ as the discrete Hodge Laplacian given by
\begin{align}
	L^I_k = (D^I_k)^TS^I_{k+1}D^I_k + S^I_kD^I_{k-1}(S^I_{k-1})^{-1}(D^I_{k-1})^TS^I_k,
\end{align}
where the operators are considered to be null for $k<0$ or $k>m$.

\subsection{Discrete differential forms and operators on $M$}

\begin{figure}[h]
	\centering
	\includegraphics[height=4cm]{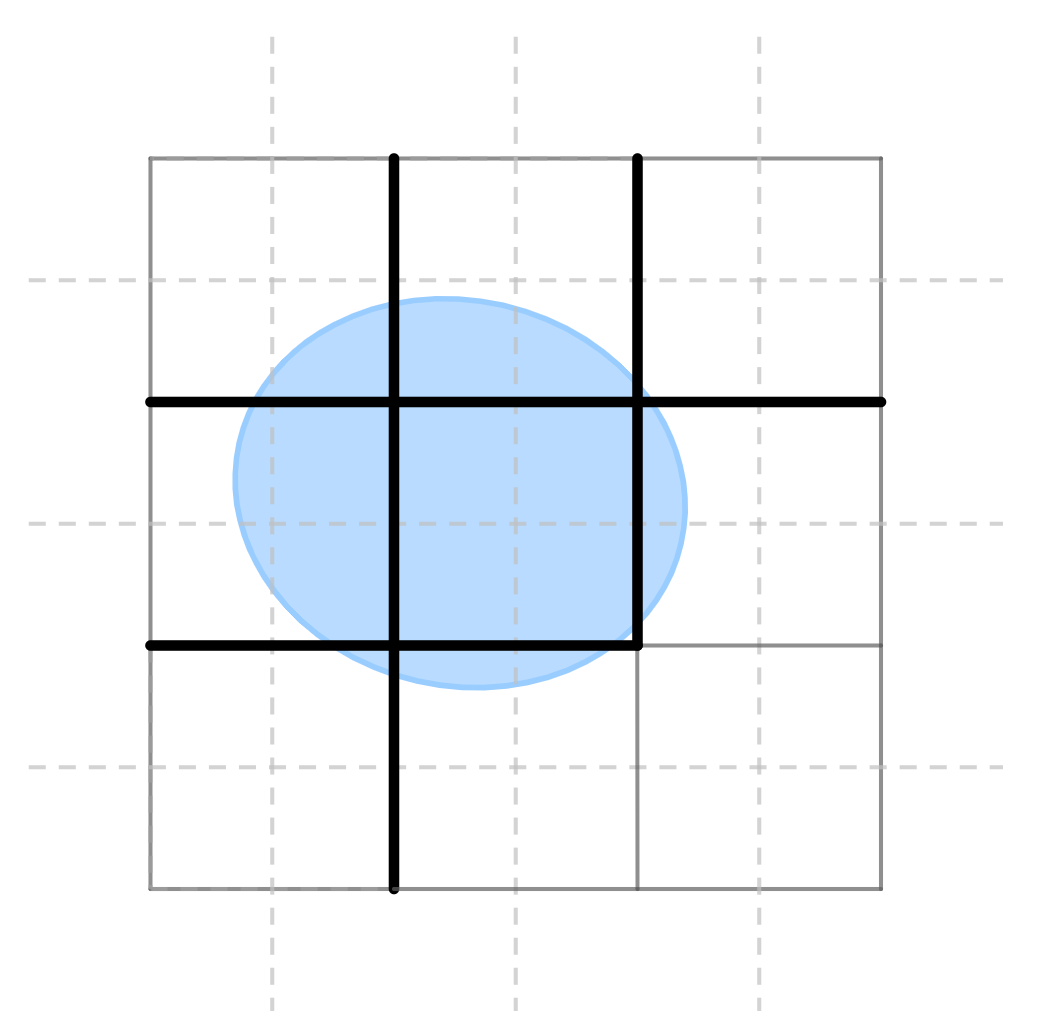}\hskip.25in
	\includegraphics[height=4cm]{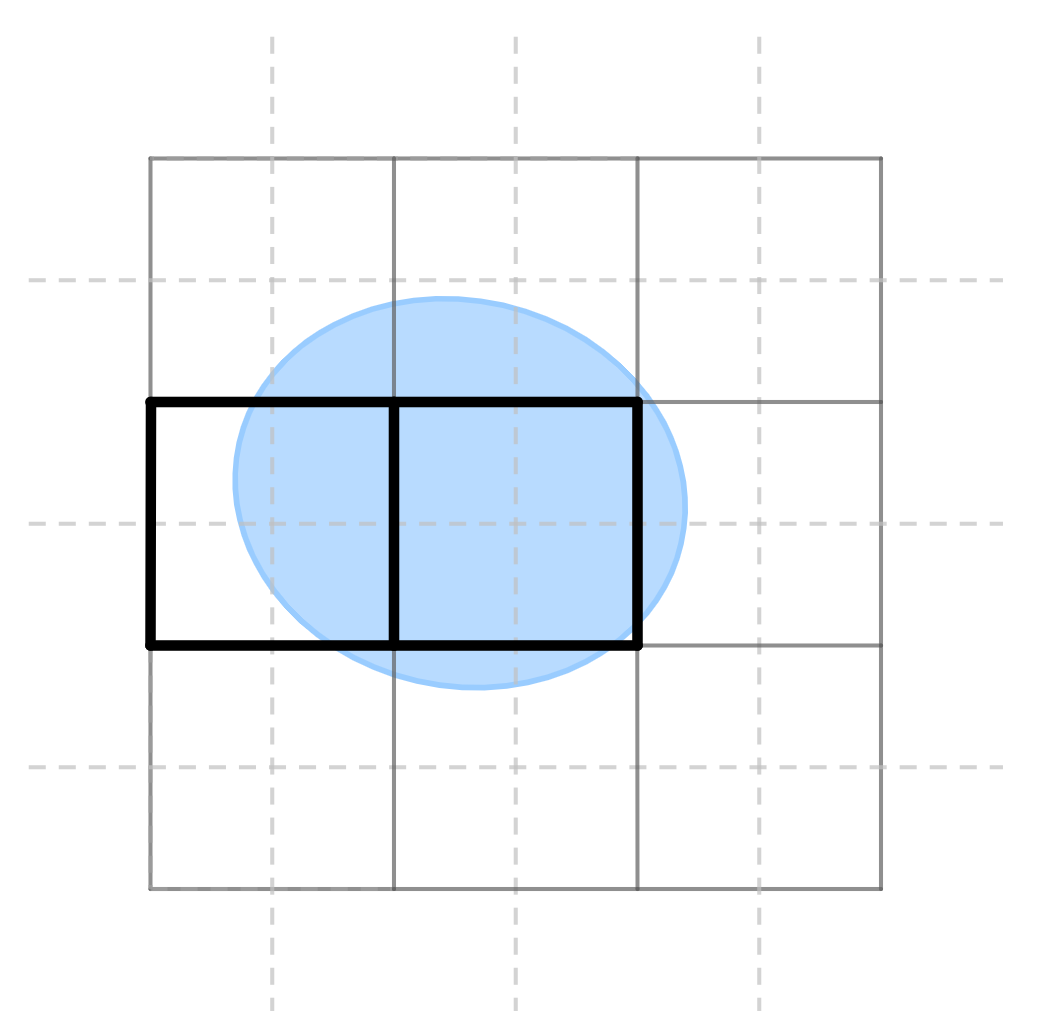}
	\caption{Distinction of normal supports (left) and tangential supports (right) for primal $1$-forms in a 2D Cartesian grid.}
	\label{fig.distinctionNTsupports}
\end{figure}
 
Compared to the case of simplicial or polygonal meshes, where the projection matrices to the interior can be straightforward to implement with the boundary elements explicitly labeled, modeling the manifold $M$ as the volume bounded by a level set surface leads to delicate computation of the projection matrices. Note that the boundary of $M$ using grid representation typically intersects with boundary $k$-cells instead of being its supersets. We restrict the computation to relevant cells by implementing the two types of boundary conditions through the inclusion or exclusion of the entire $k$-cells. We use the strategy as in \cite{su2024hodge} for the computation of projection matrices for each type of boundary condition: for the normal boundary condition, we include all cells if at least one of its vertices is inside or on the boundary of $M$, while for the tangential boundary condition, we include all cells with at least one of the vertices of the corresponding dual cells is inside or on the boundary. We refer to the former set of cells as the normal support and the latter as the tangential support. In contrast to the mesh case, it is important to note that neither the normal nor the tangential support is necessarily a superset of the other. See Fig.~\ref{fig.distinctionNTsupports} for one example showing the distinction of these two supports for $1$-forms.

In the computation of the discrete Hodge star operators, it is essential to consider and incorporate the boundary conditions. Following the procedure in \cite{su2024hodge}, we keep the dual cell volumes and adjust the primal cell volumes for normal boundary conditions, and do conversely for tangential boundary conditions with the primal cell volumes kept and the dual cell volumes changed. To be specific, when dealing with normal (resp. tangential) boundary conditions, we only compute the volume of the region of the primal (resp. dual) $k$-cells within the boundary $\partial M$ for the denominator (resp. numerator) of the ratio in the discrete Hodge star matrix, and leave the dual (resp. primal) cell volumes in the numerator (resp. denominator) unchanged. Each unaltered $k$-cell has a $k$-volume of $\ell^k$. In addition, For numerical stability, we do not alter the volume of outside primal $k$-cells, and perturb the level set function evaluated at primal/dual grid points to have an absolute value above $\epsilon=10^{-5}\ell$, which ensures well-behaved fractional $k$-volumes. We denote by $S^I_{k,n}$ and $S^I_{k,t}$ the diagonal Hodge star matrices defined on the entire grid $I^m$ corresponding to the normal and tangential boundary conditions, respectively.

The projection matrix to the corresponding support, for each type of boundary condition, can be constructed from the identity matrices by eliminating the rows corresponding to $k$-cells outside the support. Denote by $P_{k,n}$ the projection matrix for $k$-cells onto the normal support and by $P_{k,t}$ the one onto the tangential support. We then obtain a new set of differential and Hodge star operators for $M$:
\begin{align}
	D_{k,n} &= P_{k+1,n}D_kP_{k,n}^T,\quad S_{k,n} = P_{k, n}S^I_{k,n}P_{k,n}^T\\
	D_{k,t} &= P_{k+1,t}D_kP_{k,t}^T,\quad S_{k,t} = P_{k, t}S^I_{k,t}P_{k,t}^T
\end{align}
The nilpotent property $D_{k+1,n}D_{k,n} = 0$ and $D_{k+1,t}D_{k,t} = 0$ still holds for both boundary conditions due to $D^I_{k+1}D^I_k=0$ and the following observations
\begin{align}\label{eq:projectionRelation}
	P_{k+1,n}^T P_{k+1, n}D^I_k P_{k, n}^T = D^I_k P_{k, n}^T,\quad
	P_{k+1, t} D^I_k P_{k, t}^T P_{k,t} = P_{k+1,t} D^I_k.
\end{align}
The discrete Hodge $L^2$-inner products of the two types of discrete $k$-forms on the manifold $M$ for these two boundary conditions are then given by
\begin{align}
	(\xi^k,\, \zeta^k)^n &= (\xi^k)^TS_{k,n}\zeta^k \label{eq.l2innerProd.discrete.n}\\
	(\xi^k,\, \zeta^k)^t &= (\xi^k)^TS_{k,t}\zeta^k, \label{eq.l2innerProd.discrete.t}
\end{align}
whose domains are the discrete $\Omega^{k}_n(M)$ and the discrete $\Omega^{k}_t(M)$ respectively. Finally, we assemble the two types of discrete Hodge Laplacians as in the mesh case:
\begin{align}\label{eq.HodgeLaplacian.discrete}
	L_{k, n} &= D_{k, n}^TS_{k+1, n}D_{k, n} + S_{k, n}D_{k-1, n}S_{k-1, n}^{-1}D_{k-1, n}^TS_{k, n}\\
	L_{k, t} &= D_{k, t}^TS_{k+1, t}D_{k, t} + S_{k, t}D_{k-1, t}S_{k-1, t}^{-1}D_{k-1, t}^TS_{k, t}.
\end{align}
The null spaces of these discrete Hodge Laplacians, as in the continuous case, are fully determined by the topology of the underlying manifold $M$, since they only depend on the differential and projection matrices. The dimension of the kernel of $L_{k, n}$ is given by the Betti number $\beta_{m-k}$, while the dimension of the kernel of $L_{k, t}$ is given by $\beta_{k}$. Here the Betti number $\beta_k$ presents directly the number of $k$-dimensional holes on the manifold $M$. For instance, $\beta_0$ gives the number of connected components, $\beta_1$ gives the number of tunnels, and $\beta_2$ provides the number of closed cavities, respectively. The spectra of these Laplacians, in addition, could be used to study the geometric information of the manifold. It is known that the non-zero eigenvalues of the Laplacians provide rich insights into the shape of a manifold. For instance, the Fiedler value, defined as the smallest non-zero eigenvalue of a graph Laplacian, describes connectivity. As another example, the multiplicity of eigenvalues can reveal certain symmetries of the shape.

\begin{remark}
	The two types of discrete Hodge Laplacians \eqref{eq.HodgeLaplacian.discrete} not only provide rich geometrical and topological information of the underlying manifold, but also play a central role in the computation of the discrete Hodge decomposition \eqref{eq.hd.5subspaces} of differential forms for compact domains in 2D and 3D Euclidean spaces. In particular, they can be utilized, by resolving the rank deficiencies, to compute the potentials of the decomposed components in Hodge decomposition on normal or tangential support satisfying the corresponding boundary conditions. In addition, as the kernel sizes of Laplacians are finite, their eigenvectors corresponding to $0$ eigenvalues, for each $k$, form a basis for the space of normal or tangential harmonic fields.
\end{remark}

Note that the discrete Hodge stars in the Eulerian setting are almost identical to rescaled identity matrices. Therefore, the computations of the Hodge Laplacian can be further simplified by replacing the Hodge stars with identity matrices, leading to the definition of the Boundary-Induced Graph (BIG) Laplacians as follows:
\begin{align}\label{eq.BIGLaplacians}
	L^B_{k, n} &= D_{k, n}^TD_{k, n} + D_{k-1, n}D_{k-1, n}^T\\
	L^B_{k, t} &= D_{k, t}^TD_{k, t} + D_{k-1, t}D_{k-1, t}^T.
\end{align}
The BIG Laplacians were introduced in \cite{ribando2024graph} for bounded domains to facilitate the comparison and contrast of the Hodge Laplacians and the combinatorial Laplacians. They preserve the Hodge Laplacian's capability to perform differential calculus but also retain the discrete nature of combinatorial Laplacians. The convergence of the spectra of the BIG Laplacians to Hodge Laplacians has been discussed in \cite{ribando2024graph}, showing that the spectra of \eqref{eq.BIGLaplacians} converge to those of Hodge Laplacians up to a scaling value $\ell^{-2}$ when enforcing the boundary conditions. This scaling value $\ell^{-2}$ is exactly the ratio between the missing scaling factor $\ell^{m-2(k+1)}$ in $L_k$ and the missing factor $\ell^{m-2k}$ of $S_k$. As the BIG Laplacians produce results similar to those obtained from the discrete Hodge Laplacians with less computation, they can also be used to study the geometric and topological information of the underlying manifolds.

Note that the dual grid is also a Cartesian grid staggered with the primal grid by a replacement of $\ell/2$ in all three axial directions of the Cartesian coordinates. For the study of the spectra of these Laplacians, one only needs to implement one type of boundary condition, for instance, the normal boundary condition, as $L_{k,n}$ defined on the primal grid with normal boundary conditions is equivalent to $L_{m-k,t}$ defined on its dual grid with tangential boundary conditions.

\subsection{Topology-preserving construction of  Laplacians}
\label{sec.spectraLaplacians}

Note that, on the grid, the Hodge Laplacians and the BIG Laplacians are of the same sparsity patterns. For simplicity in exposition when discussing the spectrum analysis of the Laplacians, we let $L_k$ be a generic Laplacian matrix of the form
\begin{align}
	L_{k} &= D_{k}^TS_{k+1}D_{k} + S_{k}D_{k-1}S_{k-1}^{-1}D_{k-1}^TS_{k}.
\end{align}
Here the Laplacian $L_k$ can be interpreted, under choices of boundary conditions and Hodge star accuracy, as either a Hodge Laplacian, or BIG Laplacian ( with $S_k$ set to identity), under tangential or normal boundary condition. The eigenvalues and eigenvectors of $L_k$ can be solved by considering the generalized eigenvalue problem
\begin{align}
	L_k W =\lambda S_k W,
\end{align}
where $\lambda$ is an eigenvalue, and $W$ is the associated eigenvector.
To analyze the results, we perform the following transformation in the space of discrete forms: $\bar D_k = S_{k+1}^{1/2}D_kS_k^{-1/2}$, $\bar L_k = S_k^{-1/2}L_k S_k^{-1/2}$ and $\bar W = S_k^{1/2}W$. Rewriting the formulas above yields a simplified form of the Laplacian
\begin{align}
	\bar L_k = \bar{D}_k^T\bar{D}_k + \bar{D}_{k-1}\bar{D}_{k-1}^T,
\end{align}
and a regular eigenvalue problem:
\begin{align}
	\bar L_k \bar W = \lambda \bar W.
\end{align}
Note that the property $\bar{D}_{k}\bar{D}_{k-1} = 0$ is preserved.
As the non-zero eigenvalues of $\bar{D}_k^T\bar{D}_k$ and $\bar{D}_{k}\bar{D}_{k}^T$ for each $k$ are the same, given by the squared non-zero singular values of the discrete differential $\bar{D}_k$, and each Laplacian $\bar L_k$ is just the combination of $\bar{D}_k^T\bar{D}_k$ and $\bar{D}_{k-1}D_{k-1}^T$, the entire spectrum of the Laplacians can thus be studied through the singular values of discrete differentials. Let 
\begin{align}
	\bar D_k = U_{k+1}\Sigma_kV_k^T
\end{align}
be the singular value decomposition of $\bar D_k$, where $U_{k+1}$ and $V_k$ are orthogonal matrices and $\Sigma_k$ is a rectangular diagonal matrix with diagonal values given by the singular values of $\bar D_k$. It follows immediately from $\bar{D}_k\bar{D}_{k-1} = 0$ that 
\begin{align}
	\Sigma_kV_k^TU_k\Sigma_{k-1} = 0.
\end{align}
Therefore, the columns of $V_k$ corresponding to non-zero singular values of $\bar{D}_k$ are orthogonal to columns of $U_k$ associated with non-zero singular values of $\bar{D}_{k-1}$. In addition, it follows from
\begin{align}
	L_k = V_{k}\Sigma_k^2V_k^T + U_{k}\Sigma_{k-1}^2U_k^T
\end{align}
that the spectrum of $\bar{L}_k$ is given by the union of squared non-zero singular values of $\bar{D}_k$ and $\bar{D}_{k-1}$, and $0$, with the multiplicity of $0$ given by the $k$-th Betti numbers. The columns of $U_k$ and $V_k$ corresponding to non-zero singular values, together with the set of harmonic forms, span the entire space of differential $k$-forms. 

In the case that $\dim(M) = 3$, for each type of boundary condition, we have four Laplacians of different degrees in total $k = 0,1,2,3$:
\begin{align}
	\bar{L}_{0} &= \bar{D}_{0}^T\bar{D}_{0}\\
	\bar{L}_{1} &= \bar{D}_{1}^T\bar{D}_{1} + \bar{D}_{0}\bar{D}_{0}^T\\
	\bar{L}_{2} &= \bar{D}_{2}^T\bar{D}_{2} + \bar{D}_{1}\bar{D}_{1}^T\\
	\bar{L}_{3} &= \bar{D}_{2}\bar{D}_{2}^T.
\end{align}
Due to the aforementioned discussion on the spectrum of Laplacians and the duality of the normal and tangential boundary conditions, the spectral analysis of all Laplacians can be reduced to the singular spectra analysis of the three discrete differentials $\bar D_0$, $\bar D_1$, and $\bar D_2$ with one type of boundary conditions. Note that the numerical evaluation of the singular values of these differentials, in the simplicial mesh case, may differ for the two types of boundary conditions, as the DoF for normal $k$-forms and tangent $m-k$ forms are different. However, in the Cartesian representation, they are strictly equivalent to each other by shifting the grid in all directions of the axis by $\ell/2$, so long as $M$ is at least one grid spacing away from the boundary of the grid. 

For the computation of the spectra of the Laplacians, we choose the normal boundary condition. The spectra of all Laplacians $\bar{L}_{k,n}$ for compact domains in $\RR^3$ can be finally decomposed into three distinct parts: the squared singular values of the gradient of tangential scalar fields, denoted by $T$, the squared singular values of the gradient of normal scalar fields, denoted by $N$, and the squared singular values of the curl of tangential curl fields, denoted by $C$.

\section{Persistent de Rham-Hodge Laplacians}

In this section, we present the construction of the persistent de Rham-Hodge Laplacian on differentiable manifolds, which is based on the filtration of manifolds induced by varying a single parameter (the filtration parameter). The spectra of Laplacians carry rich topological and geometric information of a manifold. However, a single manifold might not provide enough information in applications like feature extraction for machine learning analysis. As such, instead of studying just a single manifold, one could examine the spectra of a family of manifolds by adjusting the filtration parameter. The spectra of the Laplacians from this family of manifolds could provide much more information than by considering just one, as the topology and geometry could change for different parameters. This single-parameter family of manifolds, called the evolution of manifolds, was first introduced in \cite{chen2019evolutionary} based on tetrahedral meshes. We briefly recap the background. 

The formal definition of the evolving manifold is given by a one-parameter family of immersions $F_c = F(\cdot, c)$ with $F: B\times [a,b]\to N$ being a smooth map, where $B$ is called the base manifold, $N$ is the ambient manifold, and $c\in[a,b]$ is a real parameter within the interval. In practice, the most common way to define the evolution of manifolds without specifying $B$ is through a level set function by adjusting the isovalues. Given a function $f: N\to [a,b]$, then in our case, we consider the sublevel sets $M = \{x\in N\,|\, f(x)\leq c\}$ with boundary given as $\partial M = \{x\in N\,|\, f(x) = c\}$ for $c\in [a,b]$. A sequence of manifolds can then be obtained by considering evenly distributed isovalues of the function $f$ with the inclusion map
\begin{align}
	M_0\xhookrightarrow{} M_1\xhookrightarrow{} M_2\xhookrightarrow{} \cdots \xhookrightarrow{} M_{s-1} \xhookrightarrow{} M_s,
\end{align}
where each $M_l$ is given as the sublevel set corresponding to $c_l$ with $a\leq c_0<c_1<\cdots < c_s \leq b$. To ensure that $M_l$ is a manifold, we assume that the function $f$ is a Morse function on $N$, and none of the $c_l$'s corresponds to a critical value of the function $f$, i.e., $f^{-1}(c_l)$ does not contain any critical points. This is always possible as the set of Morse functions on a compact manifold is dense in the space of smooth functions, and their critical points are isolated, nondegenerate, and finite for compact manifolds. Thus we can always perturb any input function slightly to avoid critical isovalues in $\{c_l, l =0, 1, \cdots s\}$. In addition, we assume that for each $l$, $M_{l, l+1} = \overline{M_{l+1}\backslash M_{l}}=\{x\in N|\, f(x)\in[c_l,c_{l+1}]\}$ contains at most one critical point, which can be realized by refining the parameter sequence. Note that both $M_l$ and $M_{l, l+1}$ are compact. By Morse theory, if $M_{l, l+1}$ contains no critical points, $M_l$ is diffeomorphic to $M_{l+1}$. The retraction from $M_{l+1}$ to $M_l$ can be easily constructed by considering a flow along the gradient of the function. As $M_{l+1}$ is homotopic to $M_l$ in this case, there is no topological change happening between $(c_l, c_{l+1})$. For the other case when there is exactly one critical point in $M_{l, l+1}$, the manifold $M_{l+1}$ is homotopic to $M_l$ with a $k$-cell attached, where $k$ is the index of the critical point, defined to be the dimension of the largest subspace on which the Hessian $\operatorname{Hess}(f)(x)$ is negative definite. The topological change of the sublevel sets occurs precisely at the critical values of the level set function. Depending on the type of the critical points, i.e., local minimum, saddle points, and local maximum, the topology changes in different ways. In general, a local maximum has the full index $m$, a local minimum has index $0$, while saddle points have indices strictly between 0 and $m$. In the case of $\RR^3$, the occurrences of minima and maxima correspond to the birth of the $0$-th generators and the death of the $2$nd homology generators respectively, while the occurrences of $1$-saddle points correspond to the birth of $1$st homology generators or the death of the $0$-th homology generators, and those of $2$-saddle points correspond to the birth of $2$nd homology generators or the death of $1$st homology generators.

\subsection{Persistent harmonic forms}
As the de Rham complex depends on the topology, it can also be extended to the filtration of manifolds. Due to the duality of the normal and tangential boundary conditions, without loss of generality, one may focus on the space of normal differential forms. Given $M_l\xhookrightarrow{} M_{l+1}$, we then need to construct a map from the space of normal $k$-forms $\Omega^k_n(M_l)$ to the space of normal $k$-forms $\Omega^k_n(M_{l+1})$, that extends each normal $k$-form on $M_l$ to a normal $k$-form on $M_{l+1}$. Let $\omega\in\Omega^k_n(M_l)$. The idea is to utilize the boundary condition of $\omega$ on $M_l$ and extend the forms $\omega|_{\partial M_l}$ to \emph{exact} normal forms on the domain $M_{l, l+1}$ with certain boundary conditions on $\partial M_{l, l+1} = \partial M_l\cup \partial M_{l+1}$. Then the combination $\overline{\omega}$ defines a normal $k$-form on the manifold $M_{l+1}$. Note however that $\delta \overline{\omega}$ is only $0$ in $M_l,$ so the extension of $\omega \in \ker \delta$ may no longer be in $\ker \delta$ on $M_{l+1}.$

To be specific, we consider the biharmonic equation $\Delta^2\zeta=\Delta(\Delta\zeta) = 0 $ on $M_{l, l+1}$ with both Dirichlet and Neumann boundary conditions to ensure the smoothness of $d \zeta$ with $\omega$ through $\partial M_{l}$. Note that $d \zeta$ satisfies the normal boundary condition on $M_{l,l+1}$. Let $\overline{\omega}$ be the extension of $\omega$ on $M_{l+1}$ with $\overline{\omega} = \omega $ on $M_{l}$ and $\overline{\omega} = d \zeta$ on $M_{l, l+1}$. It follows that $\overline{\omega}\in\Omega^k_n(M_{l+1})$ as it satisfies the normal boundary condition $\overline{\omega}|_{\partial M_{l+1}} = \zeta|_{\partial M_{l+1}}=0$. 

While the biharmonic equation produces a smooth extension, in practice, it is more efficient to consider the harmonic extension with the boundary condition $\Delta\zeta = 0$ with the boundary condition of $\star d\zeta = \star \omega$ on $\partial M_l$ and the typical normal form boundary condition on $\partial M_{l+1}$. The solution, by~\cite[Theorem 3.4.10]{schwarz2006hodge}, is unique. The resulting $\bar\omega$ is continuous but nonsmooth as $\delta \bar\omega$ may lead to a Dirac distribution on $\partial M_l$ when $M_{l,l+1}$ induces a topological change. For instance, for a harmonic normal 1-form $\omega$ on a spherical shell $M_l$ with $M_{l+1}$ turning into a solid ball, the biharmonic extension would create a uniform divergence $\delta\bar\omega$ in $M_{l,l+1},$ whereas the harmonic extension creates a thin layer of nonzero divergence $\delta\bar\omega$ near the part of $\partial M_l$ around the cavity in the middle. Thus, the harmonic extension serves the same purpose in reducing the kernel of $\delta$.

Denote by $\mathcal{I}_{l, 1}$ the map from $\Omega^k_n(M_l)$ to $\Omega^k_n(M_{l+1})$ sending $\omega$ to $\overline{\omega}$.  Note that $(d\circ\mathcal{I}_{l, 1})(\omega)$ is $0$ on $M_{l,l+1}$ and thus the same as the extension of the differential of a normal form $d\omega$ on $M_{l, l+1}$, i.e., $d\circ\mathcal{I}_{l, 1} =\mathcal{I}_{l, 1}\circ d$. It follows that there is a commutative diagram
\[ 
\begin{tikzcd}
	\Omega^{0}_n(M_0) \arrow{r}{d} \arrow{d}{\mathcal{I}^0_{0,1}} & \Omega^1_n(M_0) \arrow{r}{d} \arrow{d}{\mathcal{I}^1_{0,1}} & \Omega^{2}_n(M_0) \arrow{r}{d} \arrow{d}{\mathcal{I}^2_{0,1}}  & \Omega^{3}_n(M_0) \arrow{d}{\mathcal{I}^3_{0,1}} \\
	\Omega^{0}_n(M_1) \arrow{r}{d} \arrow{d}{\mathcal{I}^0_{1,1}} & \Omega^1_n(M_1)  \arrow{r}{d} \arrow{d}{\mathcal{I}^1_{1,1}} & \Omega^{2}_n(M_1)  \arrow{r}{d} \arrow{d}{\mathcal{I}^2_{1,1}}  & \Omega^{3}_n(M_1)  \arrow{d}{\mathcal{I}^3_{1,1}} \\
	\Omega^{0}_n(M_2)  \arrow{r}{d} \arrow{d}{\mathcal{I}^0_{2,1}} & \Omega^1_n(M_2) \arrow{r}{d} \arrow{d}{\mathcal{I}^1_{2,1}} & \Omega^{2}_n(M_2) \arrow{r}{d} \arrow{d}{\mathcal{I}^2_{2,1}}  & \Omega^{3}_n(M_2) \arrow{d}{\mathcal{I}^3_{2,1}} \\
	\cdots  & \cdots  & \cdots  & \cdots
\end{tikzcd}
\]
where the horizontal direction gives the de Rham complex and the vertical direction shows the filtration-induced extensions.

Next, we introduce the $p$-persistent Hodge Laplacian. Let $\mathcal{I}_{l, p} = \mathcal{I}_{l+p-1, 1}\circ ...\circ \mathcal{I}_{l, 1}$, which then gives an extension map from the space of normal forms on $M_l$ to the space of normal forms on $M_{l+p}$. We have the following commutative diagram
\[\begin{tikzcd}
	{ } && {\Omega_n^{k}(M_{l})} && {\Omega_n^{k+1}(M_{l})} \\
	\\
	{\Omega_n^{k-1}(M_{l+p})} && {\Omega_n^{k}(M_{l+p})} && { }
	\arrow["{d_{l+p}^{k-1}}", shift left, from=3-1, to=3-3]
	\arrow["{R_{l,p}}", shift left, from=3-3, to=1-3]
	\arrow["{\delta_{l+p}^k}", shift left, from=3-3, to=3-1]
	\arrow["{d_{l}^k}", shift left, from=1-3, to=1-5]
	\arrow["{\tilde{\delta}_{l,p}^k}", shift left=-2, from=1-3, to=3-1]
	\arrow["{I_{l,p}}", shift left, from=1-3, to=3-3]
	\arrow["{\delta_{l}^{k+1}}", shift left, from=1-5, to=1-3]
	\arrow["{\tilde{d}_{l,p}^{k-1}}", shift left=4, from=3-1, to=1-3]
\end{tikzcd}.\]

Here $d_l, \delta_l$ denotes the differential and codifferential on $\Omega^k(M_l)$, $d_{l+p}, \delta_{l+p}$ denotes the differential and codifferential on $\Omega^k(M_{l+p})$, respectively, and $\mathcal{R}_{l, p}$ is the projection of differential forms in $\Omega^k_n(M_{l+p})$ to the space spanned by the harmonic extensions followed by the restriction to $M_{l}$. Let $\tilde{\delta}_{l,p} = \delta_{l+p}\circ\mathcal{I}_{l, p}$ and $\tilde{d}_{l,p} = \mathcal{R}_{l, p}\circ d_{l+p}$. By the construction of the extension, we have 
$(\tilde{\delta}_{l,p}\omega, \eta) = (\omega, \tilde{d}_{l,p}\eta)$, i.e., $\tilde{\delta}_{l,p}$ are adjoint to $\tilde{d}_{l,p}$. 
We then define the $p$-persistent Hodge Laplacian operator $\Delta^p_{n, l}: \Omega^k_n(M_l)\to\Omega^k_n(M_{l})$ as follows
\begin{align}
	\Delta^p_{n, l} = \tilde{d}_{l,p}\tilde{\delta}_{l,p} + 
	\delta_{l}d_{l}.
\end{align}
It is easy to see that when $p=0$, the $p$-persistent Hodge Laplacian gives exactly the usual Hodge Laplacian $\Delta_{n, l}: \Omega^k_n(M_l)\to\Omega^k_n(M_{l})$ restricted to the space of normal forms. We then define the $p$-persistent normal harmonic fields as the kernel of the $p$-persistent Hodge Laplacian $\mathcal{H}_n^{k,p} = \ker\Delta^p_{n, l}$, which can be identified with the space $\ker\tilde{\delta}_{l,p}\cap\ker d_{l}$. Note that by the extension construction and $\mathcal{R}_{l, p}\circ\mathcal{I}_{l, p} = \operatorname{Id}$, one can see that $\ker\tilde{\delta}_{l,p}\subset\ker\delta$ gets smaller as $p$ increases, which confirms that fewer cohomology generators persist longer. 

\subsection{Discretization of $p$-persistent de Rham cohomology}

The regular Cartesian grid allows one to define persistent graph Laplacian on manifolds in the same way as persistent graph Laplacian~\cite{wang2020persistent}. It also allows defining persistent Hodge Laplacian in a consistent way, with the inclusion of nontrivial Hodge stars.

Recall that the discrete differential $k$-forms can be seen as a $k$-cochain, i.e., a linear mapping from the chain space $\mathcal{C}_k$ to $\RR$ that sends a $k$-chain $c_k = \sum_i a_i\sigma_i$ to $\int_{c_k}\omega = \sum_i a_iW_i$, where $W_i = \int_{\sigma_i}\omega$ is the integral of a smooth $k$-form $\omega$ over the $k$-cell $\sigma_i$.

By varying the isovalue of the level set function $f$, we can get a sequence of cell complexes given as nested sequences of sub-cell complexes of $K$ satisfying the normal boundary conditions.

\begin{figure}[t]
	\centering
	\includegraphics[height=3.2cm]{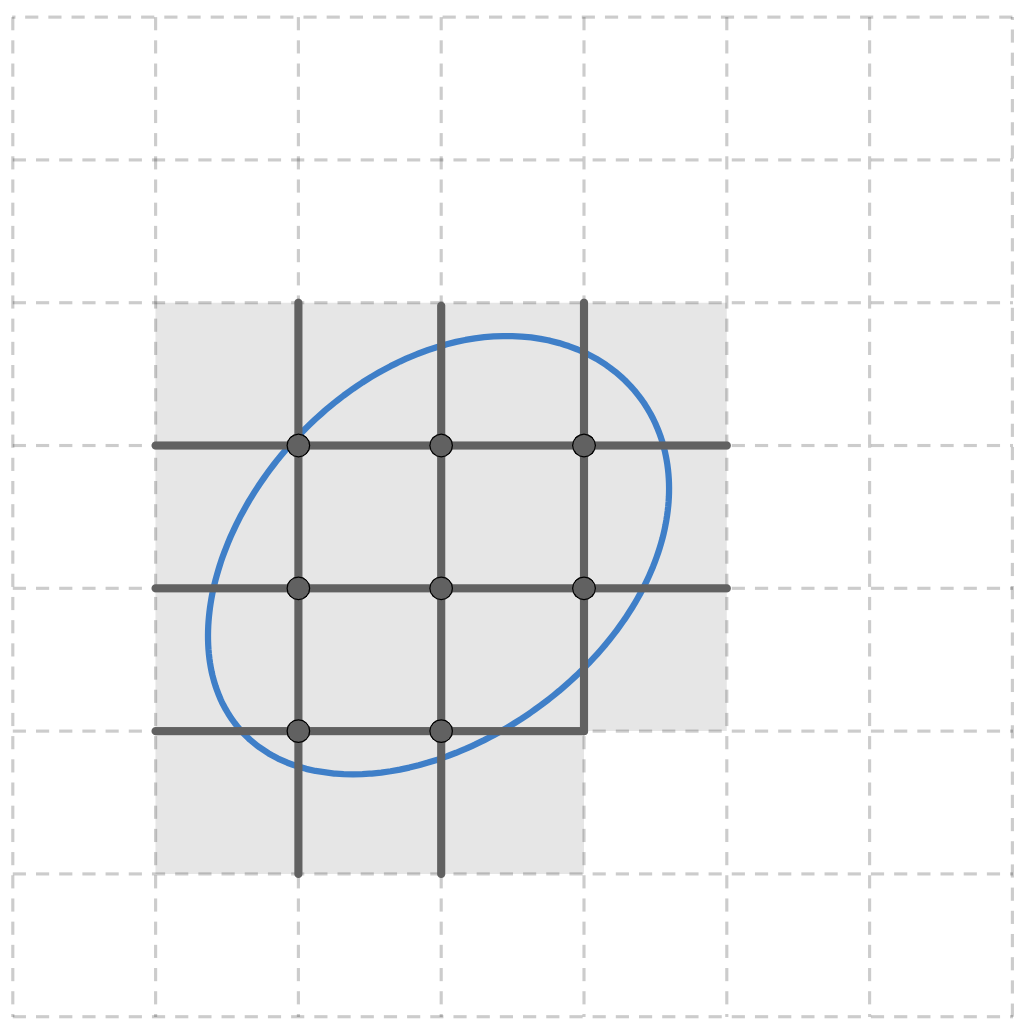}\hskip.2in
	\includegraphics[height=3.2cm]{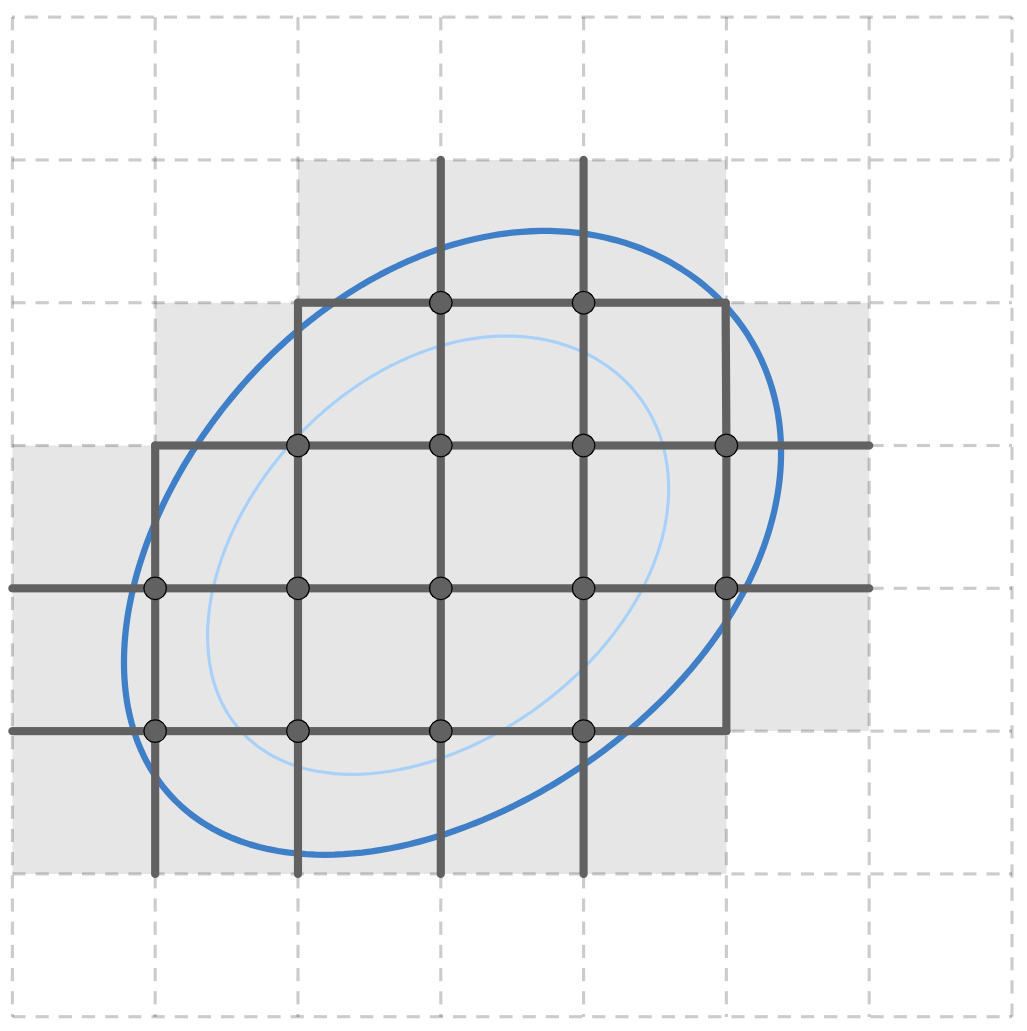}\hskip.2in
	\includegraphics[height=3.2cm]{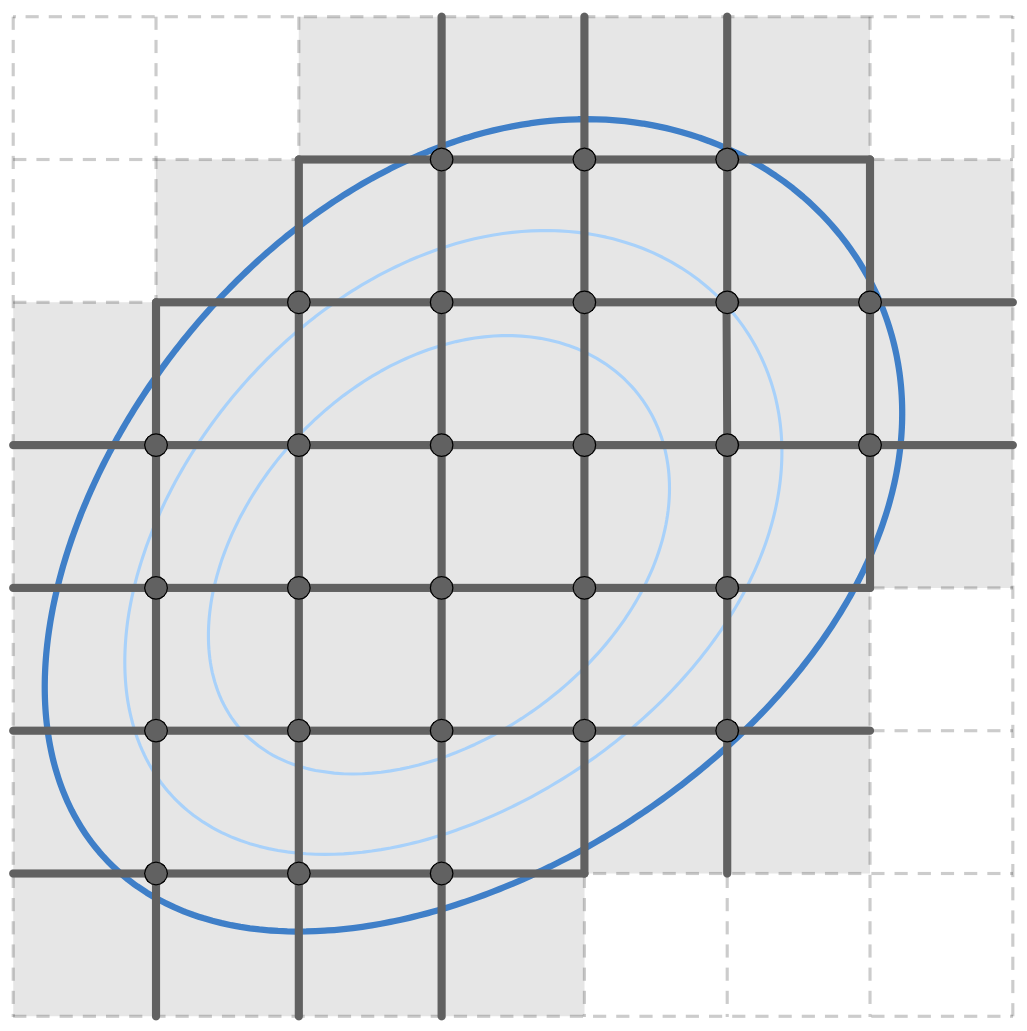}
	\caption{An example of a nested sequence of sub-cell complexes in a 2D Cartesian grid under the normal boundary condition, illustrating the inclusion of normal supports for $0$, $1$, and $2$ discrete differential forms for an evolution of manifolds. Here the manifolds are represented by the bounded regions of the blue isocurves of a level set function.}
	\label{fig.complex_lvf}
\end{figure}

\begin{align}
	\emptyset = K_0\subset K_1 \subset\cdots \subset K_{s-1} \subset K_s = K. 
\end{align}
See Fig.~\ref{fig.complex_lvf} for an example of such a nested sequence of sub-cell complexes in a 2D Cartesian grid. Denote by $\mathcal{C}^k(K_l)$ the space of discrete $k$-forms on subcomplex $K_l$ with $0\leq l\leq s$. Note that $K_l\subset K_{l+1}$. A discrete $k$-form on $K_l$ can be easily extended to $K_{l+1}$ by solving the discrete Laplace equation with the above boundary conditions for values on every $k$-cells in $K_{l,l+1}=\mathrm{Cl} (K_{l+1}\backslash K_l)$, the closure of the difference complex. We denote this extension map as $I_{l, 1}: \mathcal{C}^k(K_l)\to \mathcal{C}^k(K_{l+1})$ and by $I_{l, p} = I_{l+p-1, 1}\circ I_{l+p-2, 1}\circ \cdots \circ I_{l, 1}: \mathcal{C}^k(K_{l})\to\mathcal{C}^k(K_{l+p})$ the extension mapping from the space of discrete $k$-forms on $K_l$ to the space of discrete $k$-forms on $K_{l+p},$ which may also be constructed directly by solving the Laplace equation on $K_{l,l+p}=\mathrm{Cl}(K_{l+p}\backslash K_l)$. With this extension mapping, the space of discrete $k$-forms on $K_{l}$ can be seen as a subspace of discrete $k$-forms on $K_{l+p}$.

A sequence of the discrete de Rham cochain complexes can be defined as follows:
\[\begin{tikzcd}
	{\mathcal{C}^0(K_0)} & {\mathcal{C}^1(K_0)} & \cdots & {\mathcal{C}^{k}(K_0)} & {\mathcal{C}^{k+1}(K_0)} & \cdots \\
	{\mathcal{C}^0(K_1)} & {\mathcal{C}^1(K_1)} & \cdots & {\mathcal{C}^{k}(K_1)} & {\mathcal{C}^{k+1}(K_1)} & \cdots \\
	{\mathcal{C}^1(K_2)} & {\mathcal{C}^1(K_2)} & \cdots & {\mathcal{C}^{k}(K_2)} & {\mathcal{C}^{k+1}(K_2)} & \cdots \\
	\cdots & \cdots && \cdots & \cdots
	\arrow["{D_0^0}", shift left, from=1-1, to=1-2]
	\arrow["{I_{0,1}}", shift left, from=1-1, to=2-1]
	\arrow["{\delta^1_0}", shift left, from=1-2, to=1-1]
	\arrow["{D_0^1}", shift left, from=1-2, to=1-3]
	\arrow["{I_{0,1}}", shift left, from=1-2, to=2-2]
	\arrow["{\delta_0^2}", shift left, from=1-3, to=1-2]
	\arrow["{D_0^{k-1}}", shift left, from=1-3, to=1-4]
	\arrow["{\delta_0^{k}}", shift left, from=1-4, to=1-3]
	\arrow["{D_0^k}", shift left, from=1-4, to=1-5]
	\arrow["{I_{0,1}}", shift left, from=1-4, to=2-4]
	\arrow["{\delta_0^{k+1}}", shift left, from=1-5, to=1-4]
	\arrow["{D_0^{k+1}}", shift left, from=1-5, to=1-6]
	\arrow["{I_{0,1}}", shift left, from=1-5, to=2-5]
	\arrow["{\delta_0^{k+2}}", shift left, from=1-6, to=1-5]
	\arrow["{D_1^{0}}", shift left, from=2-1, to=2-2]
	\arrow["{I_{1,1}}", shift left, from=2-1, to=3-1]
	\arrow["{\delta_1^{1}}", shift left, from=2-2, to=2-1]
	\arrow["{D_1^{1}}", shift left, from=2-2, to=2-3]
	\arrow["{I_{1,1}}", shift left, from=2-2, to=3-2]
	\arrow["{\delta^2_{1}}", shift left, from=2-3, to=2-2]
	\arrow["{D_1^{k-1}}", shift left, from=2-3, to=2-4]
	\arrow["{\delta_1^{k}}", shift left, from=2-4, to=2-3]
	\arrow["{D_1^{k}}", shift left, from=2-4, to=2-5]
	\arrow["{I_{1,1}}", shift left, from=2-4, to=3-4]
	\arrow["{\delta_1^{k+1}}", shift left, from=2-5, to=2-4]
	\arrow["{D_1^{k+1}}", shift left, from=2-5, to=2-6]
	\arrow["{I_{1,1}}", shift left, from=2-5, to=3-5]
	\arrow["{\delta_1^{k+2}}", shift left, from=2-6, to=2-5]
	\arrow["{D_2^{0}}", shift left, from=3-1, to=3-2]
	\arrow["{I_{2,1}}", shift left,from=3-1, to=4-1]
	\arrow["{\delta_2^{1}}", shift left, from=3-2, to=3-1]
	\arrow["{D_2^{1}}", shift left, from=3-2, to=3-3]
	\arrow["{I_{2,1}}", shift left, from=3-2, to=4-2]
	\arrow["{\delta_2^{2}}", shift left, from=3-3, to=3-2]
	\arrow["{D_2^{k-1}}", shift left, from=3-3, to=3-4]
	\arrow["{\delta_2^{k}}", shift left, from=3-4, to=3-3]
	\arrow["{D_2^{k}}", shift left, from=3-4, to=3-5]
	\arrow["{I_{2,1}}", shift left, from=3-4, to=4-4]
	\arrow["{\delta_2^{k+1}}", shift left, from=3-5, to=3-4]
	\arrow["{D_2^{k+1}}", shift left, from=3-5, to=3-6]
	\arrow["{I_{2,1}}", shift left, from=3-5, to=4-5]
	\arrow["{\delta_2^{k+2}}", shift left, from=3-6, to=3-5]
\end{tikzcd}\]
where $D_l^k: \mathcal{C}^{k+1}(K_l)\to \mathcal{C}^{k}(K_l)$ denotes the discrete differential operator, and $\delta_l^k: \mathcal{C}^k(K_l)\to \mathcal{C}^{k-1}(K_l)$ denotes the discrete codifferential operator on $K_l$.

To define the persistent discrete Hodge Laplacian, we construct the discrete counterparts of $\tilde{d}_{l,p}$ and $\tilde{\delta}_{l,p}$ in the previous section. 

Denote by $\delta_{l,p}^{k+1,n}:\mathcal{C}^{k+1}(K_l)\to\mathcal{C}_{l,p}^{k}$  the operator given as $\delta_{l,p}^{k,n} = \delta_{l+p}^k I_{l,p}^{k,n}$, where $\delta_{l+p}^{k,n}$ is the previously defined discrete operator for $K_{l+p}$ and $I_{l,p}^{k,n}$ is the discrete harmonic extension operator defined next. Assuming $K_{l,l+p}$ contains few  $k$-cells, the harmonic extension is then constructed by the linear system $L^{k-1,n}_{K_{l,l+p}} \zeta = 0,$ and shifting all $\star d\zeta$ values in the overlap of supports of $K_l$ and $K_{l,l+p}$ to the right-hand side and replacing them with a rescaling of $\star \omega$ based on the $k$-volume within each support. More specifically, the resulting system is $\tilde{L}^{k-1,n}_{K_{l,l+p}}\tilde \zeta = -S^{k-1,n}\delta^k_{\partial K_l} \omega,$ where $\tilde{L}^{k-1,n}_{K_{l,l+p}}$ is the Laplace operator applied to a form $\tilde\zeta$ defined on $K_{l,l+p}\backslash \partial K_l,$ and $\delta^k_{\partial K_l}$ is the boundary codifferential operator that uses the values of $\omega$ on $\partial K_l$ to evaluate the neighboring $(k\!-\!1)$-cells in $K_{l,l+p}\backslash \partial K_l.$ 

The resulting extension operator 
\[I_{l,p}= 
\begin{pmatrix} 
\operatorname{Id}_{K_l}\\
-D^k_{K_{l,l+p}}(\tilde{L}^{k,n}_{K_{l,l+p}})^{-1} S^{k,n} \delta^k_{\partial K_l}
\end{pmatrix},
\]
where $\operatorname{Id}_{K_l}$ is the identity matrix in $K_l$ up to a rescaling in the boundary, 
provides the combination of $\omega$ in $K_l$ and  $d\tilde\zeta$ in $K_{l,l+p}\backslash \partial K_l,$ when applied to $\omega.$ The matrix corresponding to $I_{l,p}$ is dense for rows corresponding to cells in $K_{l,l+p}$ but diagonal for rows corresponding to cells in $K_l.$  Note that $\delta_{\partial K_l}$ is not necessarily $0$ for coclosed $\omega$, but is $0$ for coexact $\omega.$

The adjoint operator of $\delta_{l,p}^{k+1,n}$ defines $D^k_{l,p}$. In the following, we drop most of the subscripts for clarity. Recall that $(\omega, \tilde{d} \eta)=(\tilde{\delta} \omega,\eta)$ can be discretized as
\[[W]^T S [\tilde{D} E]= [S^{-1} D^T S I_{l,p} W]^T S [E]\]
with $W$ and $E$ as discrete versions of $\omega$ and $\eta$. Thus $\tilde{D}= S^{-1} I_{l,p}^T S D,$ from which we may recognize the restriction operator as $R=S^{-1}I_{l,p}^T S$. This restriction operator can be seen as the $L_2$-projection onto the space formed by all harmonic extensions from $\Omega_n^K(M_l).$ 


Note that in this case, we immediately $\delta_{l,p}^k\delta_l^{k+1} = 0$, since the extension operator will generate $\tilde\zeta=0$ for any coexact form $\omega=\delta\beta$ on $K_l$ as the right-hand side of the associated linear system essentially corresponds to $\delta\delta\beta=0$. From the adjoint version, we have  $D_l^{k}D_{l,p}^{k-1} = 0,$ and thus the following commutative diagram
\[\begin{tikzcd}
	{ } && {\mathcal{C}^{k}(K_{l})} && {\mathcal{C}^{k+1}(K_{l})} \\
	\\
	{\mathcal{C}^{k-1}(K_{l+p})} && {\mathcal{C}^{k}(K_{l+p})} && { }
	\arrow["{D_{l+p}^{k-1}}", shift left, from=3-1, to=3-3]
	\arrow["{R_{l,p}}", shift left, from=3-3, to=1-3]
	\arrow["{\delta_{l+p}^k}", shift left, from=3-3, to=3-1]
	\arrow["{D_{l}^k}", shift left, from=1-3, to=1-5]
	\arrow["{\delta_{l,p}^k}", shift left=-2, from=1-3, to=3-1]
	\arrow["{I_{l,p}}", shift left, from=1-3, to=3-3]
	\arrow["{\delta_{l}^{k+1}}", shift left, from=1-5, to=1-3]
	\arrow["{D_{l,p}^{k-1}}", shift left=4, from=3-1, to=1-3]
\end{tikzcd}.\]

The discrete $p$-persistent Hodge Laplacian is then given as follows 
\begin{align}
	L_{l,p}^k = D_{l,p}^{k-1}\delta_{l,p}^k + \delta_{l}^{k+1}D_l^k,
\end{align}
and the discrete $p$-persistent BIG Laplacian is  
\begin{align}
	L_{l,p}^k = D_{l,p}^{k-1}(D_{l,p}^{k-1})^T + (D_l^k)^TD_l^k.
\end{align}

\begin{figure}[t]
	\centering
	\includegraphics[height=3.4cm]{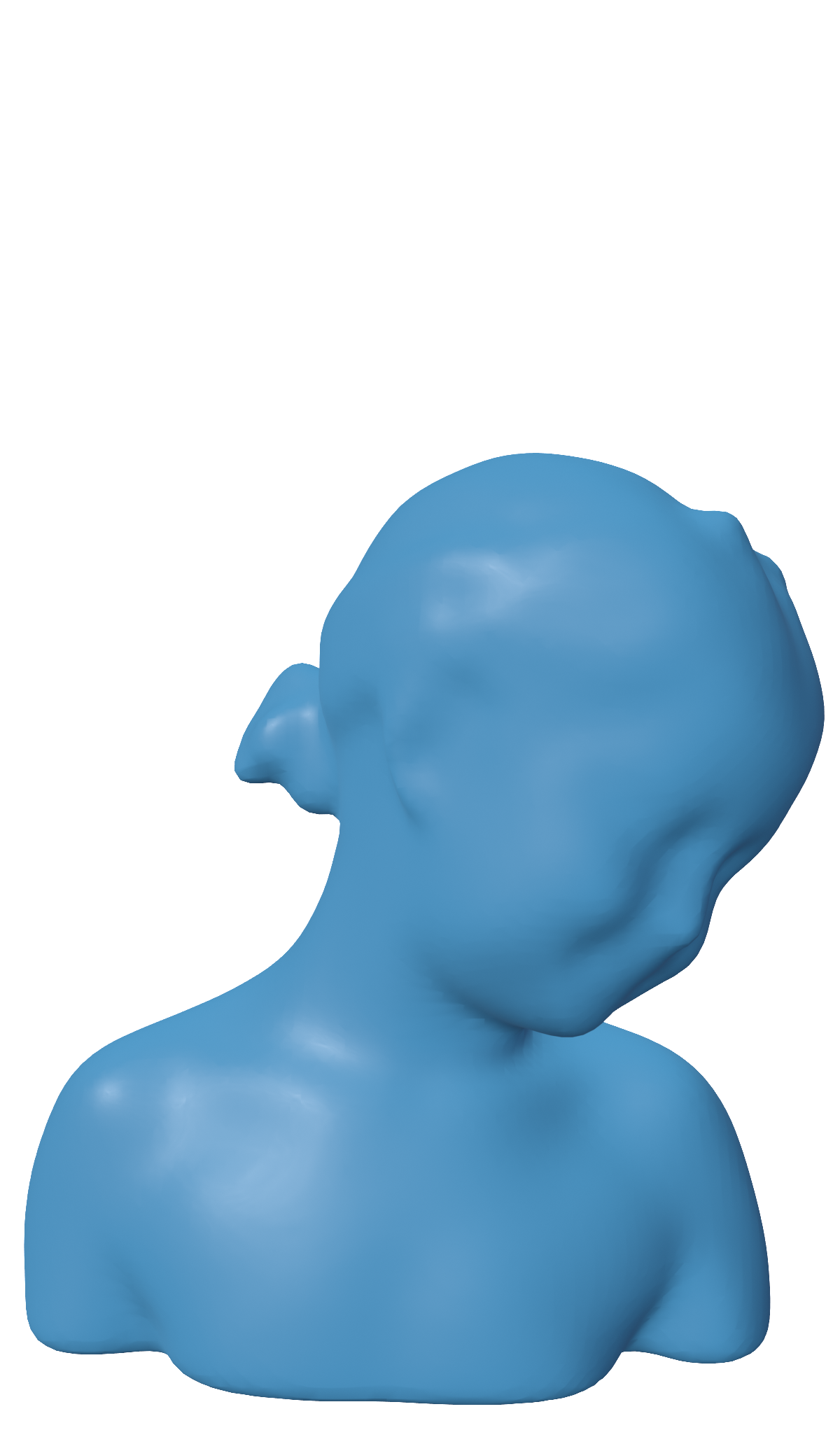}
	\includegraphics[height=3.4cm]{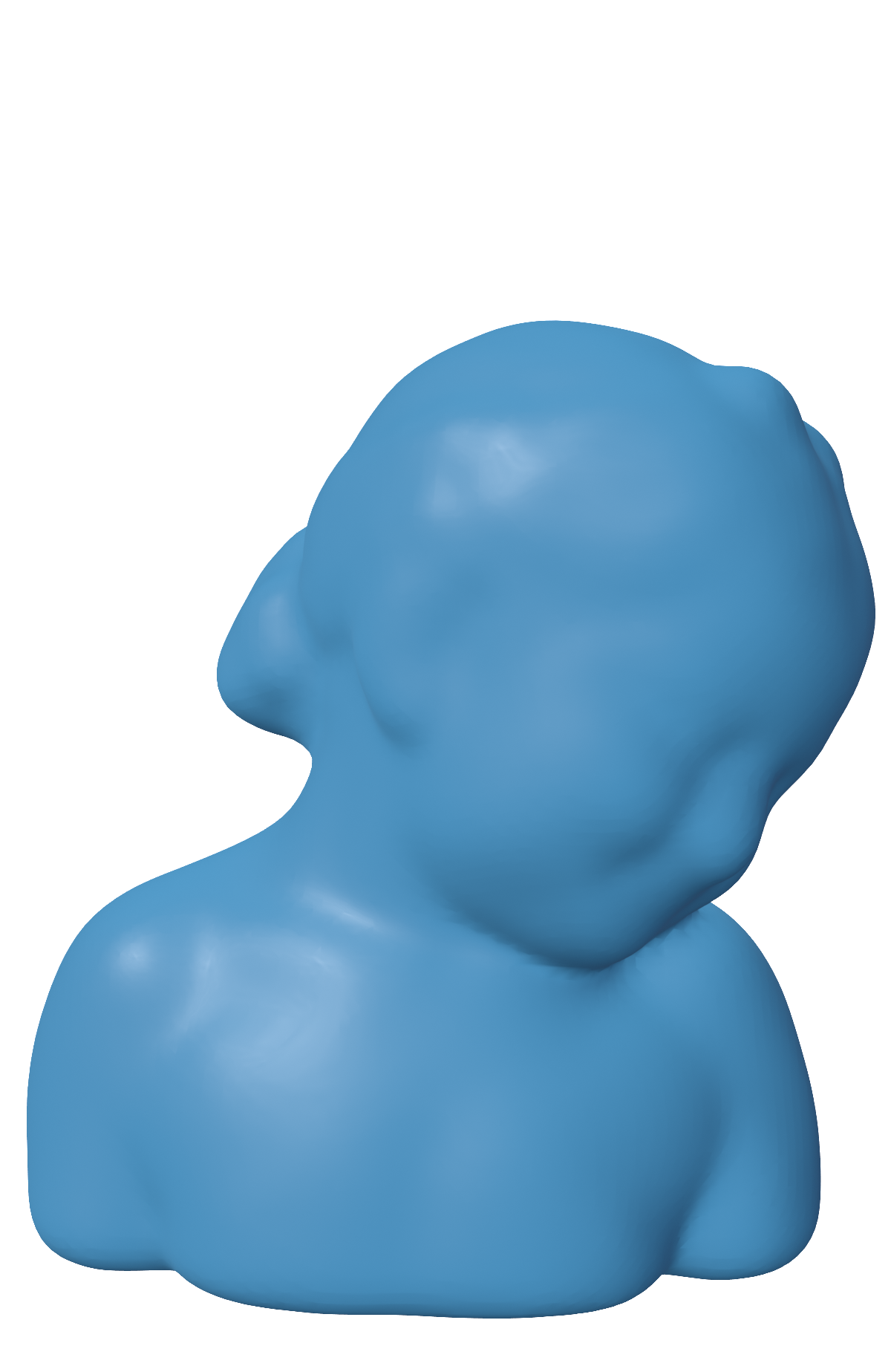}
	\includegraphics[height=3.4cm]{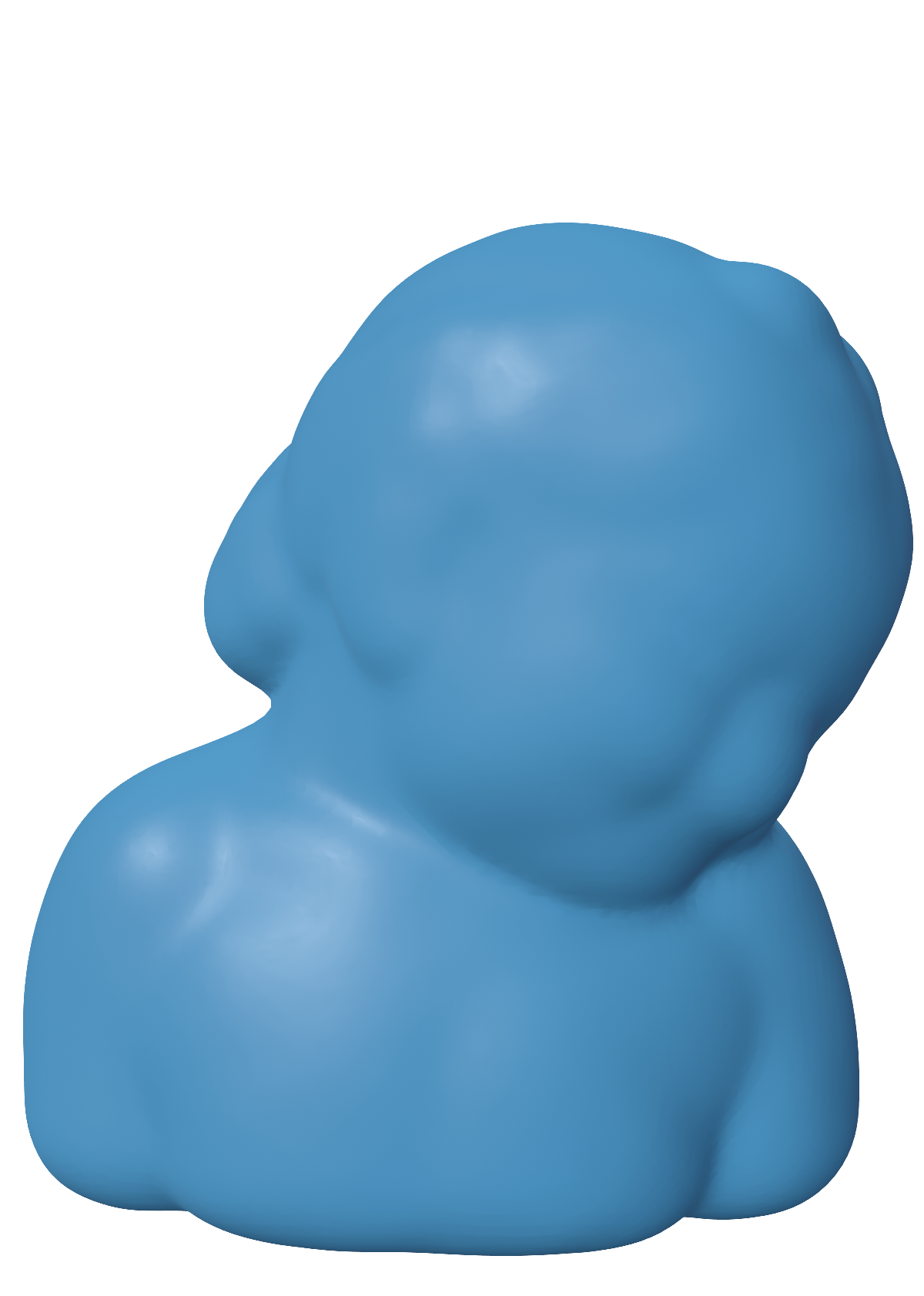}
	\includegraphics[height=3.4cm]{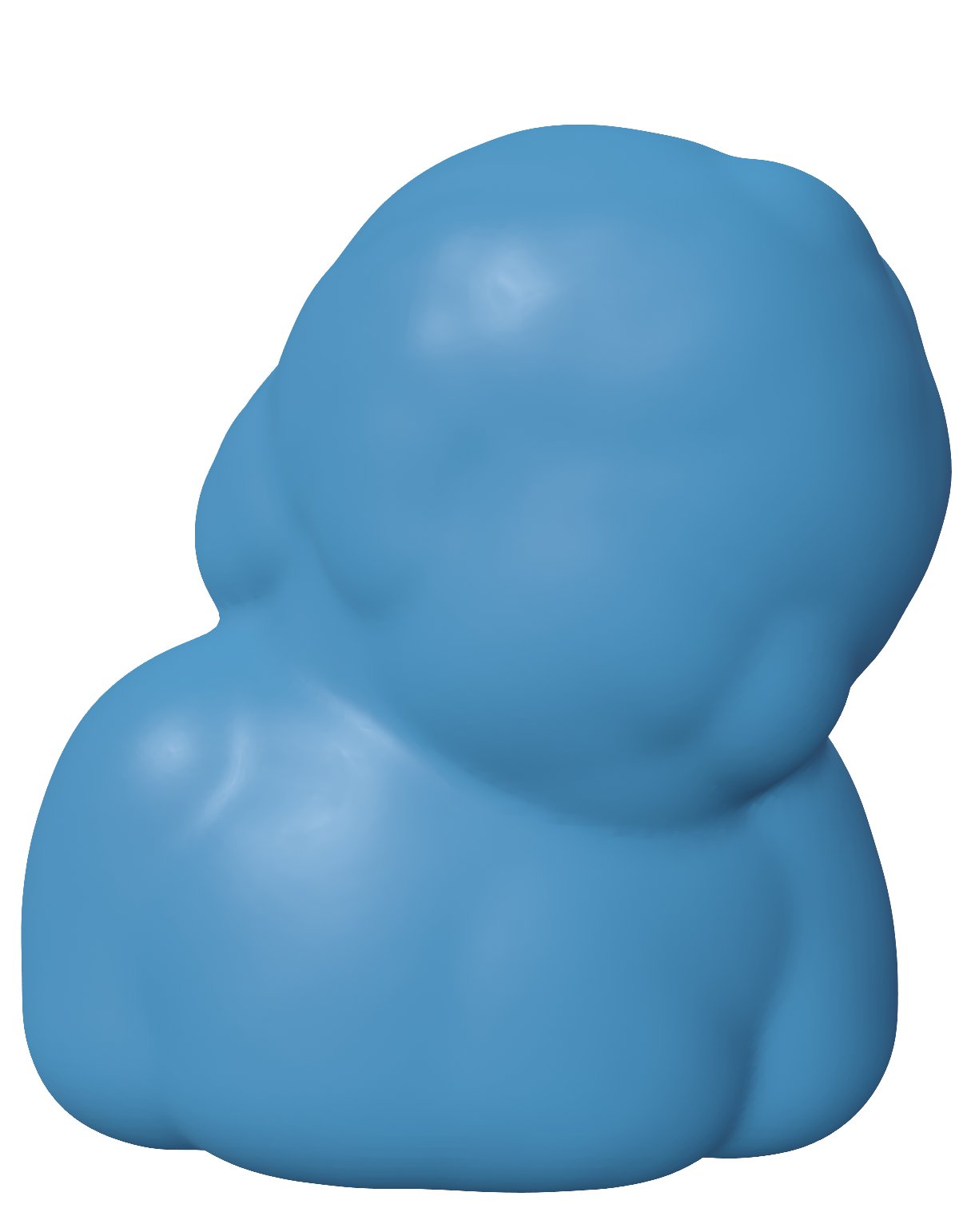}
	\includegraphics[height=3.4cm]{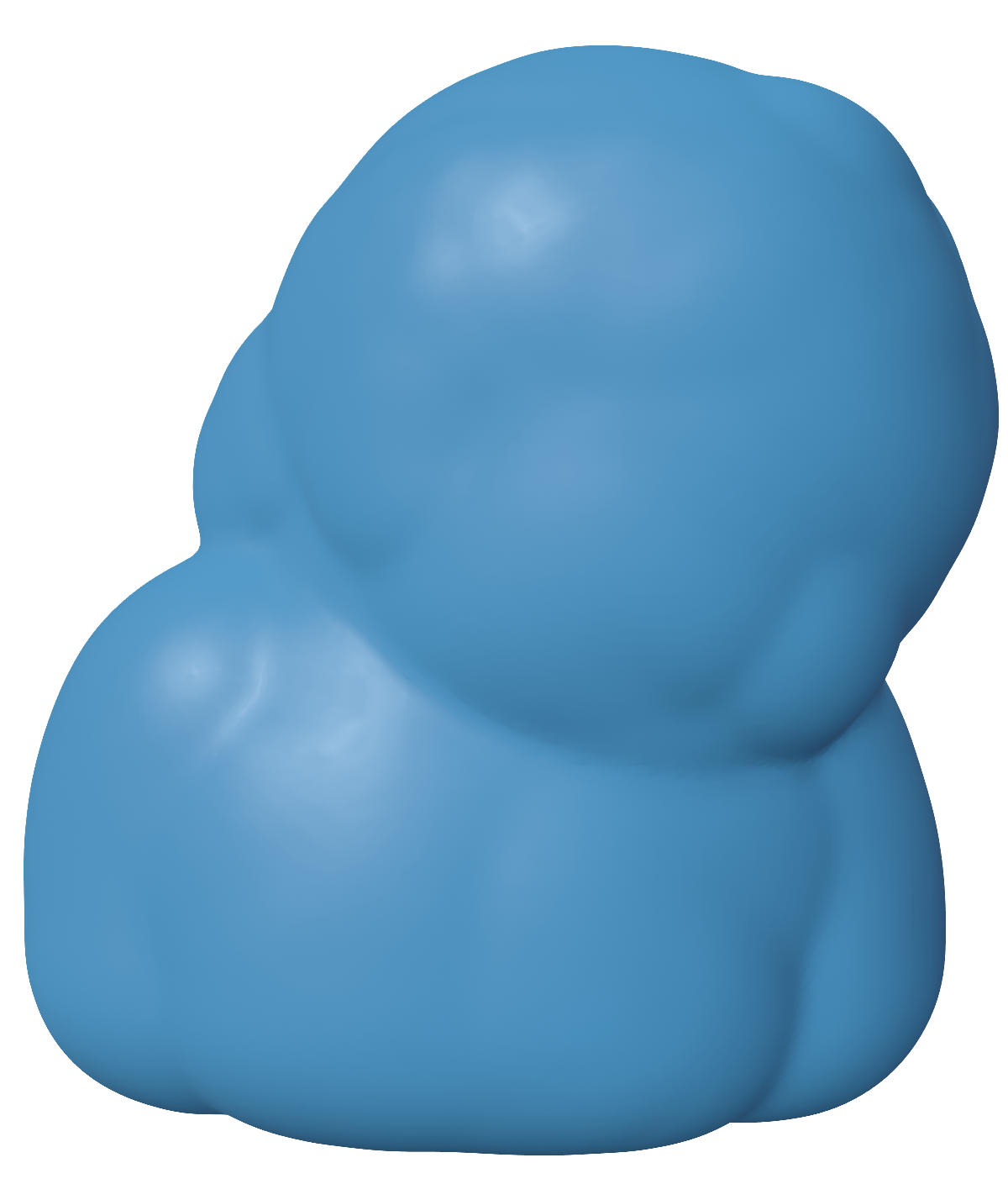}
	
	\includegraphics[height=3.425cm, trim=0.3cm 0cm 0.25cm 0,clip]{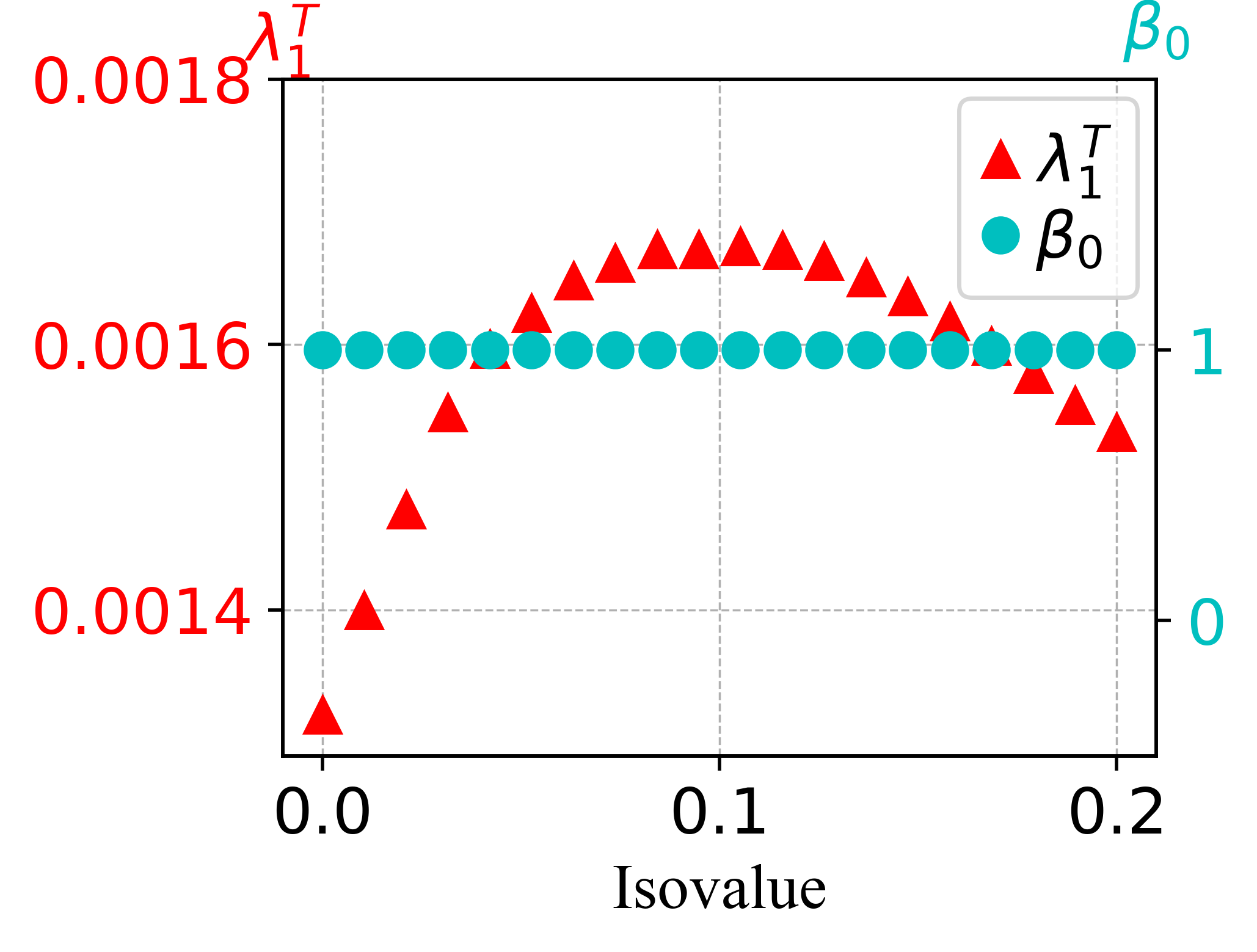}
	\includegraphics[height=3.425cm, trim=0.3cm 0cm 0.25cm 0,clip]{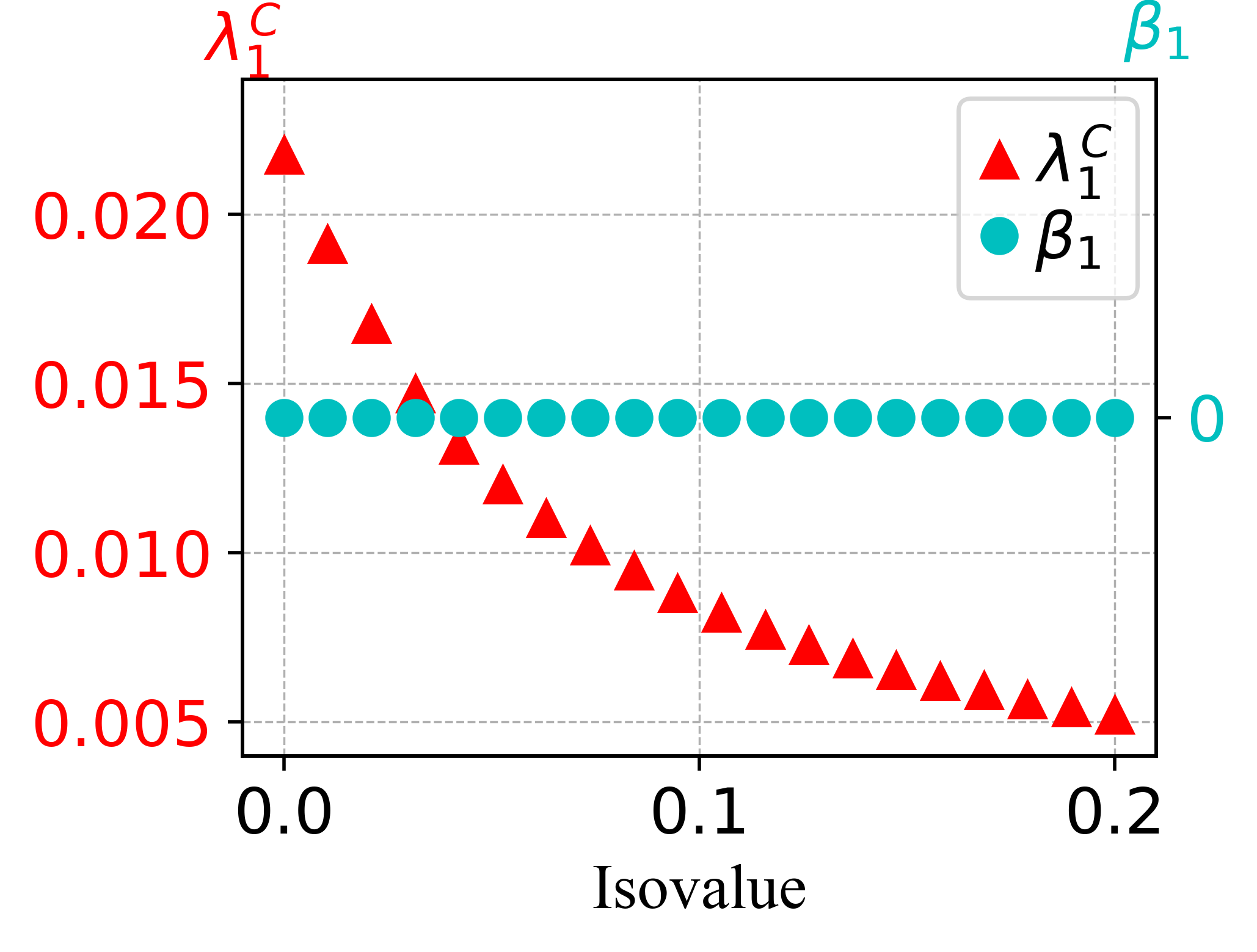}
	\includegraphics[height=3.425cm, trim=0.3cm 0cm 0.25cm 0,clip]{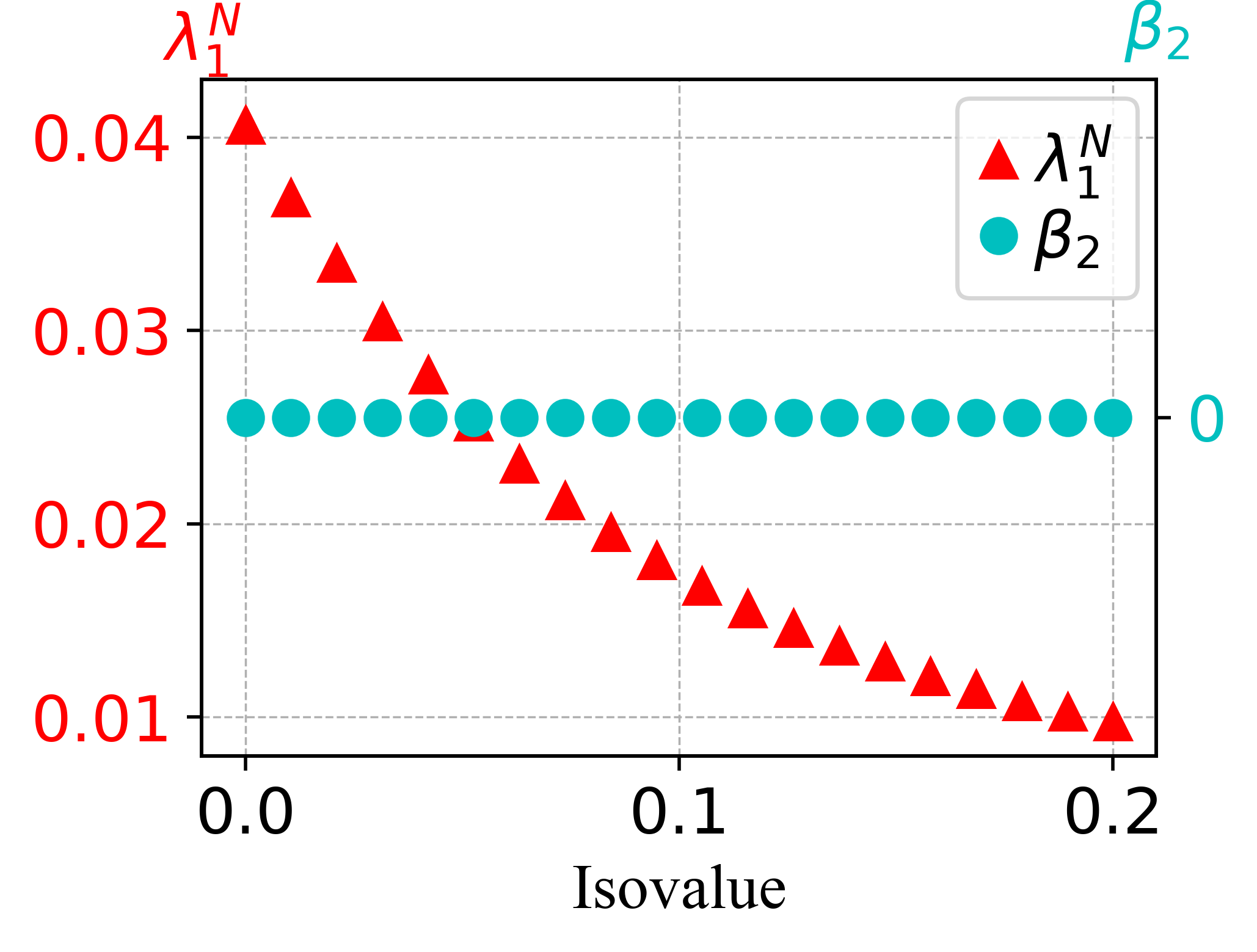}
	\caption{First row: Snapshots of evolving manifolds for the Bimba model. Second row: Changes in Betti numbers $\beta_0$, $\beta_1$, $\beta_2$ and the first non-zero eigenvalues in T, C, N along $20$ evenly spaced isovalues from $0$ to $0.2$. Here the first shape in the top first row corresponds to isovalue $0$ and the last shape in the first row corresponds to isovalue $0.2$. $\lambda^T_1$, $\lambda^C_1$ and $\lambda^N_1$ are the first non-zero eigenvalues in the set T, C, N, respectively. The signed distance function generated from the original Bimba model is used as the level set function.}
	\label{fig.emfld.bimba}
\end{figure}

\begin{figure}[h]
	\centering
	\includegraphics[height=4.1cm]{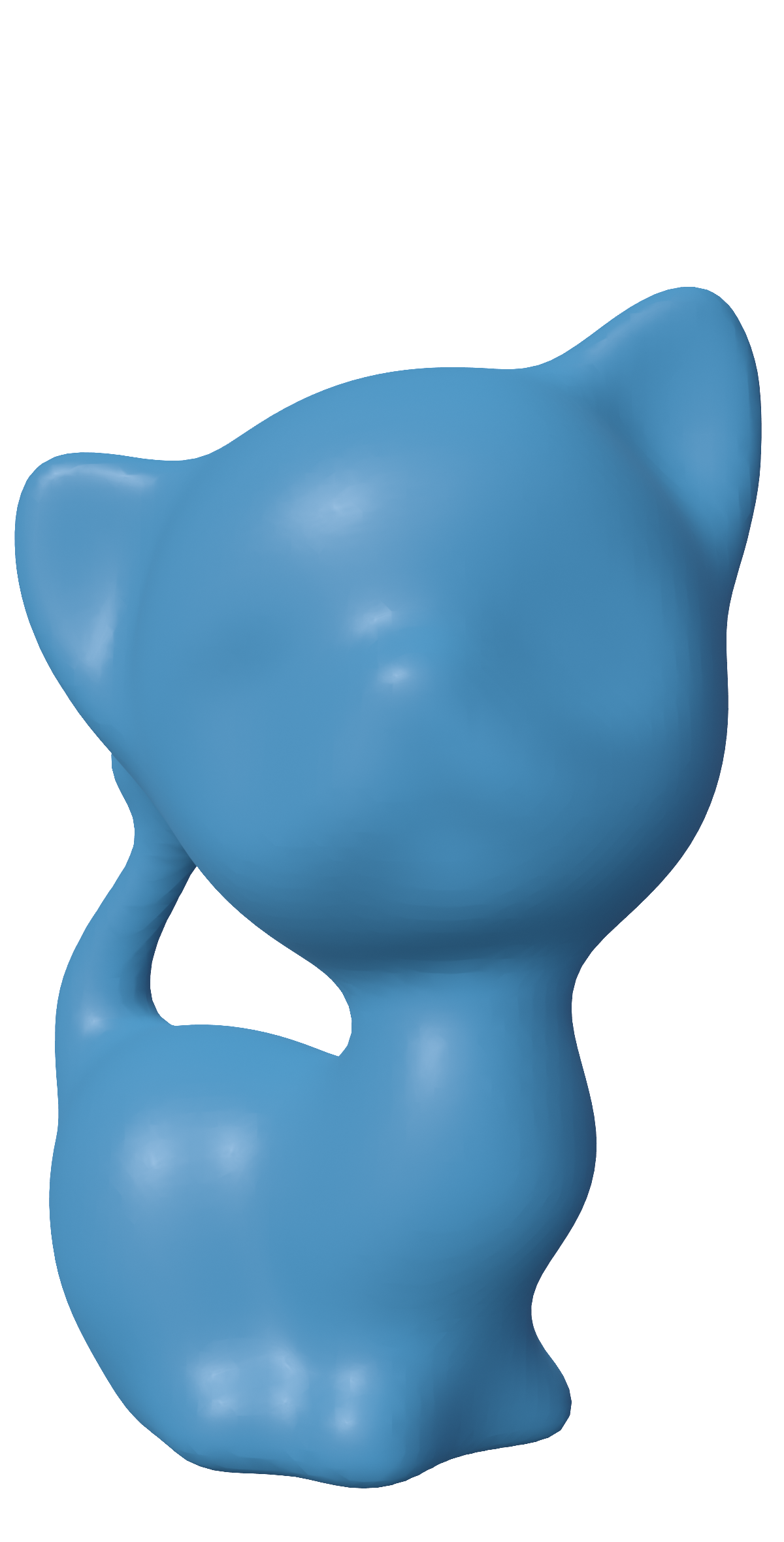}
	\includegraphics[height=4.1cm]{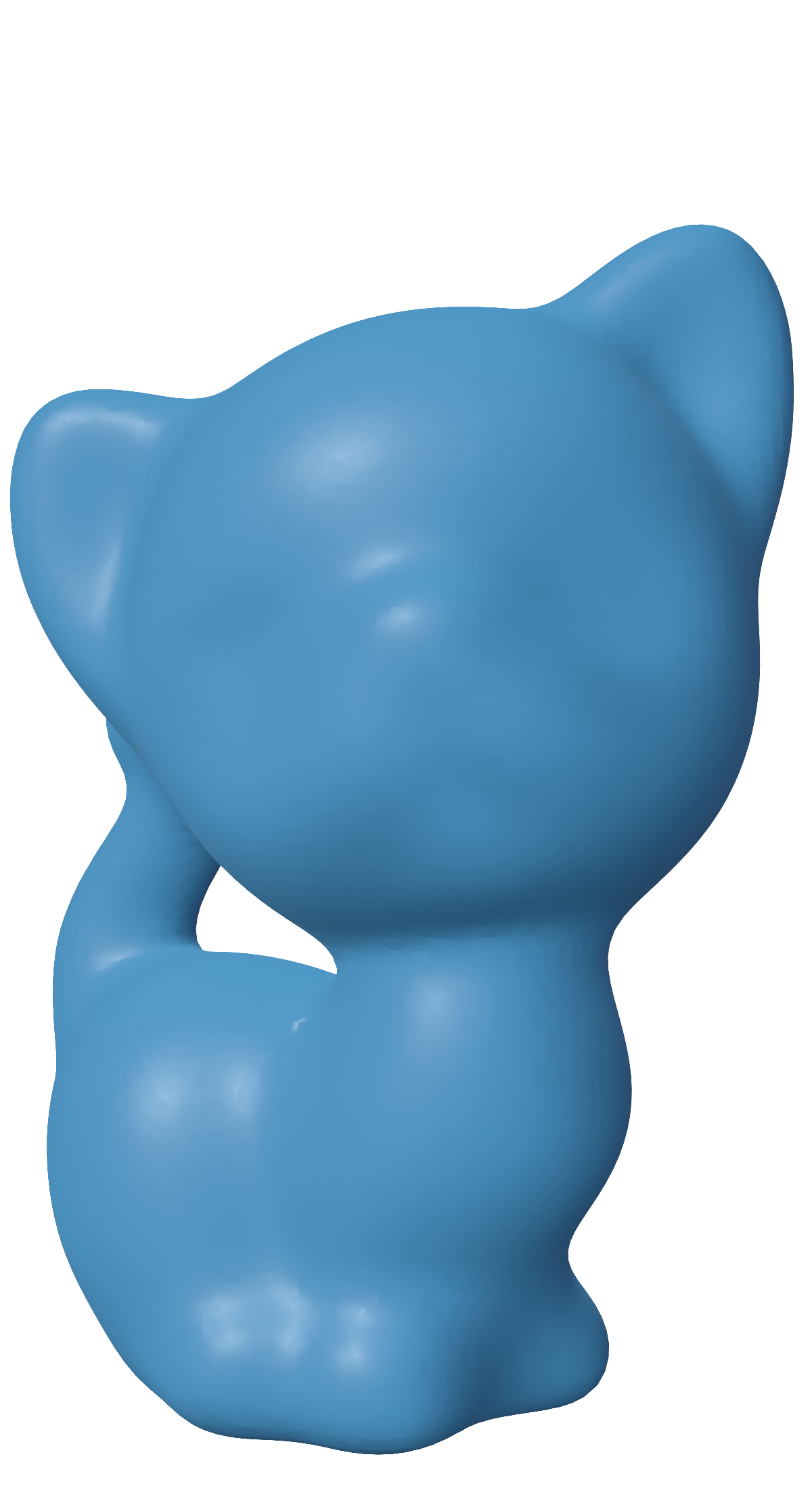}
	\includegraphics[height=4.1cm]{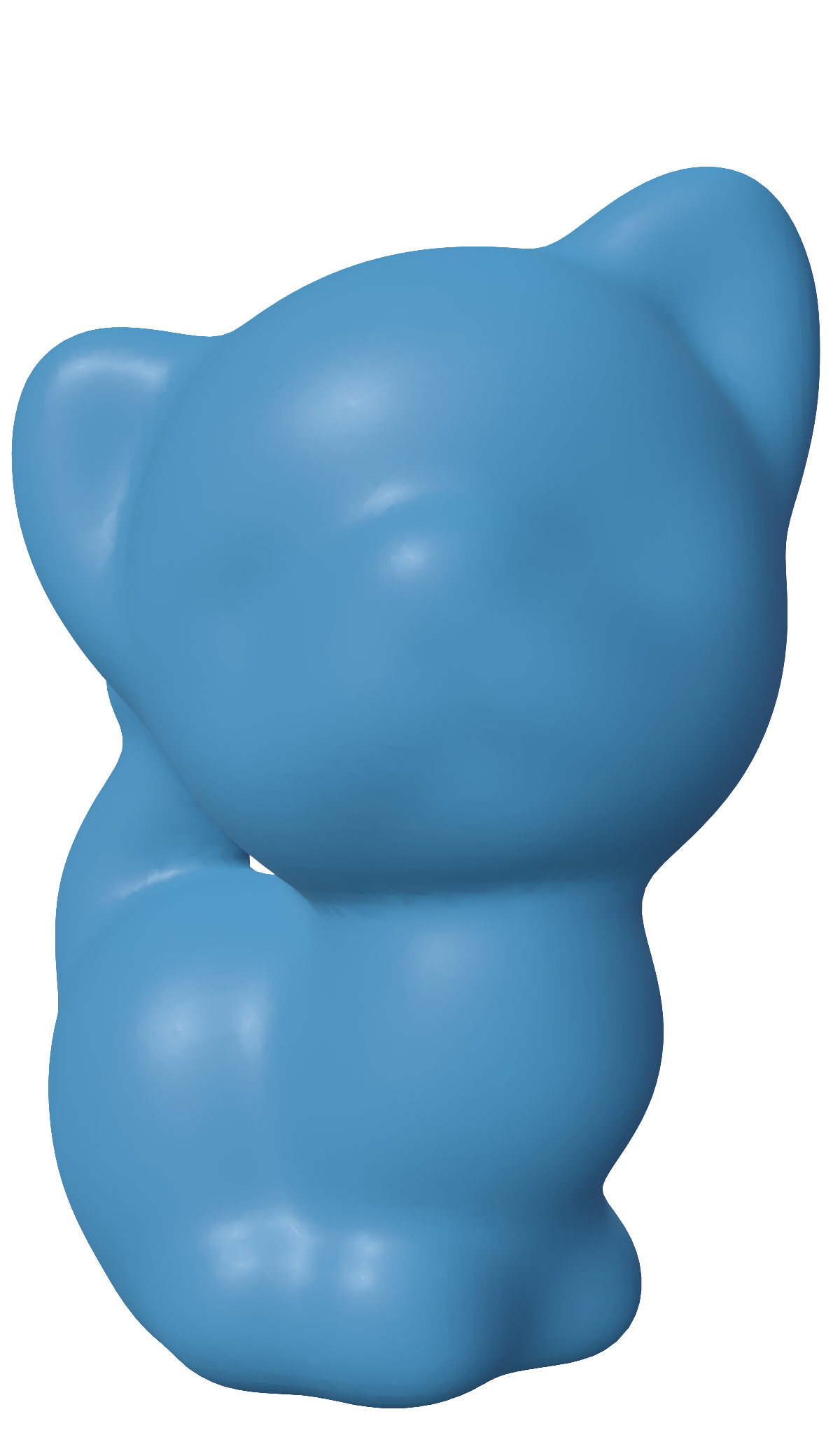}
	\includegraphics[height=4.1cm]{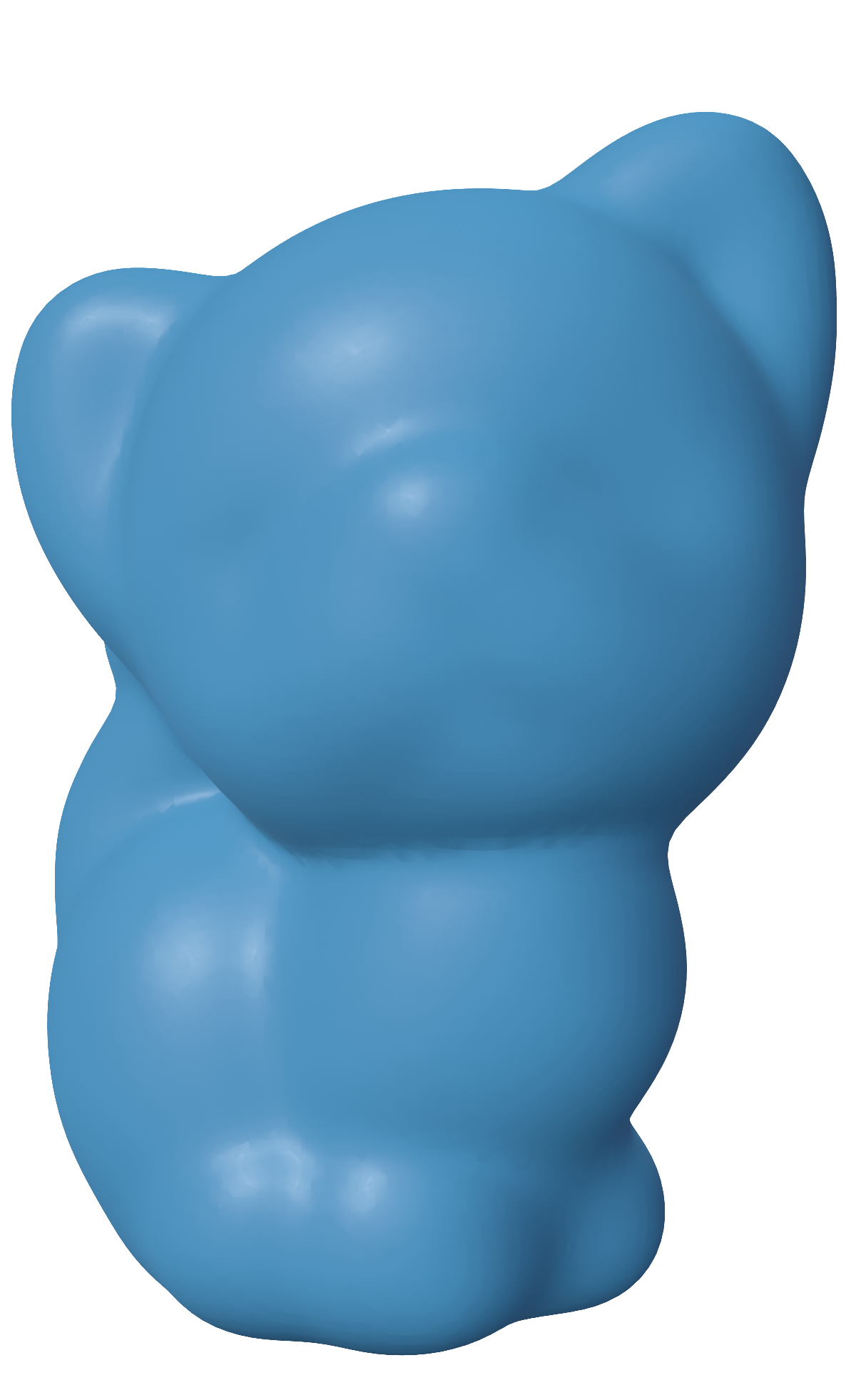}
	\includegraphics[height=4.1cm]{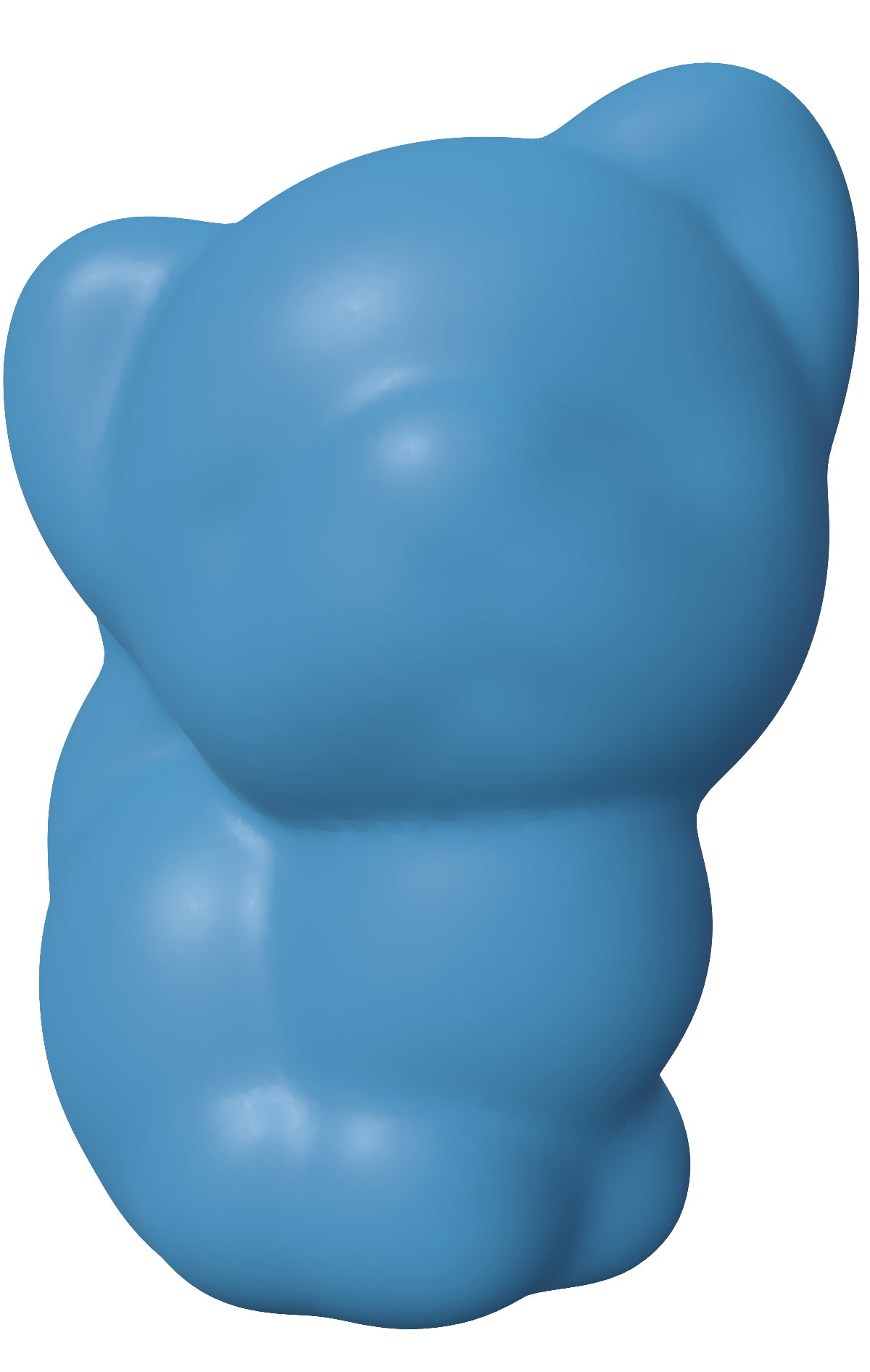}
	
	\includegraphics[height=3.425cm, trim=0.3cm 0cm 0.25cm 0,clip]{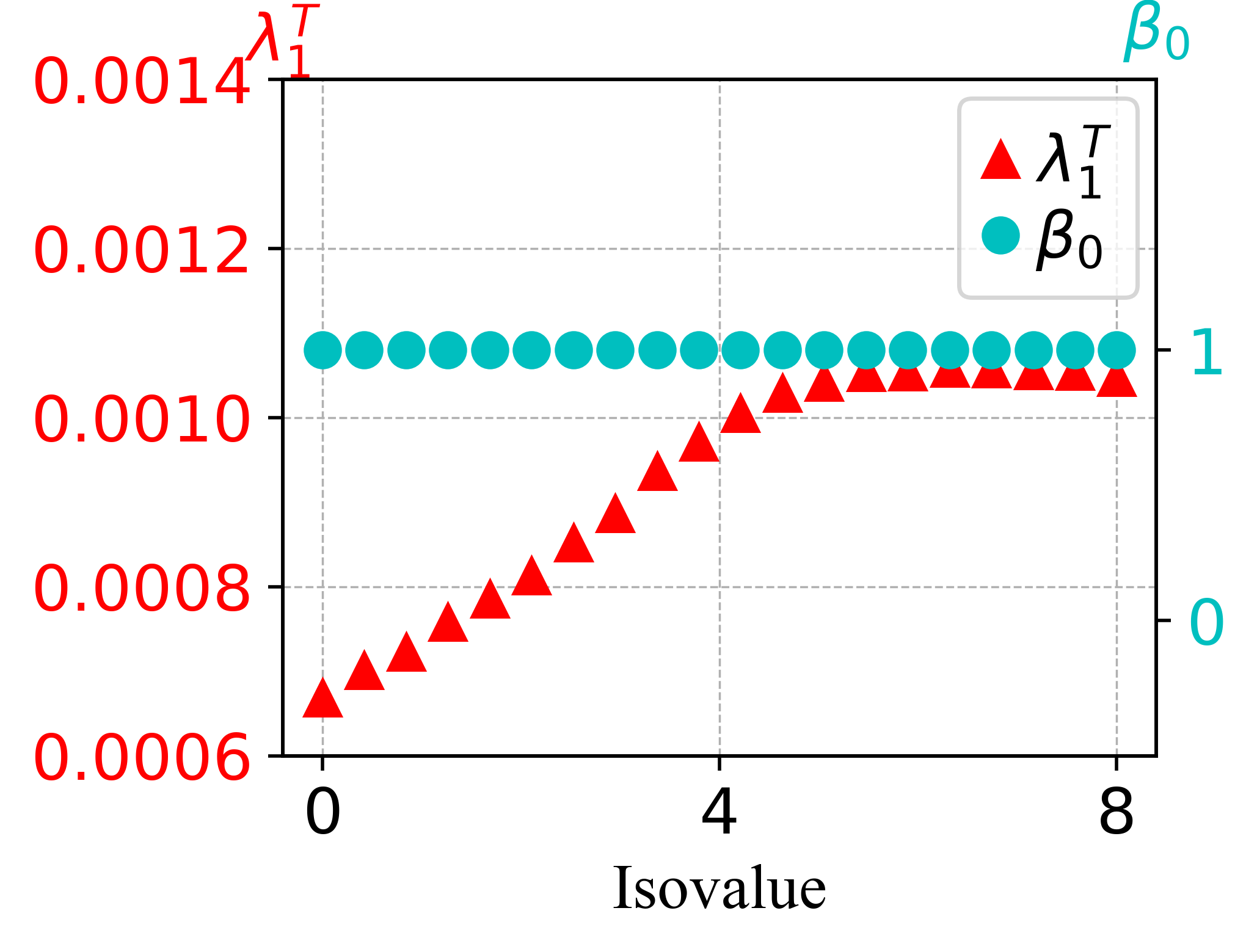}
	\includegraphics[height=3.425cm, trim=0.3cm 0cm 0.25cm 0,clip]{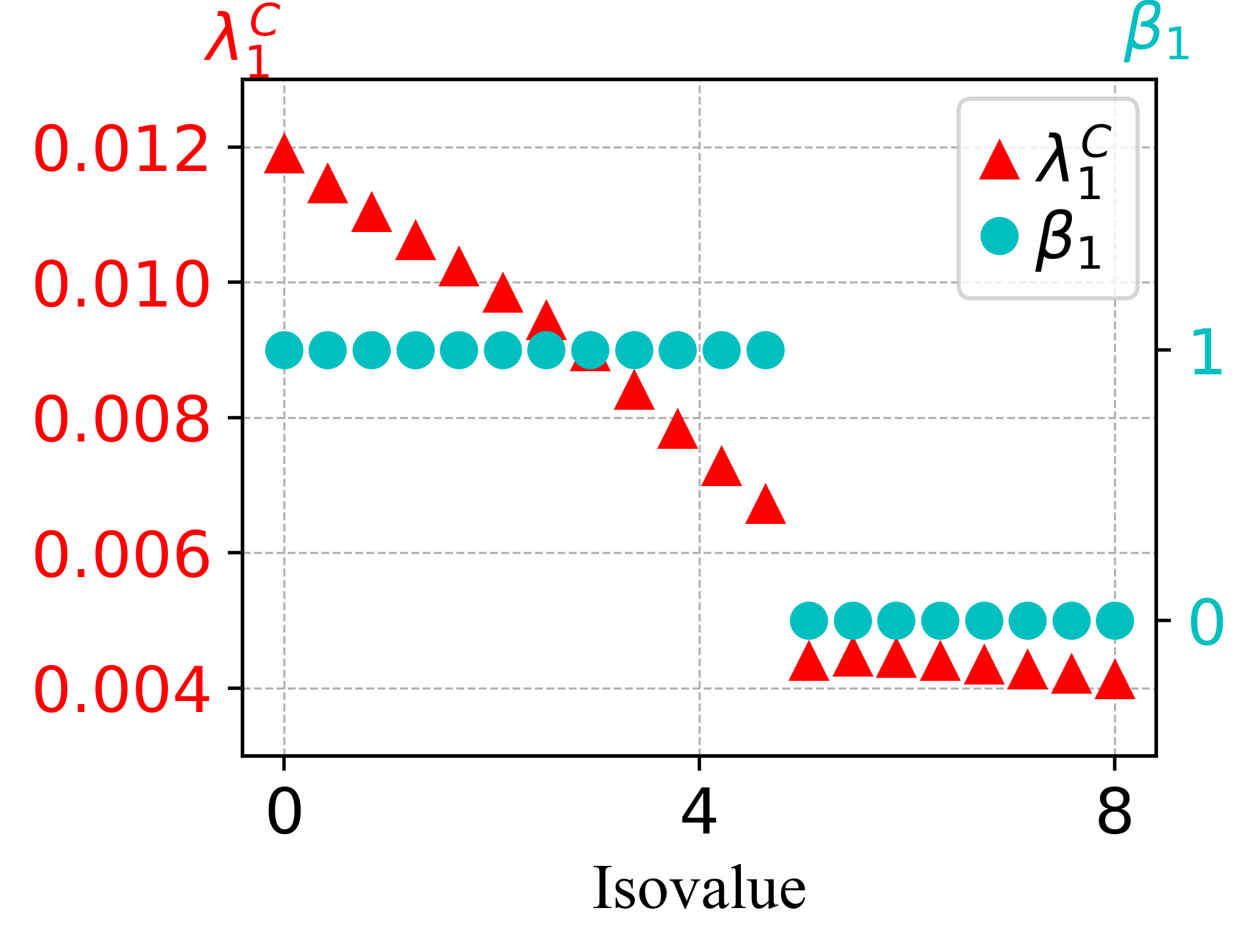}
	\includegraphics[height=3.425cm, trim=0.3cm 0cm 0.25cm 0,clip]{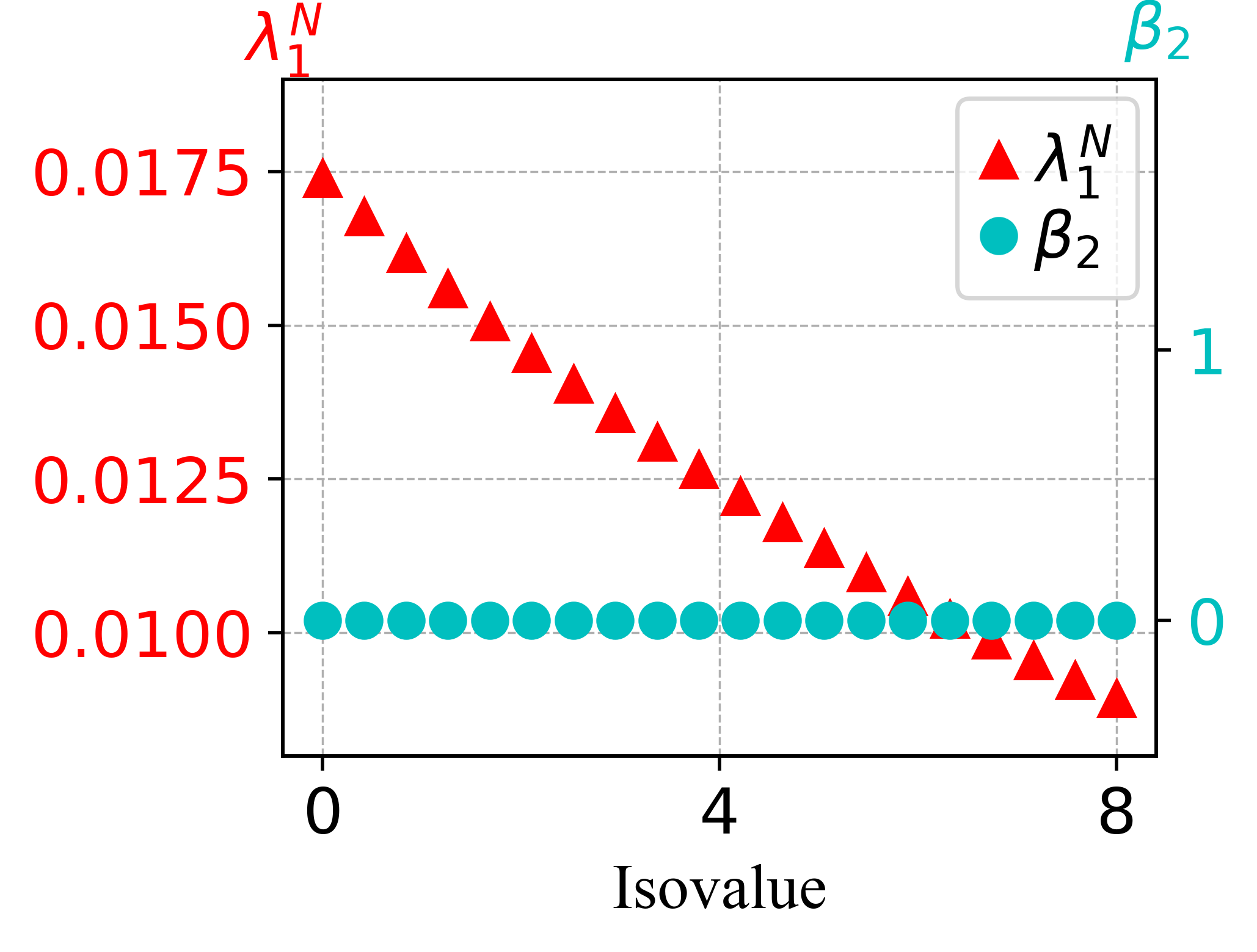}
 \caption{First row: Snapshots of evolving manifolds for the kitten model. Second row: Changes in Betti numbers $\beta_0$, $\beta_1$, $\beta_2$, and the first non-zero eigenvalues in T, C, N along $20$ evenly spaced isovalues from $0$ to $8$. Here the first shape in the top first row corresponds to isovalue $0$ and the last shape in the first row corresponds to isovalue $8$. $\lambda^T_1$, $\lambda^C_1$ and $\lambda^N_1$ are the first non-zero eigenvalues in the set T, C, N, respectively. The signed distance function generated from the original Kitten model is used as the level set function.}
	\label{fig.emfld.kitten}
\end{figure}

\begin{figure}[h]
	\centering
	\includegraphics[height=3cm]{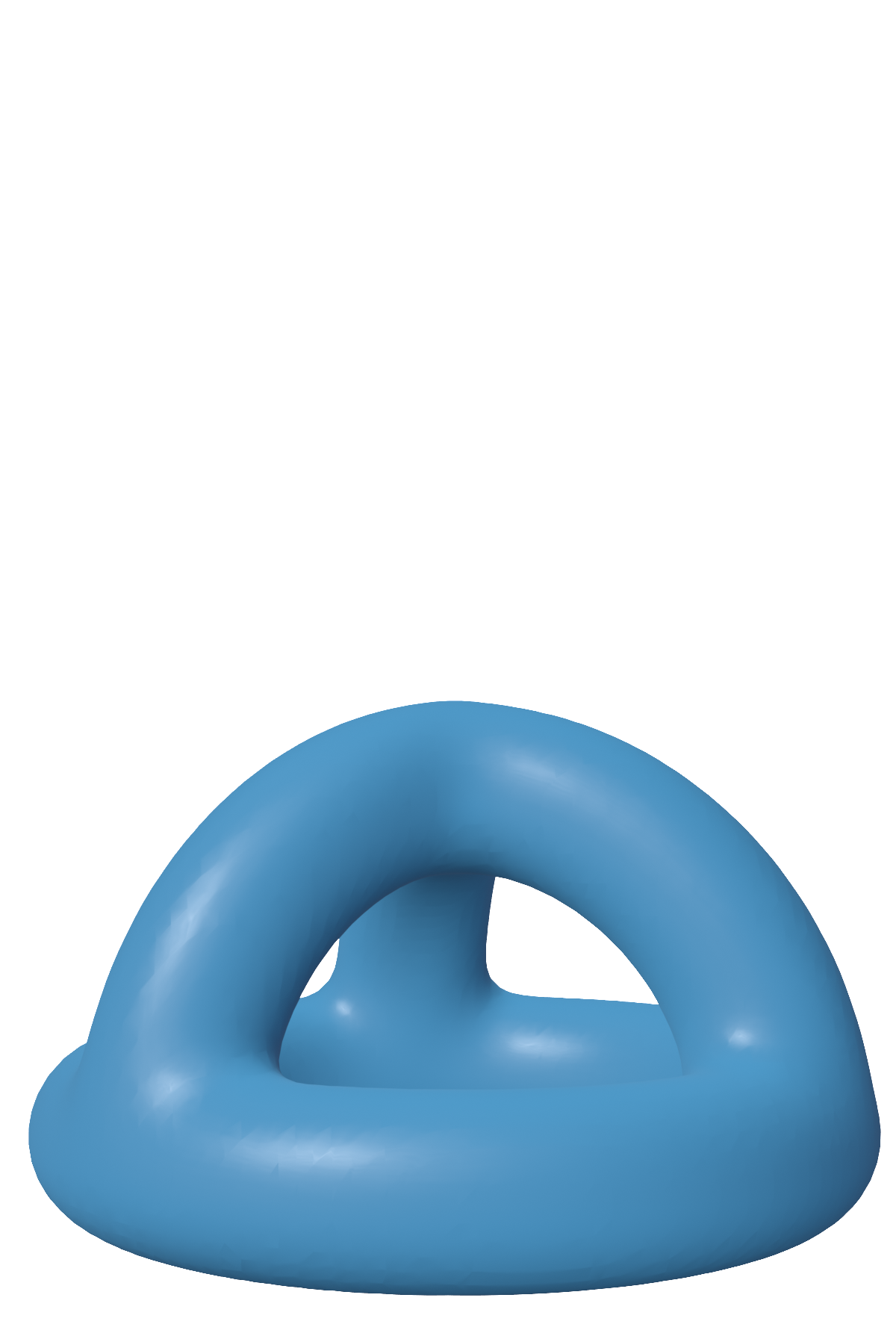}
	\includegraphics[height=3cm]{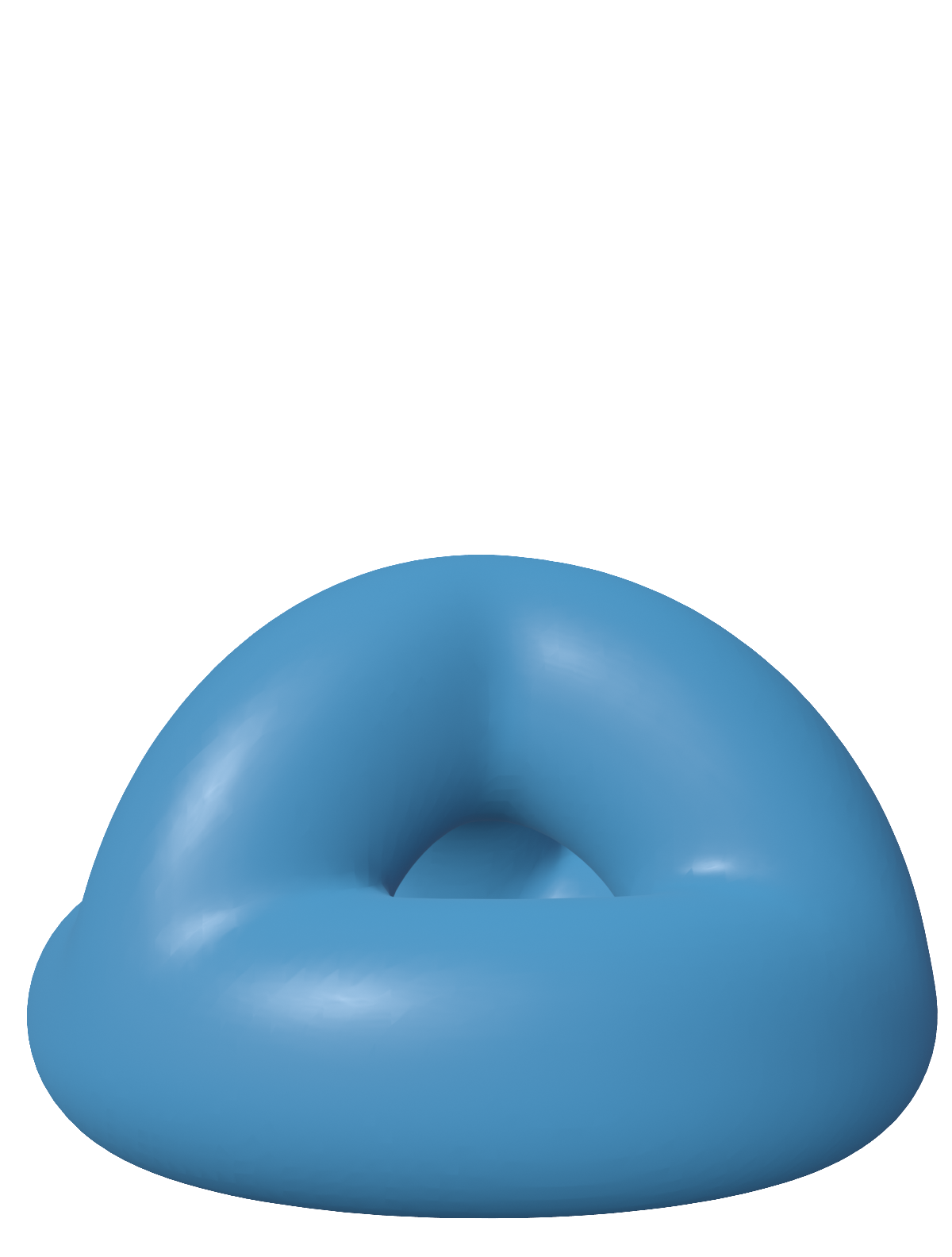}
	\includegraphics[height=3cm]{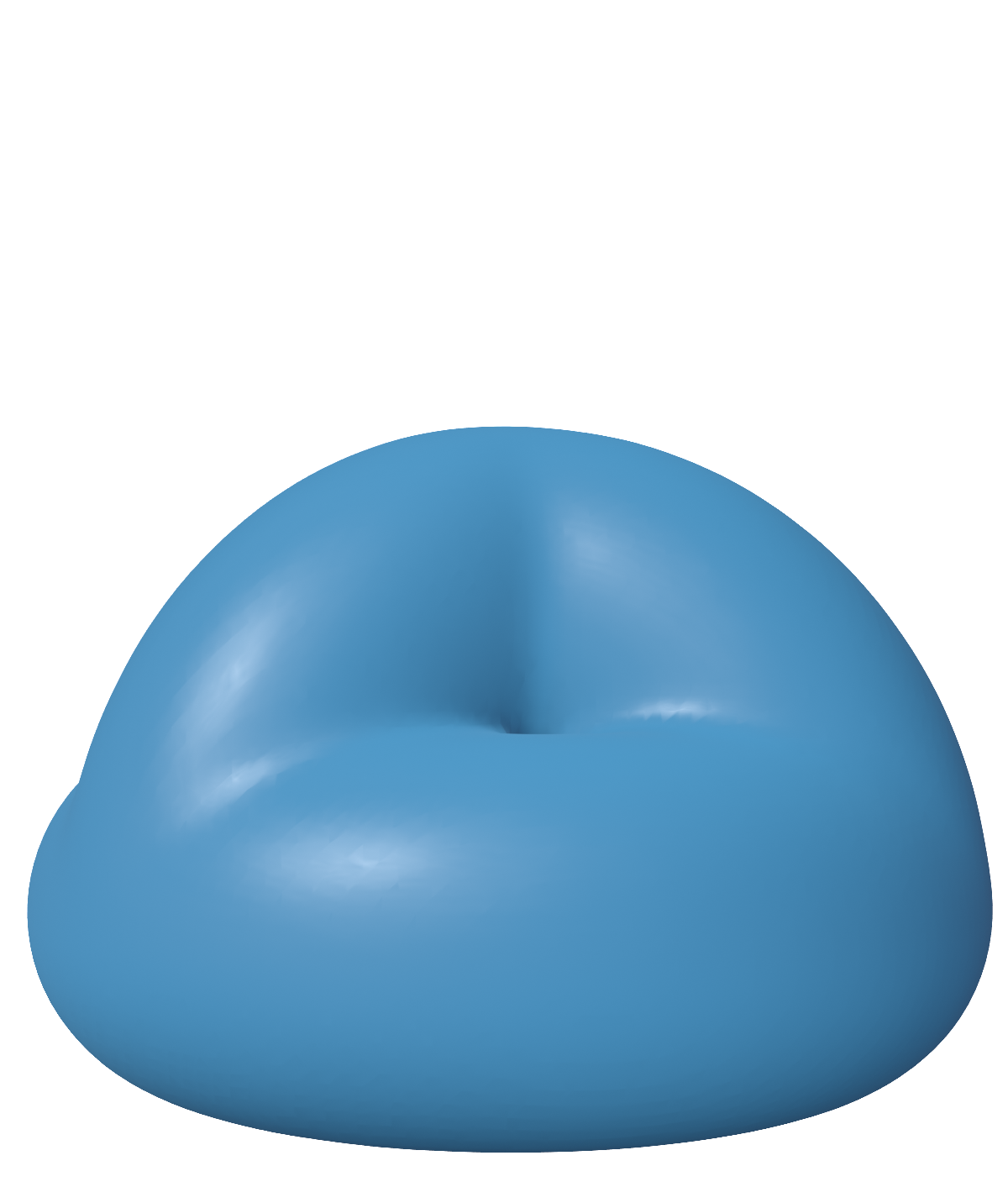}
	\includegraphics[height=3cm]{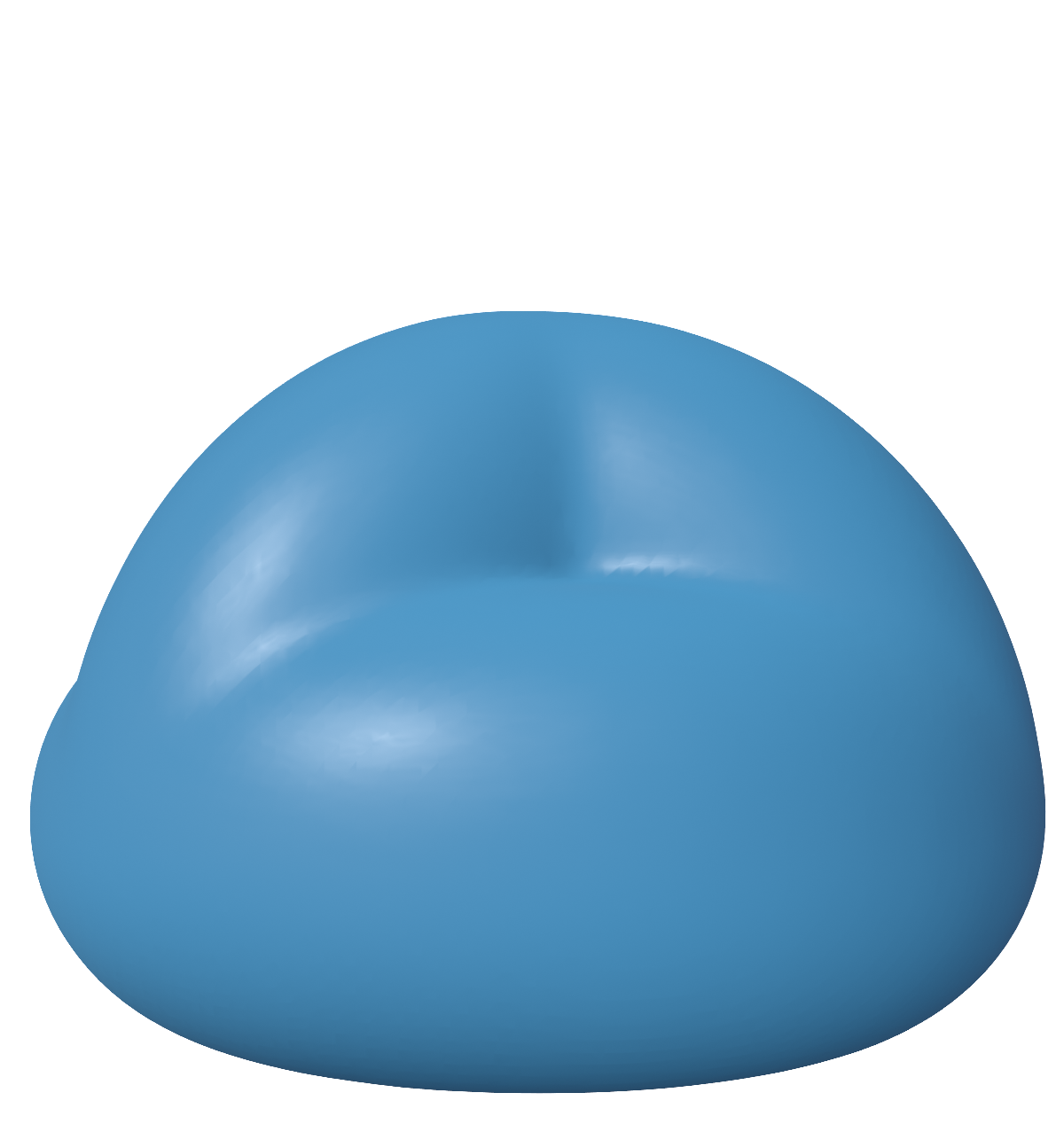}
	\includegraphics[height=3cm]{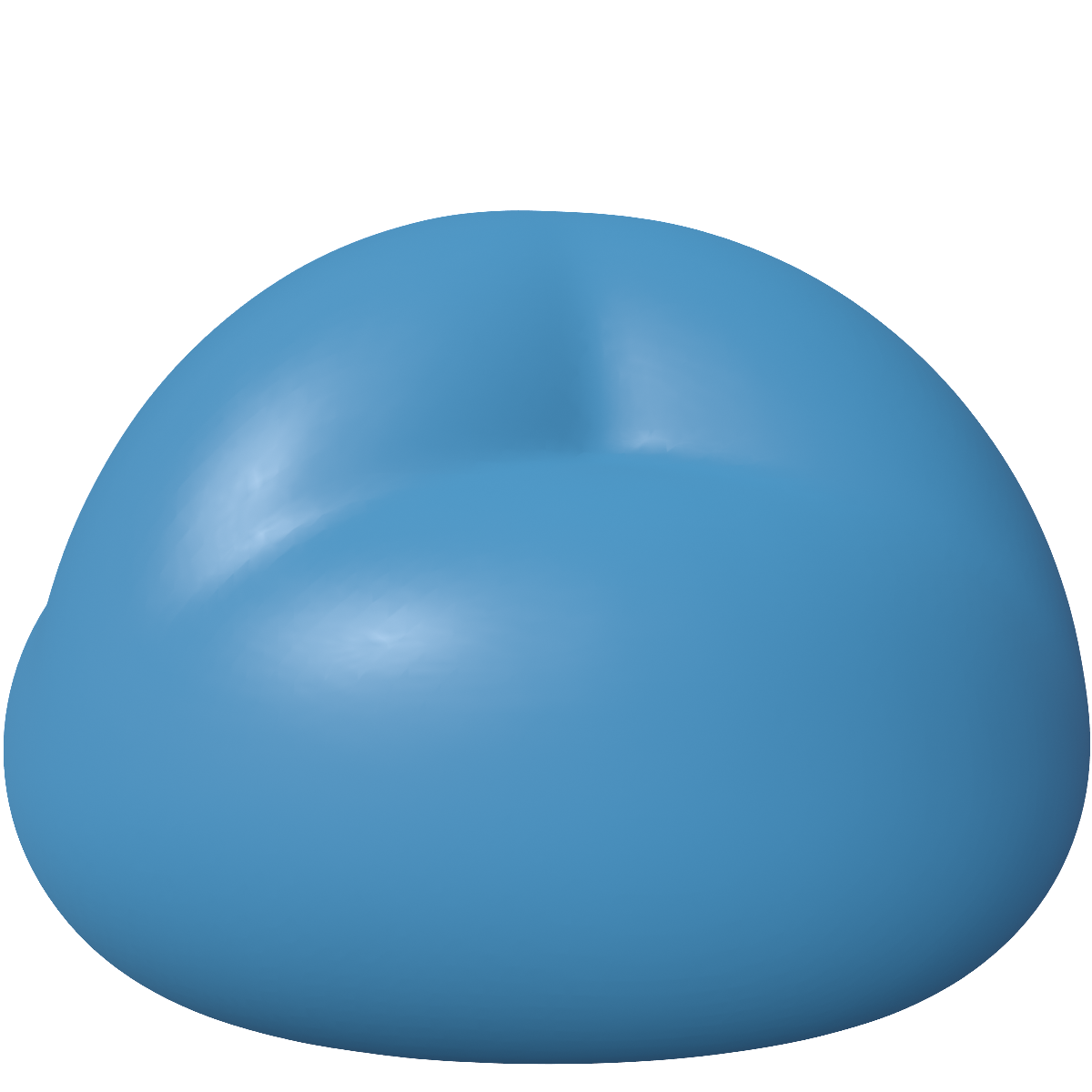}
	
	\includegraphics[height=3.425cm, trim=0.3cm 0cm 0.25cm 0,clip]{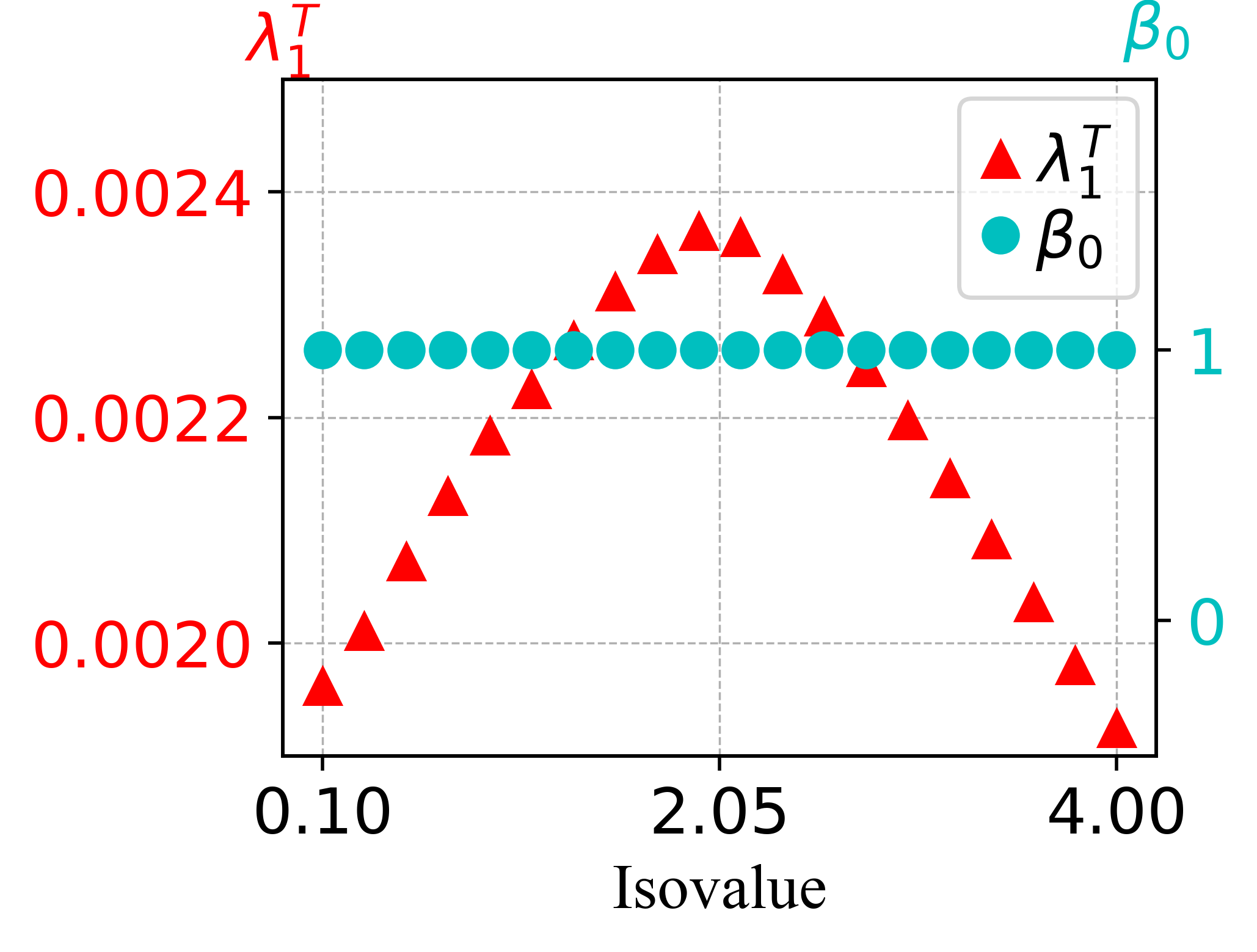}
	\includegraphics[height=3.425cm, trim=0.3cm 0cm 0.25cm 0,clip]{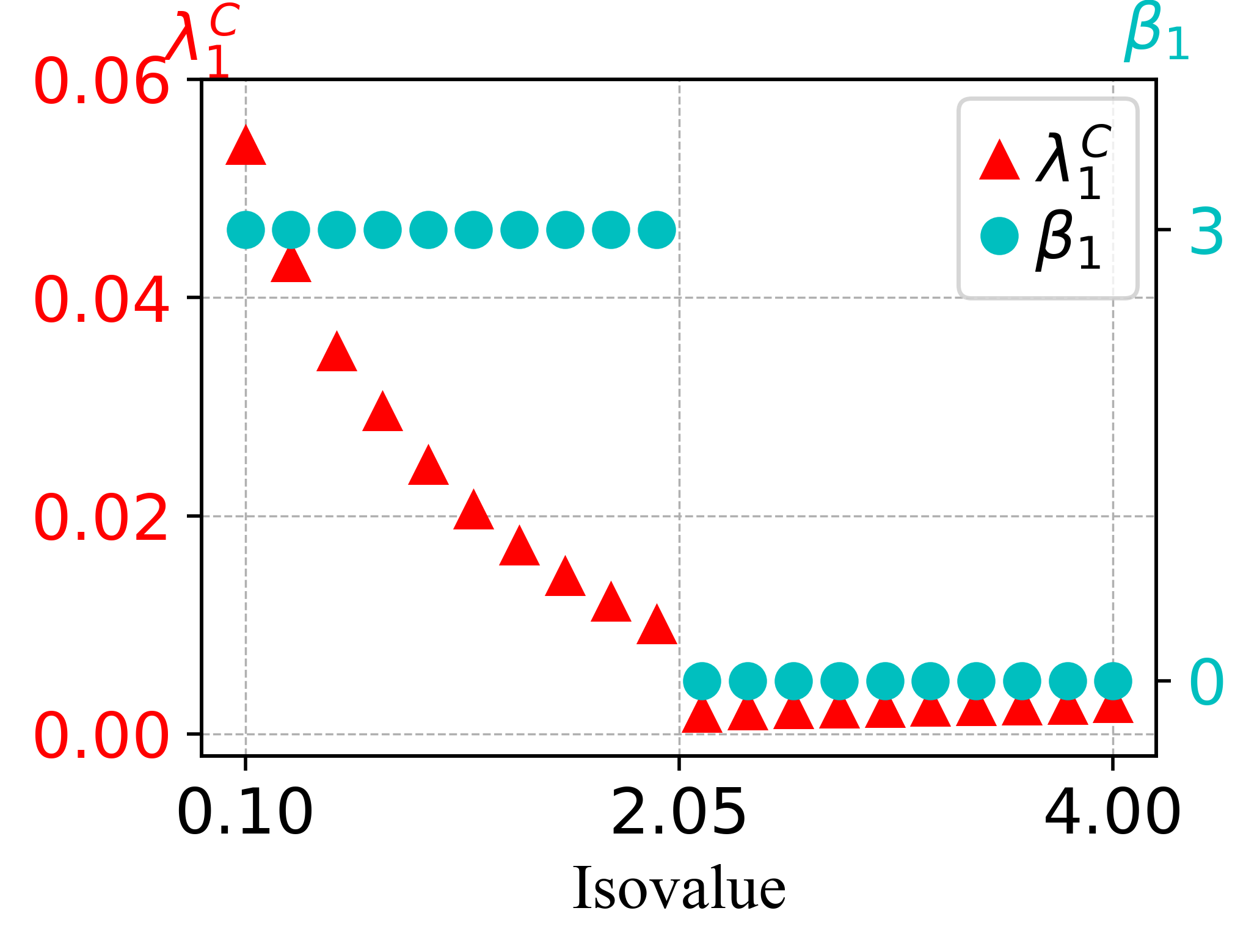}
	\includegraphics[height=3.425cm, trim=0.3cm 0cm 0.25cm 0,clip]{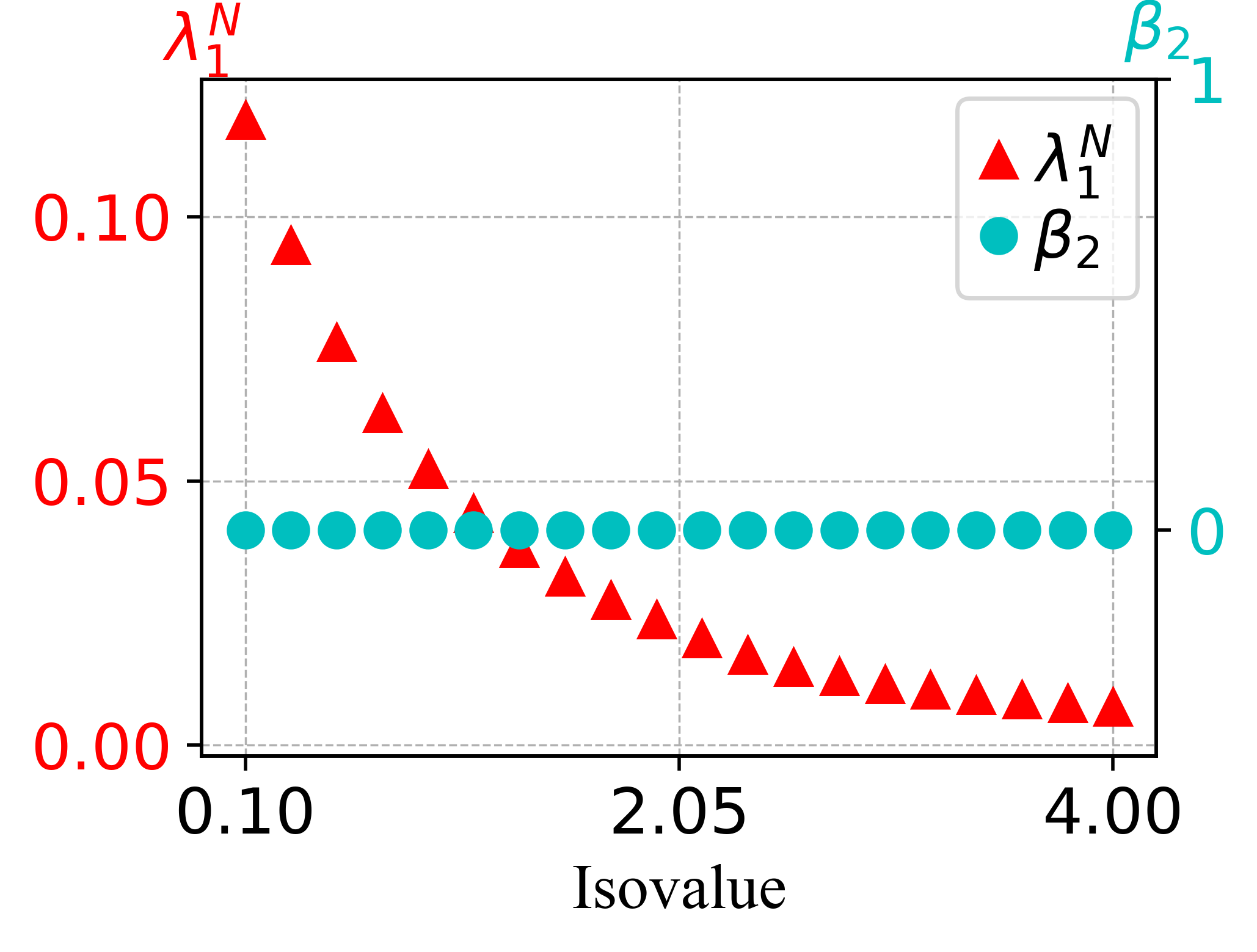}
	\caption{First row: Snapshots of evolving manifolds for a genus 3 model. Second row: Changes in Betti numbers $\beta_0$, $\beta_1$, $\beta_2$ and the first non-zero eigenvalues in T, C, N along $20$ evenly spaced isovalues from $0.1$ to $4$. Here the first shape in the top first row corresponds to isovalue $0.1$ and the last shape in the first row corresponds to isovalue $4$. $\lambda^T_1$, $\lambda^C_1$ and $\lambda^N_1$ are the first non-zero eigenvalues in the set T, C, N, respectively. The signed distance function generated from a genus 3 shape is used as the level set function.}
	\label{fig.emfld.g3}
\end{figure}

\begin{figure}[h]
	\centering
	\includegraphics[height=2.7cm]{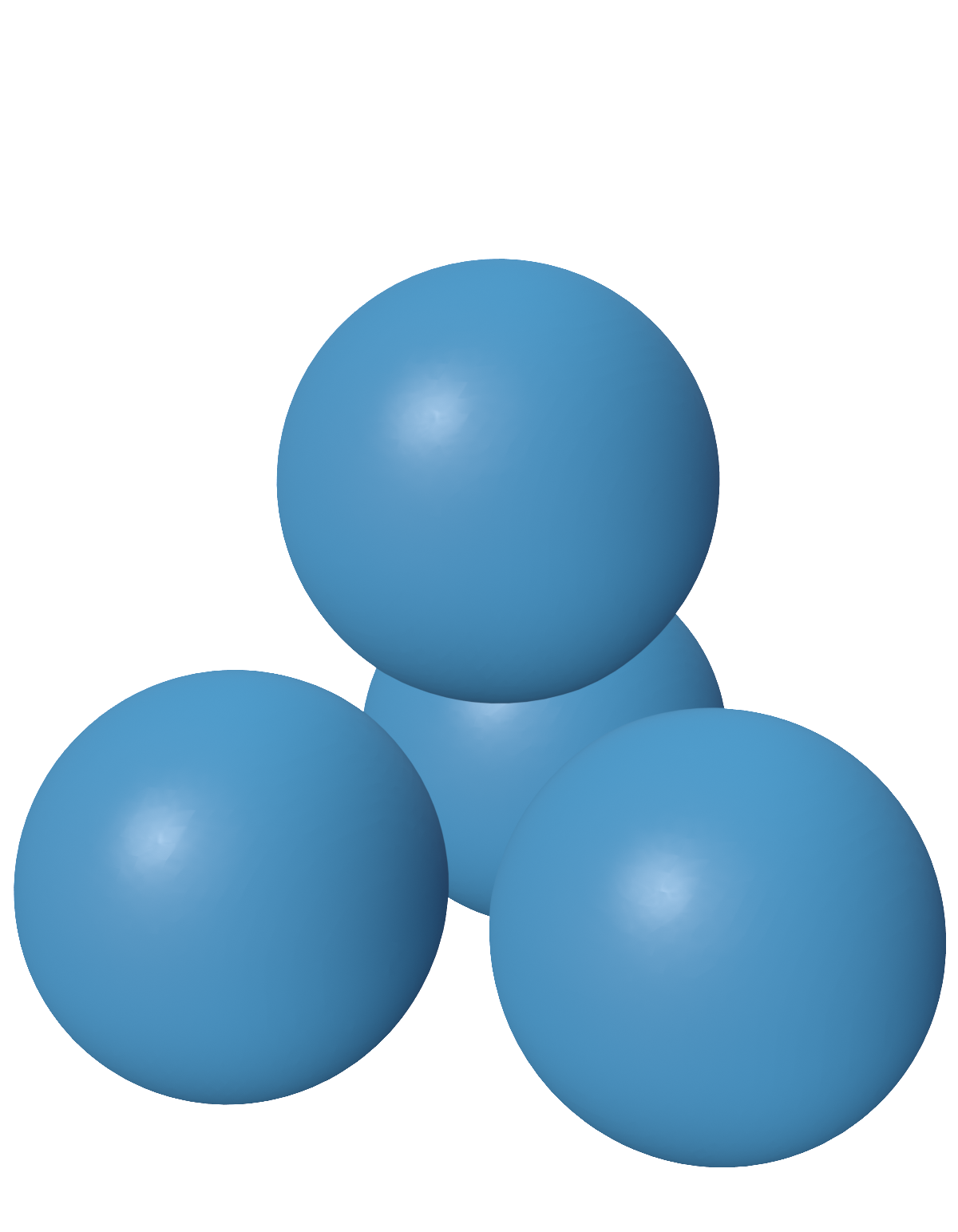}
	\includegraphics[height=2.7cm]{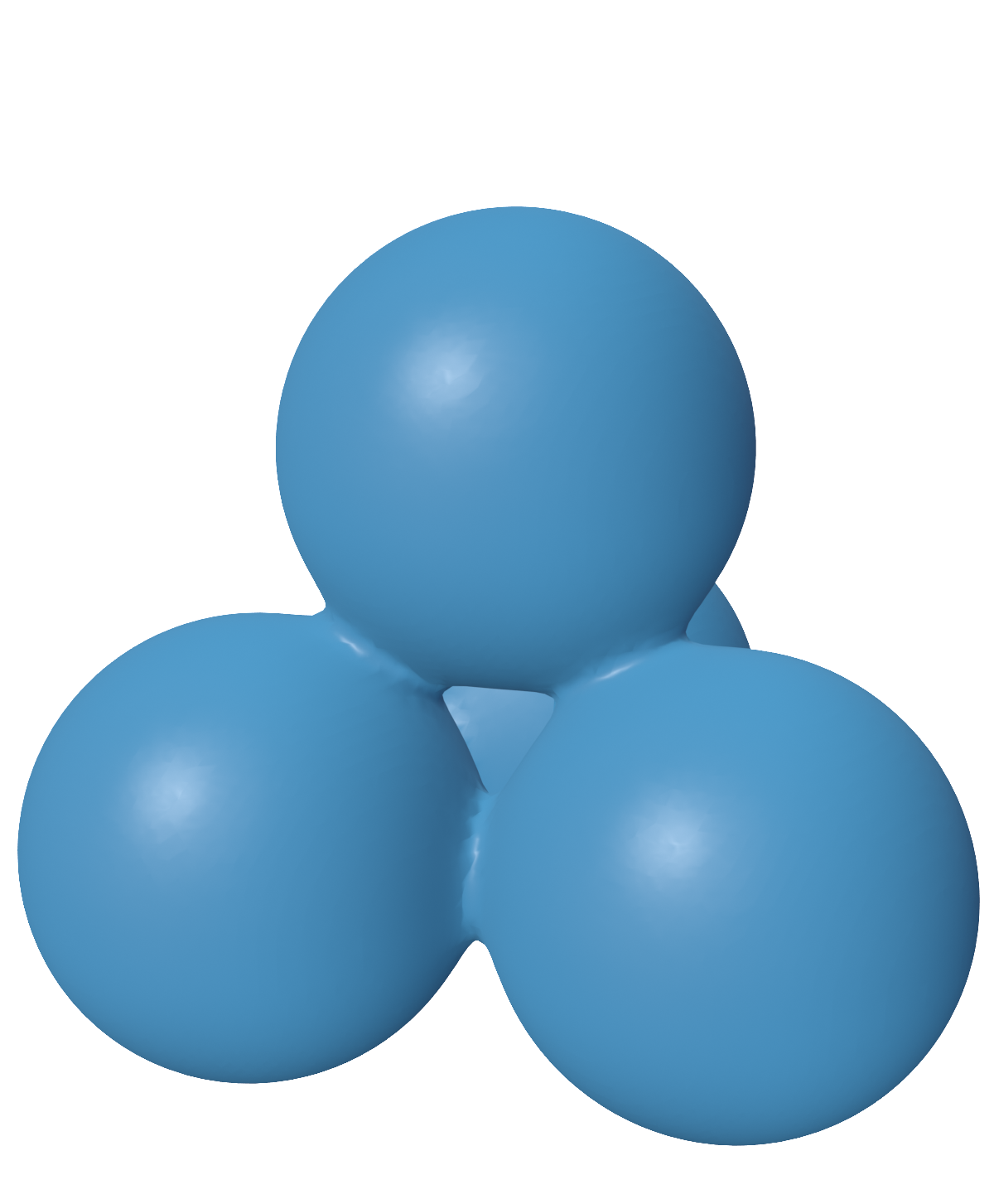}
	\includegraphics[height=2.7cm]{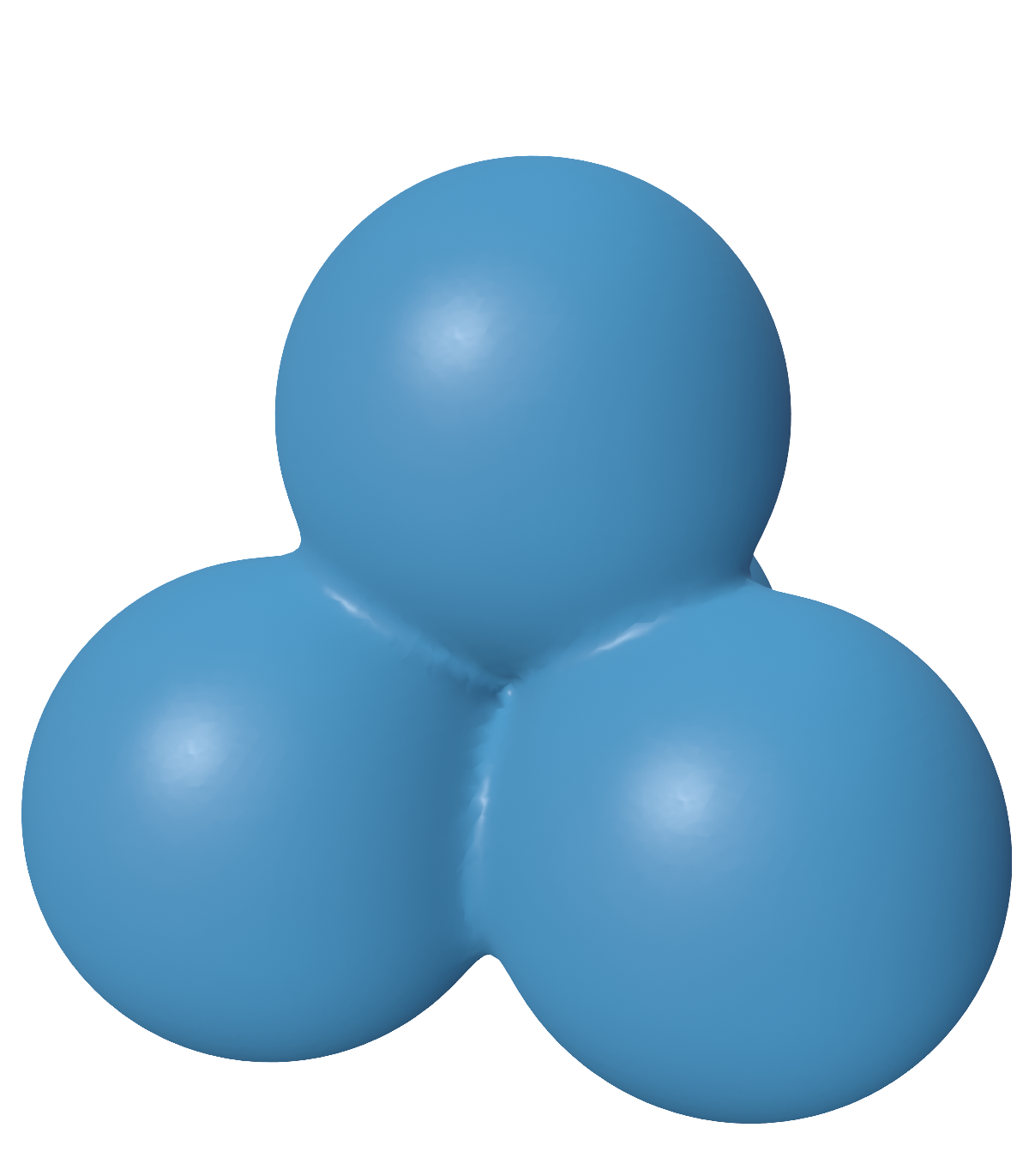}
	\includegraphics[height=2.7cm]{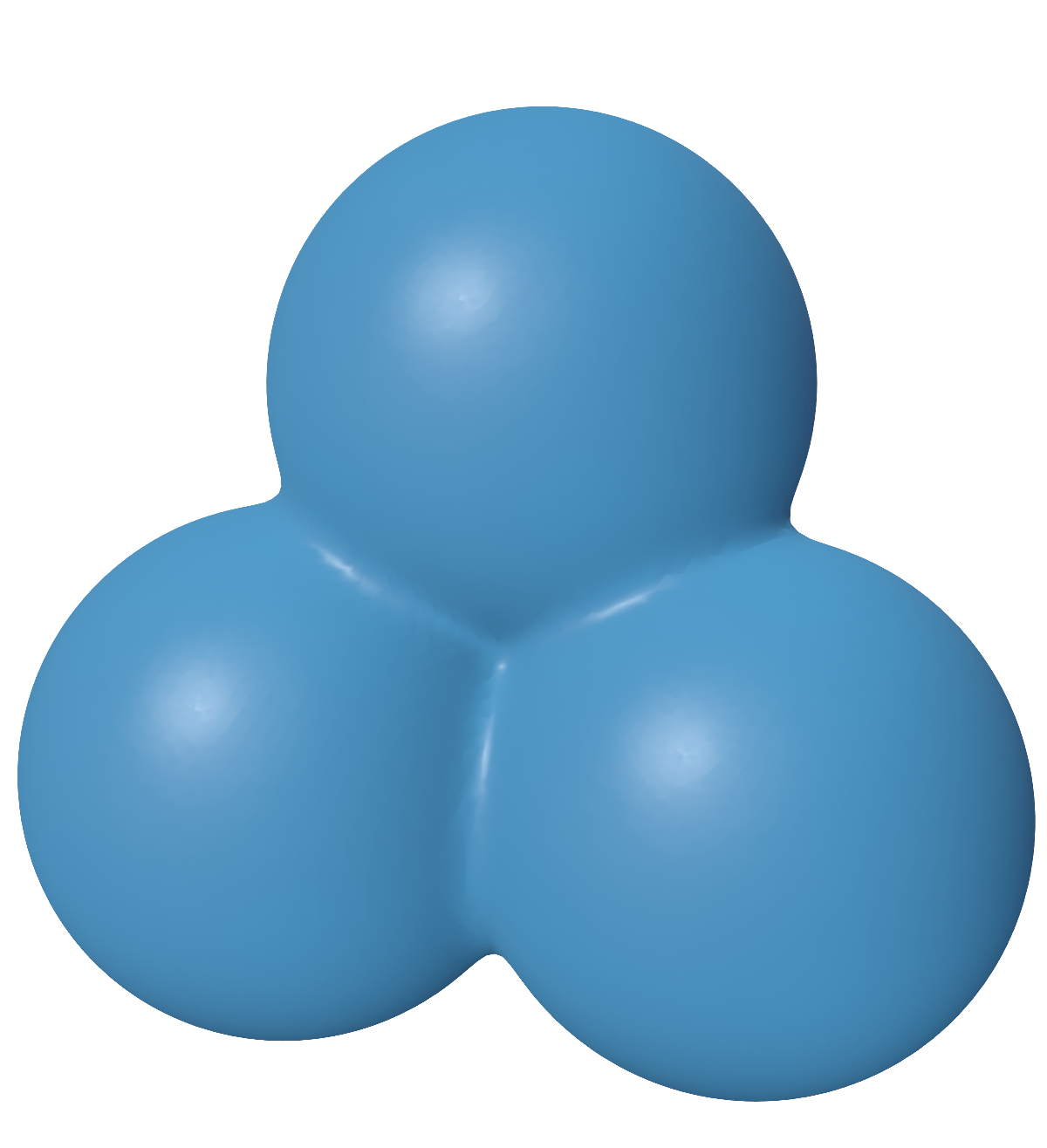}
	\includegraphics[height=2.7cm]{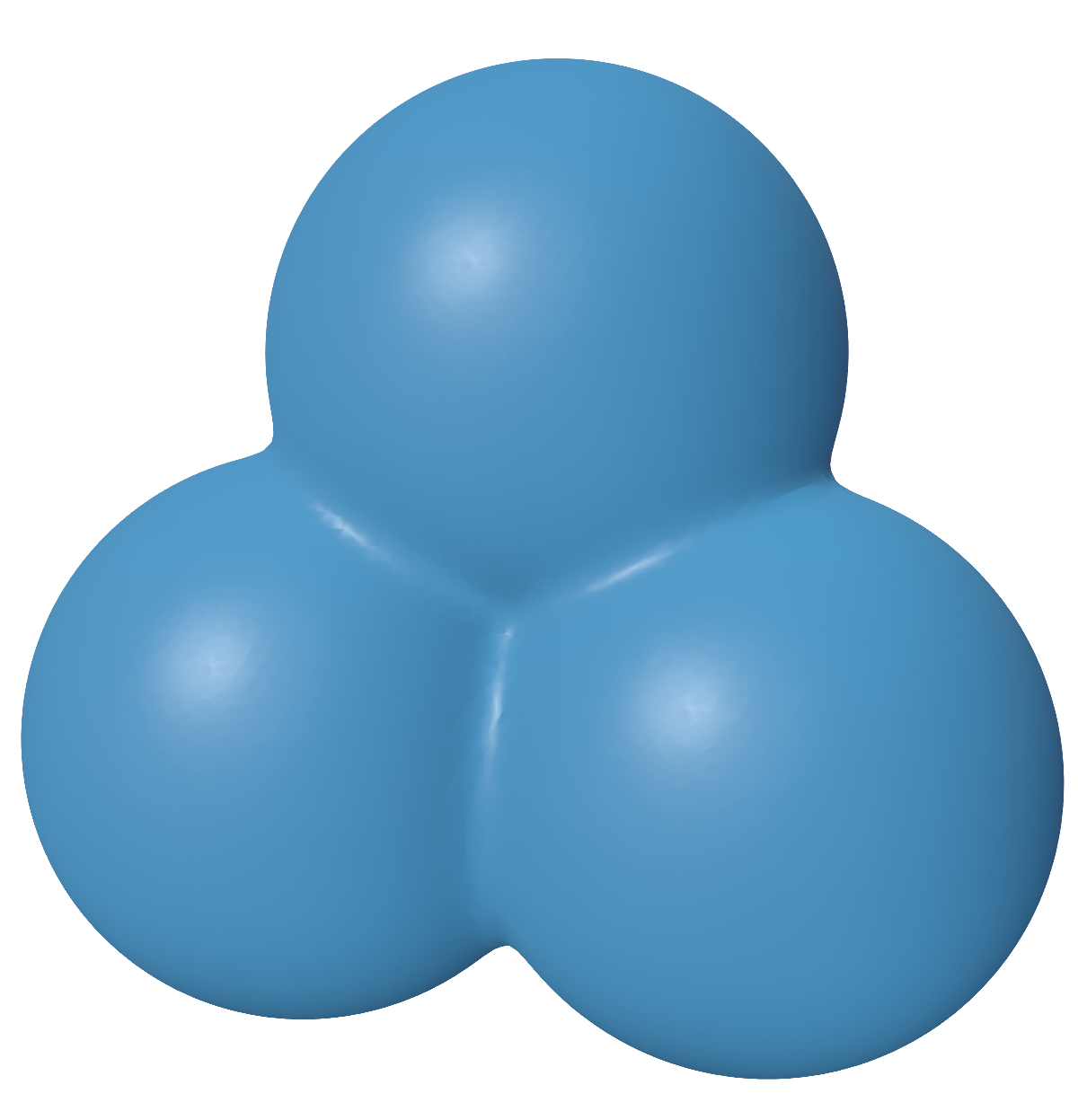}

        \includegraphics[height=3.425cm, trim=0.3cm 0cm 0.25cm 0,clip]{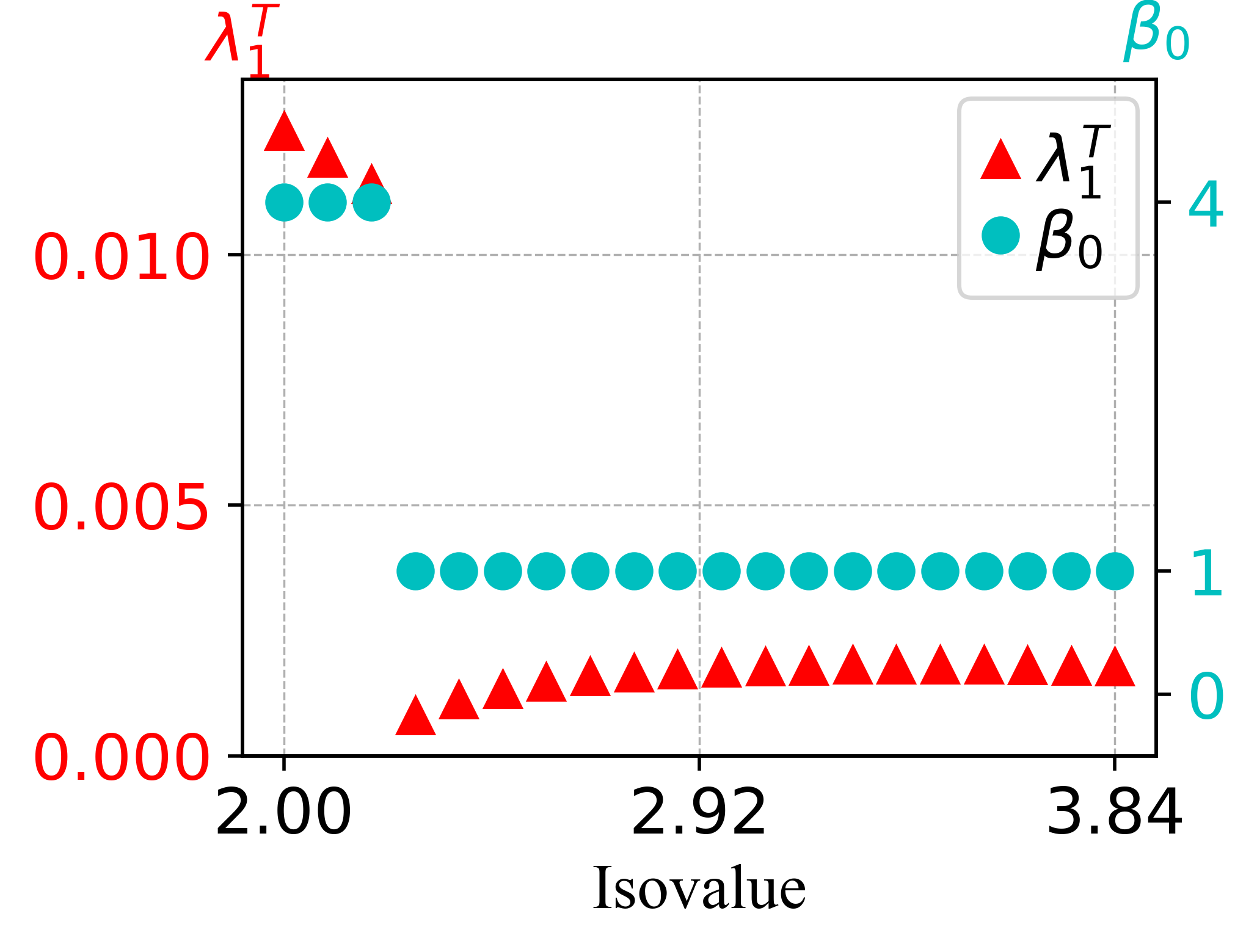}
	\includegraphics[height=3.425cm, trim=0.3cm 0cm 0.25cm 0,clip]{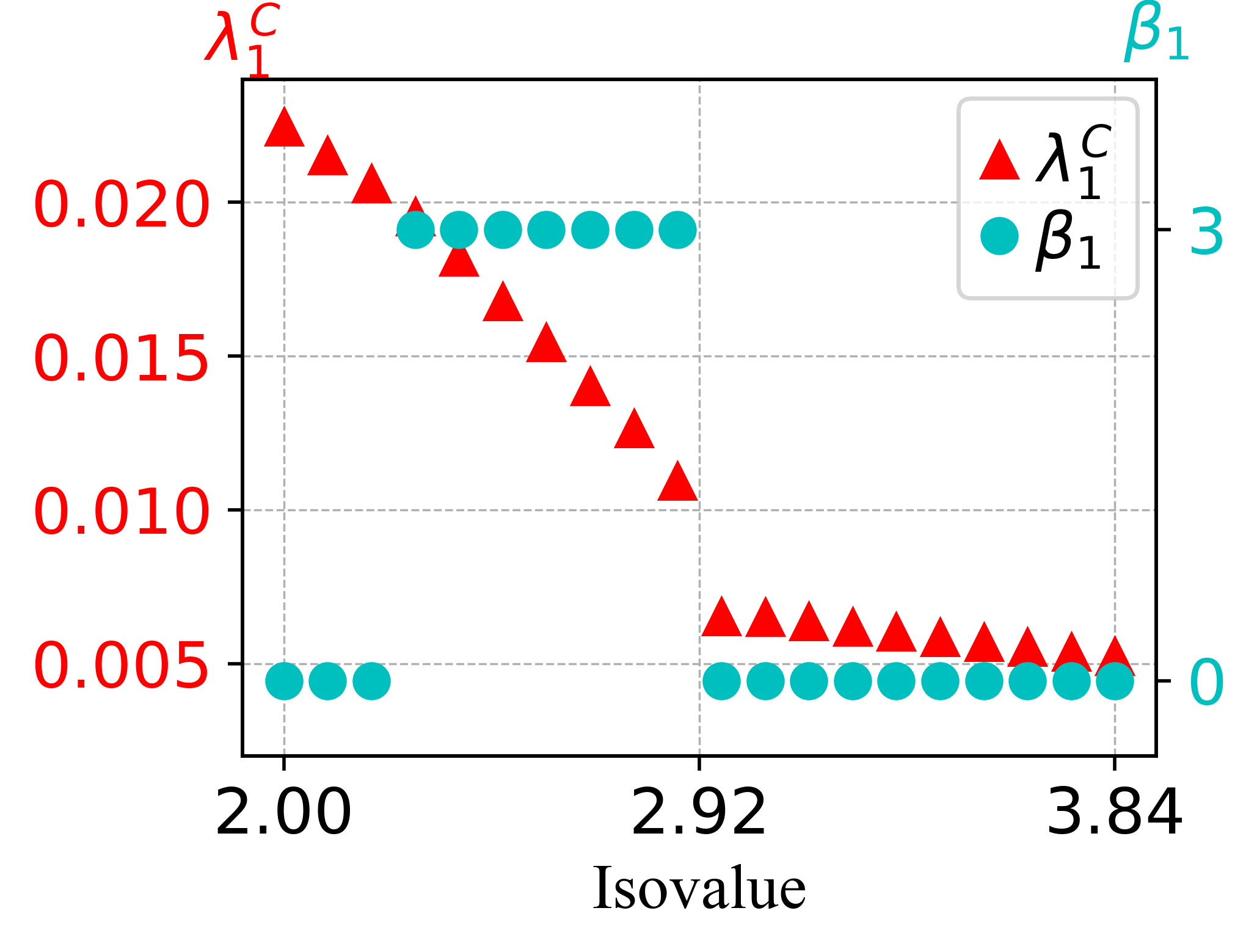}
	\includegraphics[height=3.425cm, trim=0.3cm 0cm 0.25cm 0,clip]{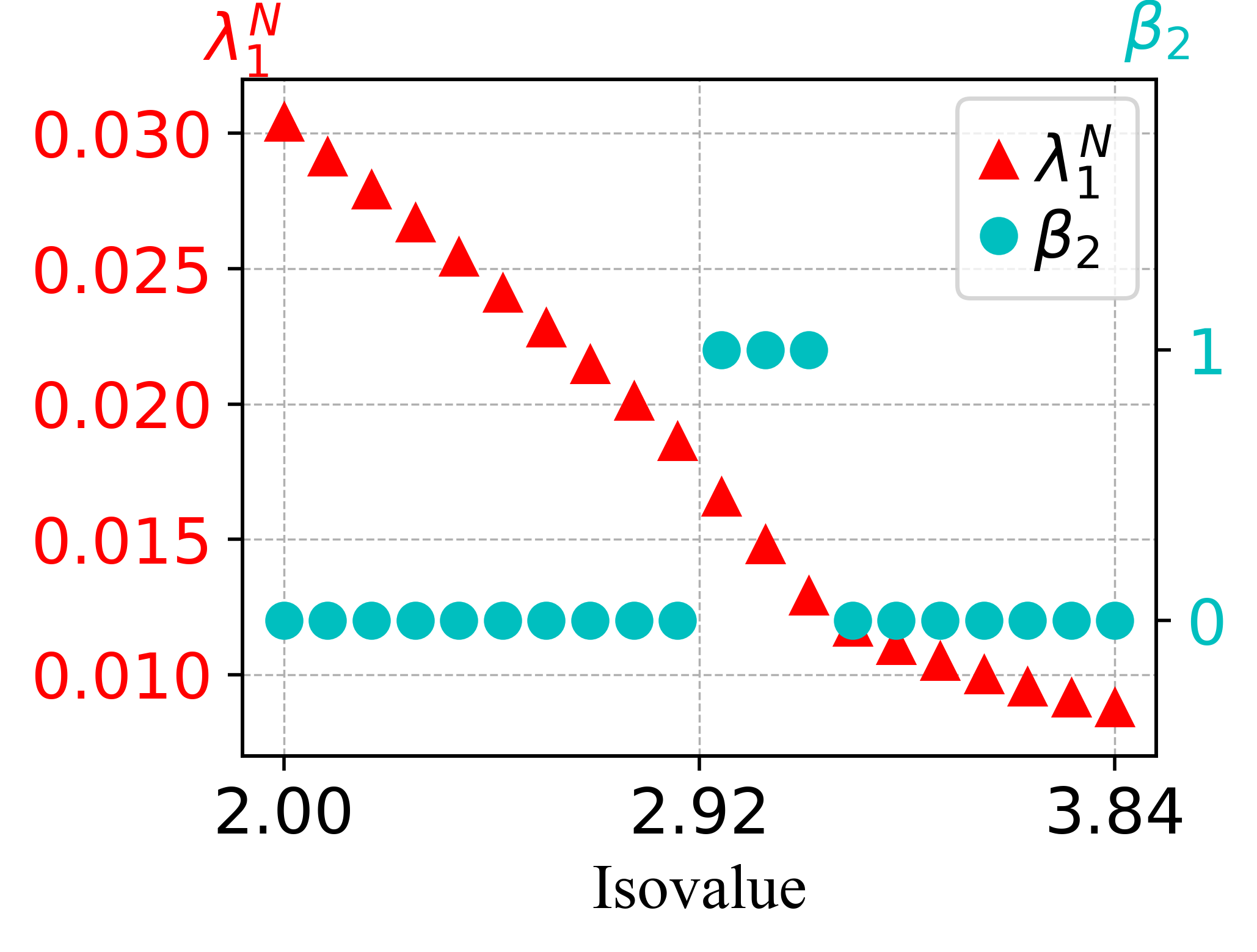}
	\caption{First row: Snapshots of evolving manifolds for a four-ball model. Second row: Changes in Betti numbers $\beta_0$, $\beta_1$, $\beta_2$ and the first non-zero eigenvalues in T, C, N along $20$ evenly spaced isovalues from $2$ to $3.84$. Here the first shape in the top first row corresponds to isovalue $2$ and the last shape in the first row corresponds to isovalue $3.84$. $\lambda^T_1$, $\lambda^C_1$ and $\lambda^N_1$ are the first non-zero eigenvalues in the set T, C, N, respectively. The signed distance function generated from four seperate balls centered at the vertices of a tetrahedron is used as the level set function.}
	\label{fig.emfld.4balls}
\end{figure}

We now present some examples of evolving manifolds and show results for the spectral analysis of their persistent Laplacians. In particular, we focus on the changes of Betti numbers $\beta_0$, $\beta_1$ and $\beta_2$ and the first non-zero eigenvalues $\lambda^T_1, \lambda^C_1$ and $\lambda^N_1$ of the 0-persistent BIG Laplacians in the set T, C and N, respectively, as introduced in Sec.~\ref{sec.spectraLaplacians}. Four models are considered, including the Bimba model, the kitten model, a genus-3 model, and a four-ball model. For each model, we show on the top row snapshots of evolving manifolds at five evenly spaced isovalues in a chosen interval, and on the bottom row the changes in Betti numbers and the first non-zero eigenvalues $\lambda^T_1, \lambda^C_1$ and $\lambda^N_1$. All the evolving manifolds are generated using isovalues of the signed distance function (SDF) from the original surface model, given as the 0-isosurface of the SDF. As we show below, these values from the evolution of manifolds provide rich information than considering just a single manifold. The discontinuity of these variables indicates the topological changes occurring during the evolution process, and the monotonicity of these non-zero eigenvalues reveals the geometric changes. 

The results for the Bimba model are presented in Fig.~\ref{fig.emfld.bimba} with an isovalue interval $[0, 0.2]$. As there is no topological change happening in the evolution process, all Betti numbers $\beta_0$, $\beta_1$ and $\beta_2$ remain constant, and $\lambda^T_1, \lambda^C_1$ and $\lambda^N_1$ are continuous throughout the whole process. Both $\lambda^C_1$ and $\lambda^N_1$  decrease as the isovalue increases. 

Fig.~\ref{fig.emfld.kitten} illustrates the results for the kitten mode with one tunnel formed by its tail. The isovalue interval $[0, 8]$ is considered. One can see all variables are continuous during the evolution process except that $\beta_1$ and $\lambda^C_1$ both drop at the same isovalue, where $\beta_1$ changes from $1$ to $0$. This happens due to the disappearance of the tunnel. In addition, $\lambda^T_1$ increases at the beginning, and then slows down its rate of increase at the isovalue after the tunnel disappears, and $\lambda^C_1$ and $\lambda^N_1$ decrease during the evolution process.

Note that there are also tunnels in the evolving manifolds for the genes-3 model, as we expected, a similar phenomenon can also be observed in Fig.~\ref{fig.emfld.g3} for the change of the Betti numbers and the first non-zero eigenvalues. The isovalue interval $[0.1, 4]$ is considered for this model. The disappearance of the three tunnels leads to a drop of $\beta_1$ from $3$ to $0$ and also a drop of $\lambda^C_1$. $\lambda^T_1$ initially increases, and then changes its behavior to decrease after the tunnels vanish. The evolution process results in a decrease in $\lambda^C_1$ and $\lambda^N_1$, just as the previous two models.

The evolving process of the four-ball model with isovalue interval $[2, 3.84]$, see Fig.~ \ref{fig.emfld.4balls}, leads to discontinuities in all Betti numbers and the first non-zero eigenvalues. As the four separate components merge in the evolution, $\beta_0$ changes from $4$ to $1$, along with a drop in $\lambda^T_1$ at the same isovalue. In addition, $\beta_1$ increases from $0$ to $3$ due to the appearance of three tunnels when the merge happens and then decreases to $0$ after the disappearance of all tunnels. The non-zero eigenvalue $\lambda^C_1$ has a drop that occurs when the tunnel vanishes, however, it is continuous when the tunnels are formed. This suggests that the continuity of $\lambda^C_1$ is only related to the death but not the birth of tunnels. One can also observe a slowdown in the rate of change of $\lambda^T_1$ following the disappearance of all tunnels. As the isolate increases further, a cavity occurs in the manifold, resulting in an increase of $\beta_2$ from $0$ to $1$ and finally a decrease from $1$ to $0$ after the cavity disappears. This topological change can also be observed in $\lambda^N_1$, where $\lambda^N_1$ becomes non-differentiable.

As illustrated by these models, changes in Betti numbers $\beta_0$, $\beta_1$ and $\beta_2$ and the first non-zero eigenvalues $\lambda^T_1, \lambda^C_1$ and $\lambda^N_1$ not only reflect the changes in topology, but also characterizes the changes in geometry for the evolution of manifolds. The rich information revealed by these variables leads to potential applications in various topological data analysis tasks.


\section{Proof-of-Principle Experimentation}

In this section, we carry out a proof-of-principle experimental demonstration of the proposed persistent de Rham-Hodge theory. In this approach, the problem is defined on manifolds with boundaries. Appropriate boundary conditions are implemented to match actual topological dimensions. The resulting persistent Hodge Laplacians are solved to deliver the corresponding series of eigenvectors and eigenvalues at various scales. In this approach, we use these eigenvalues for machine learning predictions of protein-ligand binding affinity.  
The binding affinity describes the strength of protein-ligand interactions for each protein-ligand complex.

We consider two benchmark datasets, PDBbind-v2007 and PDBbind-v2016 \cite{liu2017forging}, to demonstrate the effectiveness of our framework in capturing the topological features of protein-ligand complexes. The datasets can be downloaded from http://pdbbind.org.cn/. These two PDBbind datasets provide collections of biomolecular complexes in Protein Data Bank (PDB) with experimentally a measured binding affinity for each protein-ligand complex, and are commonly used in various studies such as drug-discovery or molecular recognition, etc \cite{cang2018integration, liu2017forging, cang2018representability, su2018comparative,francoeur2020three, meng2021persistent, liu2023persistent}. We aim to build a machine learning model, by utilizing the topological and geometric features of the protein-ligand complexes generated using our persistent Hodge Laplacian (PHL) framework as inputs, for predicting the protein-ligand binding affinities. 

The biomolecular complexes in each PDBbind dataset are organized into three sets, including a general set, a refined set and a core set, with each set being a superset of the next. In our experiments, for each dataset, we use the refined set, excluding the core set, to train the predictive model for the binding affinities of the protein-ligand complexes in the core set. The PDBbind-v2007 dataset contains a total of 1,300 complexes with 1,105 in the refined set and 195 in the core set, while the PDBbind-v2016 dataset has a total of 4,057 complexes with 3,767 in the refined set and 290 in the PDBbind core set.

\subsection{Element specific discrete to continuum mapping}

The original datasets contain atomic names and coordinates, which are the so-called point cloud data. To generate manifold representations, we carry out the discrete to continuum mapping using the flexibility and rigidity index \cite{nguyen2019dg}. To compute the topological feature of each protein-ligand complex for the machine learning model, we use the element-specific approach \cite{cang2017topologynet}.  Specifically, we consider the pairwise interactions between element types that are commonly found in proteins and ligands, including Hydrogen (H), Carbon (C), Nitrogen (N), Oxygen (O), and Sulfur (S) in proteins, and Hydrogen (H), Carbon (C), Nitrogen (N), Oxygen (O), Sulfur (S), Phosphorus (P), Fluorine (F), Chlorine (Cl), Bromine (Br), and Iodine (I) in ligands. These interactions result in a total of $50$ pairs of atom types for each protein-ligand complex \cite{cang2017topologynet}. However, due to the absence of H in most proteins, we reduce the number of atom pairs to $40$ in practice, ignoring the element $H$ in all proteins. These $40$ atom pairs, formed by atom types \{C, N, O, S\} in proteins, and atom types \{H, C, N, O, S, P, F, Cl, Br, I\} in ligands, along with their $xyz$ coordinates, are used to generate the topological features for each protein-ligand complex. In this paper, all atom-pair complexes are determined by a cutoff distance $12\mathring{\text{A}}$ from the ligand.

\begin{figure}[t]
	\centering
        \includegraphics[height=3.4cm]{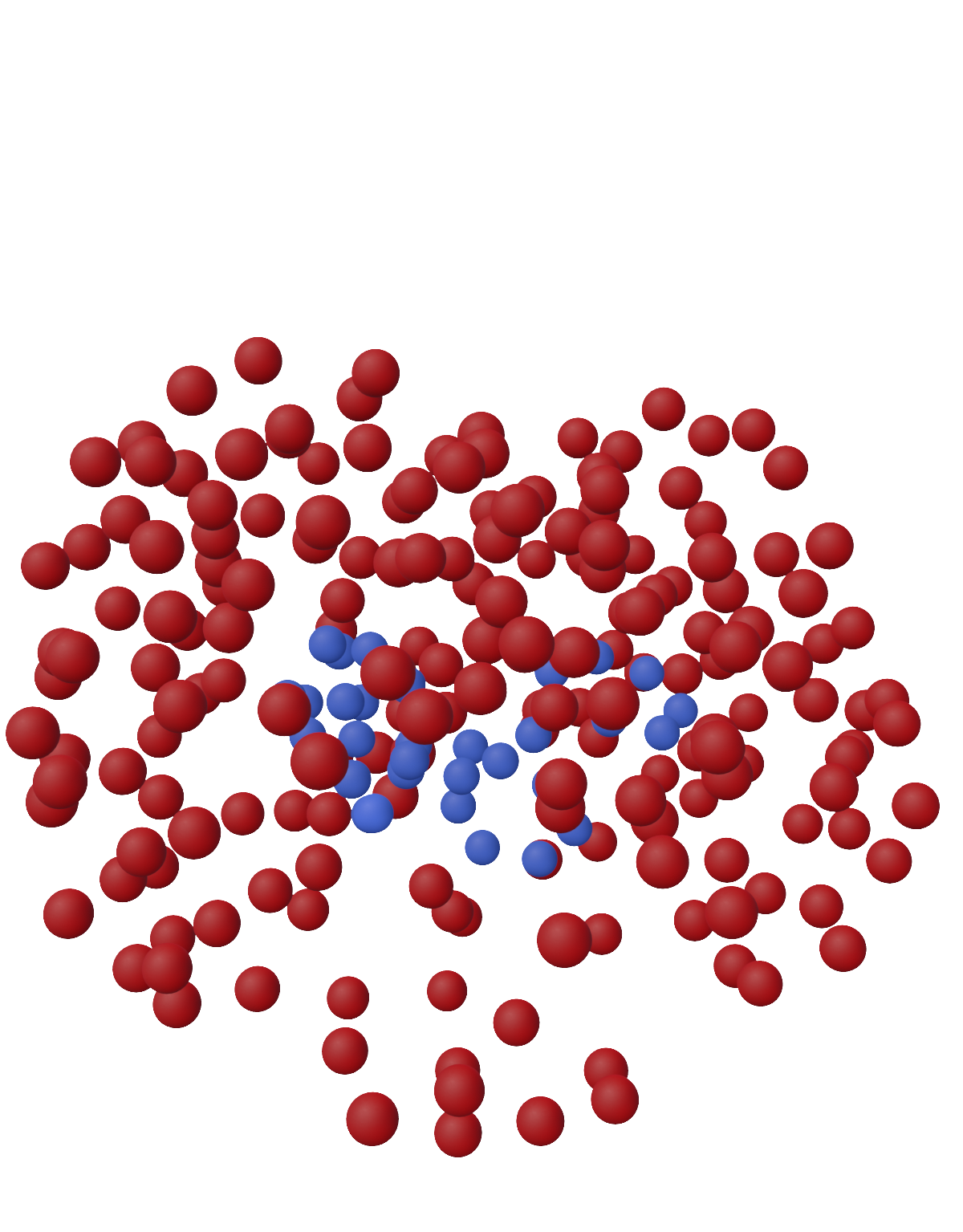}
        \includegraphics[height=3.4cm]{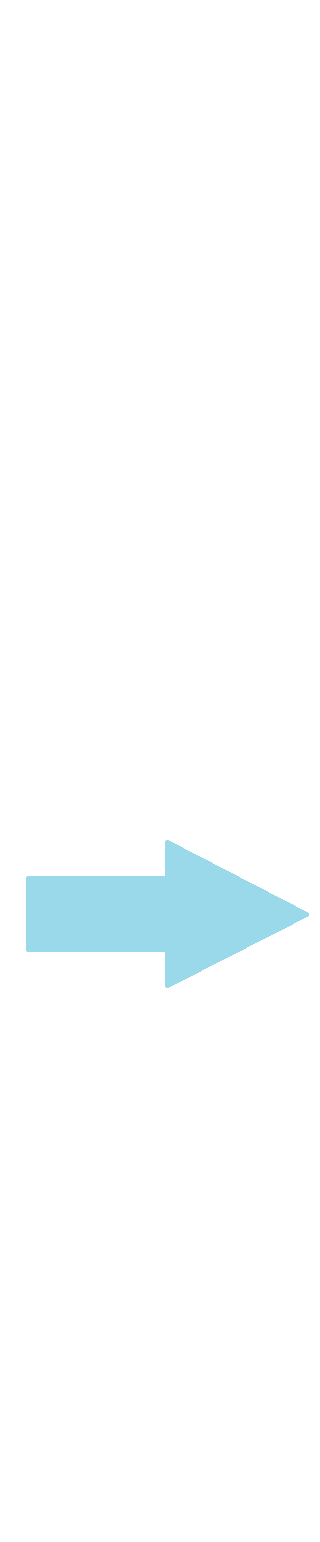}
	\includegraphics[height=3.4cm]{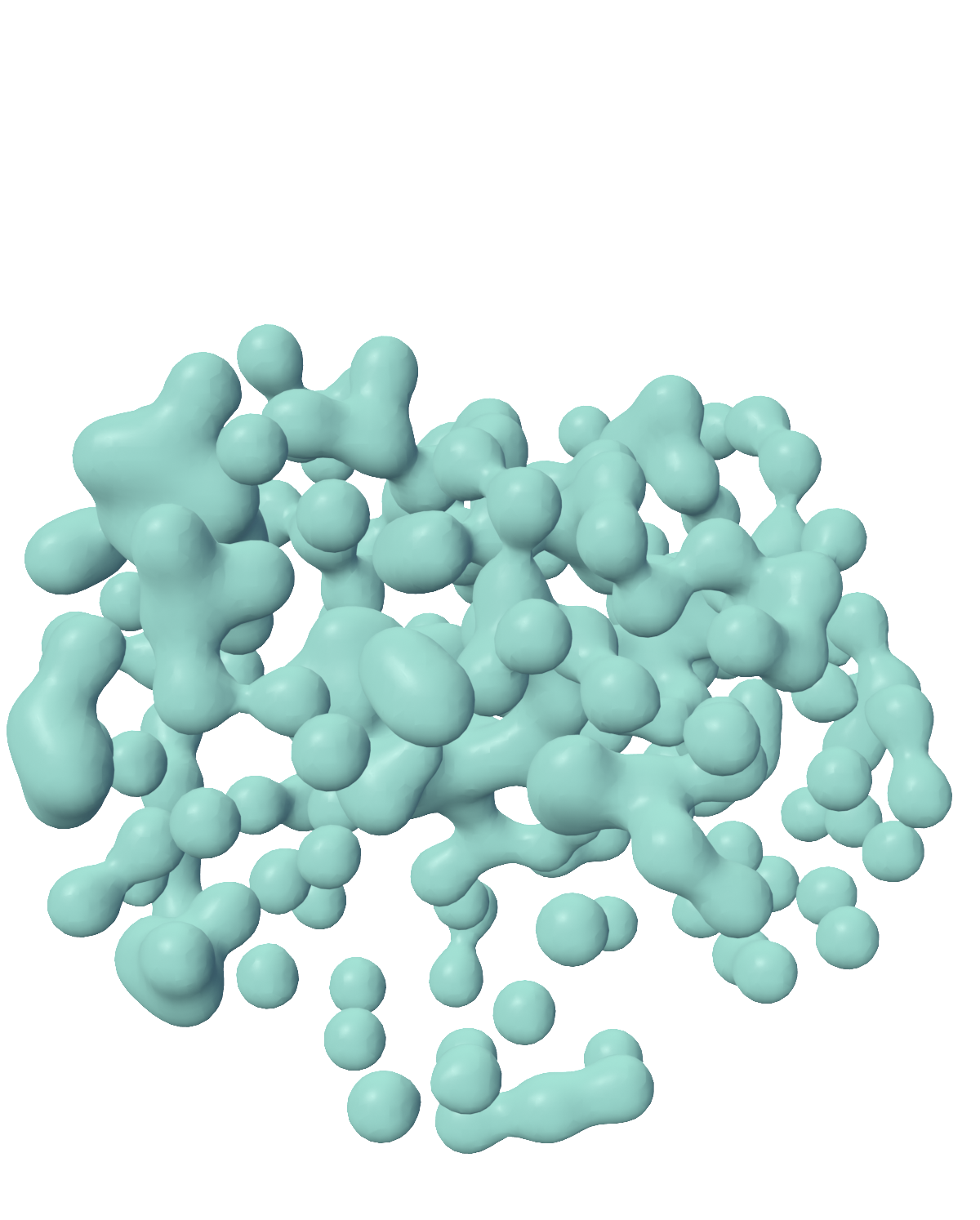}
	\includegraphics[height=3.4cm]{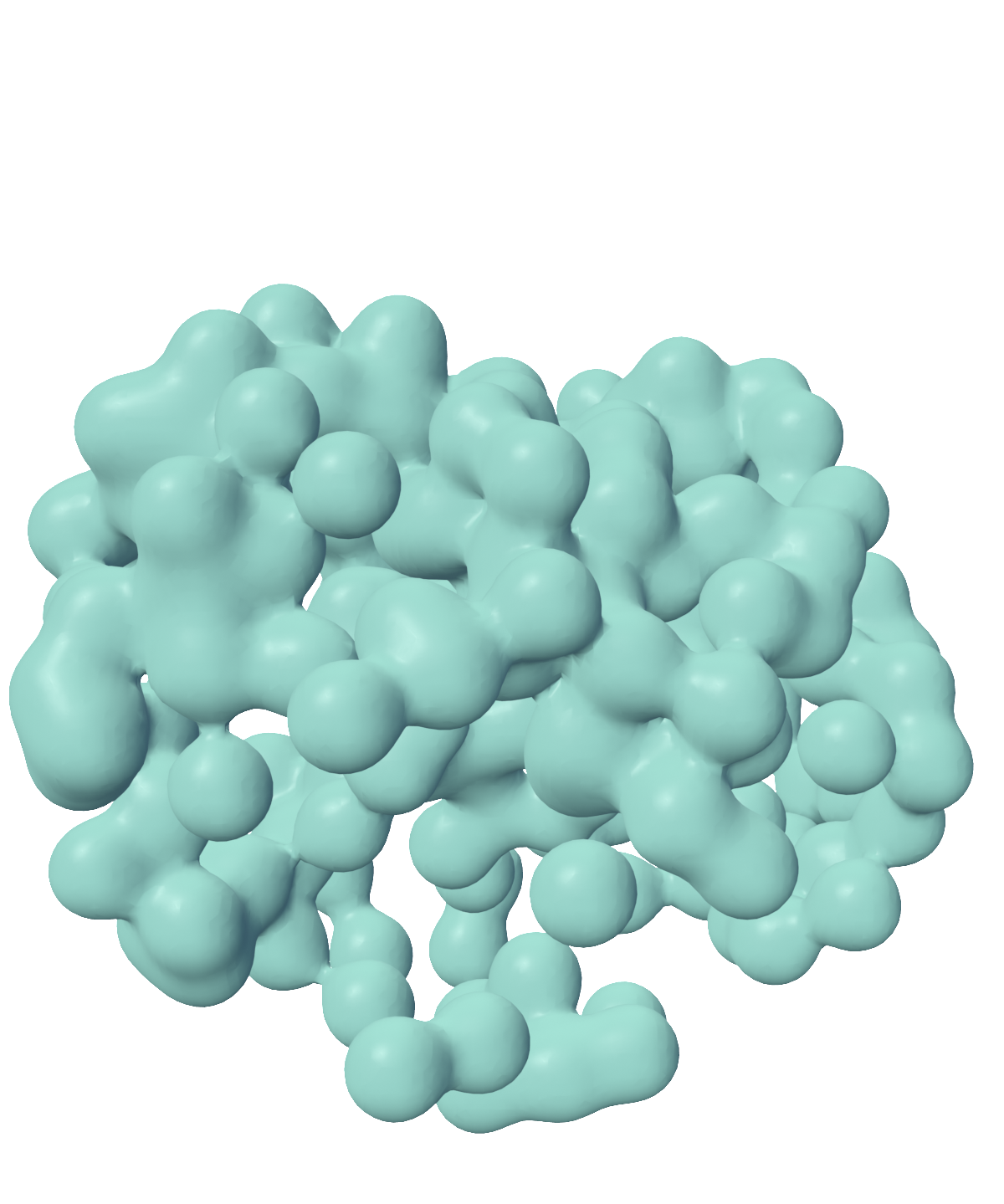}
	\includegraphics[height=3.4cm]{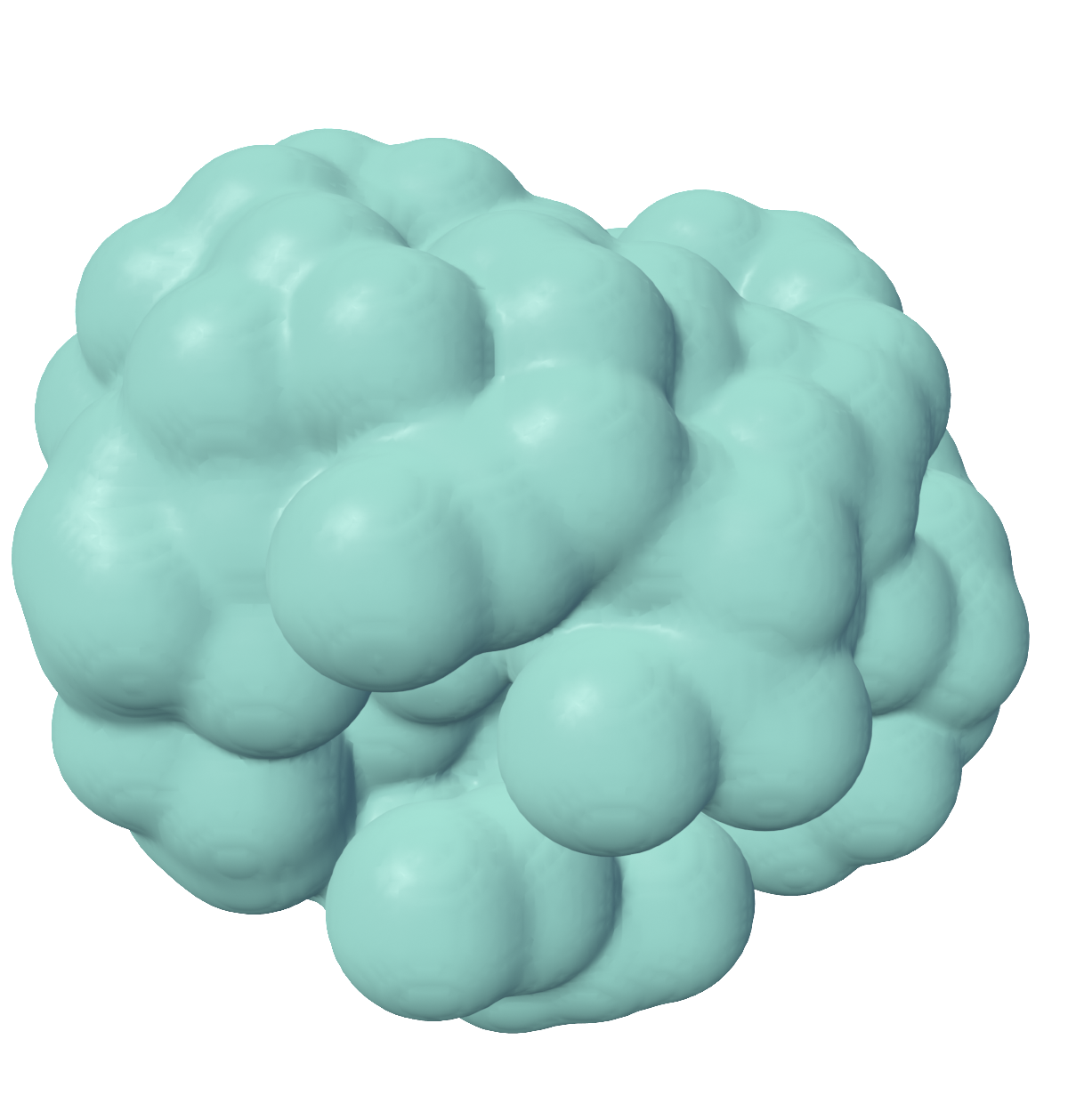}
	\caption{Left: atoms in the atom pair of type OH in protein-ligand complex 4mnt, with O shown in red and H in blue. Right: a filtration of manifold for this atom pair complex at $3$ different isovalues with level set function Eq.~\eqref{eq.lvf.rho}.}
	\label{fig.emfld.4tmn_OH}
\end{figure}


Let $\{ \mathbf{x}^{\alpha}_i, i=1,\cdots, s\}$ be the location coordinates of all $s$ atoms in an atom pair, where $\alpha$ denotes the atom type of the atom either in the protein or in the ligand. For this atom pair, a level set function can then be obtained by considering the negative sum of Gaussian density functions defined at the $xyz$ coordinates of all atoms, given as
\begin{align}\label{eq.lvf.rho}
	\rho(\mathbf{x}, \tau) = -\sum^{s}_{i = 1}\exp\left(-\left(\frac{||\mathbf{x} - \mathbf{x^{\alpha}_i}||}{\tau r^{\alpha}_i}\right)^2\right),
\end{align}
where $||\mathbf{x} - \mathbf{x^{\alpha}_i}||$ is the Euclidean distance from position $\mathbf{x}$ to the location $\mathbf{x^{\alpha}_i}$ of the $i$-th atom, $\tau$ is a scalar value, and $r^{\alpha}_i$ is the van der Waals radius of the $i$-th atom, determined by the atom type $\alpha$. Given an isovalue $c$, the sublevel set
\begin{align}
	M = \{\mathbf{x}\,|\, \rho(\mathbf{x}, \tau)\leq c \}
\end{align}
defines a compact manifold in $\RR^3$ with its boundary given by the isosurface $\partial M = \{\mathbf{x} \, | \, \rho(\mathbf{x}, \tau) = c \}$. A filtration of a manifold for the atom pair can then be obtained by choosing a list of evenly spaced isovalues of this level set function \eqref{eq.lvf.rho}. Let $c_1< c_2<\cdots < c_s$ be such isovalues. We have their corresponding sublevel sets given as follows
\begin{align}
	M_1\subset M_2\subset \cdots \subset M_s,
\end{align}
where $M_i$ is the compact manifold associated to isovalue $c_i$. In Fig.~\ref{fig.emfld.4tmn_OH} we present one example of the resulting filtration of manifolds at $3$ different isovalues for atom pair OH in protein-ligand complex 4tmn. Note that the function \eqref{eq.lvf.rho} is a special case of the flexibility rigidity index (FRI) density function \cite{nguyen2019dg}, which has been shown computationally stable in converting discrete point cloud representations to continuous embeddings, and been used for generating protein boundary surfaces \cite{chen2021evolutionary} and interactive manifolds \cite{nguyen2019dg}. Therefore, one can also make other reasonable choices of FRI density functions to generate the filtration of manifolds.

\subsection{Machine learning feature extraction}
\label{sec.ml.featureExtraction}

In the computation of the Laplacians, one can ideally choose a common Cartesian grid such that it contains all manifolds of interest for all protein-ligand complexes, which ensures that all Laplacians are computed consistently, making their spectra comparable for different complexes. However, as atoms are spread out in the space for different atom pairs, we need to use a sufficiently large grid with a fine resolution for accurate computation of Laplacians, which significantly increases the computational load. Instead, we consider, for each type of atom pairs, a fixed Cartesian grid, regardless of the types of protein-ligand complexes. This approach also ensures that the topological features are comparable for different protein-ligand complexes, as all spectra are computed in a same grid for all atom pairs of the same type. For simplicity, we choose a fixed grid spacing for all Cartesian grids across different atom pairs.

We consider $9$ evenly spaced isovalues in the interval $[-0.5, -0.001]$ for all level set functions, which provide $9$ compact manifolds for each atom pair. Note that the level set function \eqref{eq.lvf.rho} is always less than $0$ and approaches $0$ as the norm of $\mathbf{x}$ increases. This interval is chosen as isovalues greater than $-0.001$ result in no change on the $0$-th Betti number $\beta_0$ of manifolds for most atom pairs, and isovalues smaller than $-0.5$ leads to high computational cost, as finer grids are necessary to resolve those isosurfaces. To ensure that the computation of Laplacians is accurate and no topological information is missing due to numerical errors caused by low resolution, we require that at least $8$ grid cells of the Cartesian grid are contained in each connected component of a manifold. We compute, for each manifold, the BIG Laplacian $L_{3,n}$ under the normal boundary condition, for which the number of its $0$ eigenvalues gives the 0-th Betti number $\beta_0$. We then use the $0$-th Betti number $\beta_0$ and the first $k$ non-zero eigenvalues of $L_{3,n}$, as the topological feature for the manifold. These $k+1$ features for each of the $9$ compact manifolds for each atom pair, amount to $(k+1)\!\times\!9\!\times\!40$ topological features for each protein-ligand complex. While we only used $9$ isovalues within this interval for generating the manifolds in our experiments, more isovalues can be considered, which gives a filtration of more manifolds for each atom pair, and finally leads to more topological features for each protein-ligand complex.

The spectra of the $0$-th Laplacian, which in our case corresponds to $L_{3,n}$ under the normal boundary condition, have proven effective and successful in many machine learning tasks \cite{cang2018integration, liu2017forging, meng2021persistent, wang2020persistent}. While the Laplacians of other orders could also be used for generating more topological features, we utilize, in this preliminary test, only the spectra of $L_{3,n}$ as features for the protein-ligand complexes in the machine learning model due to the computation efficiency. The results, as shown in Sec.\ref{sec.ml.results}, indicate that these features are sufficient to validate our framework in the machine learning task for predicting the protein-ligand binding affinities. 

\subsection{Machine learning algorithm}

The machine learning models for predicting protein-ligand binding affinities often fall into two categories depending on the type of input data: complex-based or sequence-based models. The complex-based methods are trained using features obtained from the 3D protein-ligand complexes, while the sequence-based models learn from the one-dimensional protein sequences and the ligand simplified molecular-input line-entry system (SMILES) strings. In our experiments, besides the topological features from the 3D protein-ligand complexes, we incorporate protein-ligand features obtained from sequence-based models to build consensus models. To be specific, we make use of the recent pre-trained transformer protein language model Evolutionary Scale Modeling-2 (ESM-2) \cite{lin2022language}, and the pre-trained Transformer-CPZ model \cite{chen2021extracting} for generating the protein and ligand features, respectively, and use their concatenation as inputs for the binding affinity prediction. The residue embeddings from the last layer of the pre-trained ESM-2 model esm.pretrained.esm2\_t33\_650M\_UR50D are used as the protein features, while the embeddings from the last layer of the pre-trained Transformer-CPZ model chembl27\_pubchem\_zinc\_512 are used as the ligand features.

With the topological features and the embedding features obtained from ESM-2 and Transformer-CPZ, we employ the Gradient Boosting Regressor (GBR) module from Scikit-learn 1.4.2 for predicting the protein-ligand binding affinities. We then use the consensus prediction from these models as the final results. The GBR parameters used in our experiments are: n\_estimators=10,000, max\_depth=5, min\_samples\_split=5, learning\_rate=0.005, loss=squared\_error, subsample=0.5, and max\_features=sqrt. Changing these parameters does not result in significant differences. To address the randomness from the machine learning algorithm, we repeat each modeling process 20 times with different random seeds, and use the average of predictive results. The Pearson correlation coefficients (PCC) are used as the evaluation metric to assess the performance of our proposed models.

\subsection{Experimental Results}
\label{sec.ml.results}

\renewcommand{\arraystretch}{1.5}
\begin{table}
	\caption{Model performance on PDBbind-v2007 and PDBbind-v2016 benchmarks}
	\label{table.ml.performance}
	\begin{tabular}{lcccc}
		\hline\hline
		& \textbf{Method} & \textbf{PCC} & \textbf{RMSE (kcal/mol)}\\
		\hline
		\multirow{3}{*}{\textbf{PDBbind-v2007}} & PHL & 0.794 & 2.066 \\
		& TF & 0.795 &  2.006\\
		& Consensus & 0.826 & 1.954 \\
		\hline\hline
		\multirow{3}{*}{\textbf{PDBbind-v2016}} & PHL & 0.808 &  1.863\\
		& TF & 0.836 &  1.716\\
		& Consensus & 0.849 & 1.728\\
		\hline\hline
	
	\end{tabular}
\small
Abbreviations: PCC, Pearson correlation coefficient; RMSE, root mean squared error.

\end{table}

\begin{figure}[t]
	\centering
	\includegraphics[scale=0.34]{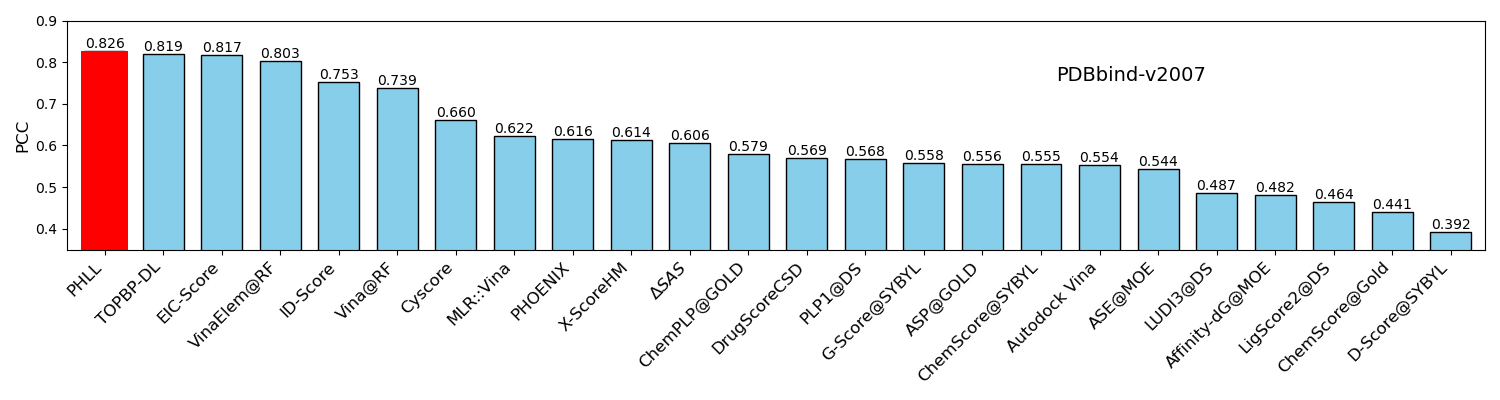}\\
	\includegraphics[scale=0.34]{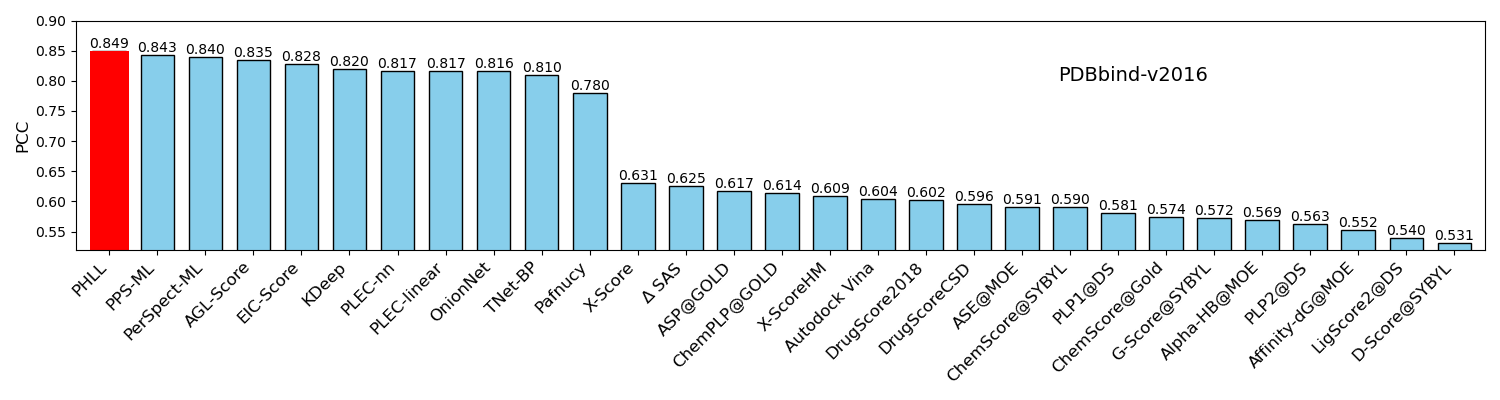}
	\caption{Performance comparison of the proposed model with other machine learning models for the two PDBbind datasets. The results of the proposed model (PHLL) are in red. The results of other methods are adapted from Refs. \cite{cang2018integration, cang2018representability, liu2017forging, su2018comparative, meng2021persistent, liu2023persistent}}
	\label{fig.pcc_comparison}
\end{figure}

The number of topological features for each protein-ligand complex, as in Sec.~\ref{sec.ml.featureExtraction}, is given by $(1+k)\!\times\!9\!\times\!40$, where $k$ denotes the number of the first $k$ non-zero eigenvalues of the Laplacians. To find the optimal parameter $k$ leading to the best performance of predictive modules, we carry out the five-fold cross-validation on the training set of each PDBbind dataset with varying values of $k$ based on the average of PCC values. 
The results indicate that the optimal PCC values for the PDBbind-v2007 and PDBbind-v2016 training sets can be achieved when $k=5$ and $k=10$, respectively. For the PDBbind-v2007 training set, the PCC value is $0.709$ and the RMSE value is $2.049$, while for the PDBbind-v2016 training set, the PCC value is $0.748$ and the RMSE value of $1.812$. 
These choices of $k$ result in a total of 2,160 topological features for each protein-ligand complex in the PDBbind-v2007 dataset and 3,960 topological features for each protein-ligand complex in the PDBbind-v2016 dataset. These topological features, along with the concatenated protein-ligand features from ESM-2 and Transformer-CPZ, are then used as inputs of the gradient-boosting regressor for binding affinity prediction. 

In Table~\ref{table.ml.performance}, we report the average PCC values and the average root mean squared error (RMSE) of our models on the test set for each PDBbind dataset using only the topological features from PHL, the model using only the transformer features (TF), and the consensus module using both types of features. With the incorporation of topological features, one can see a significant improvement in PCC values when using the proposed consensus model for each dataset, compared to the model using only TF features. The best performance is achieved when using the consensus model, yielding a PCC value of $0.826$ with RMSE given as $1.954$ for PDBbind-v2007 and $0.849$ with RMSE $1.728$ for PDBbind-v2016. In addition, we present the Pearson correlation coefficients obtained from our model and those in the previous studies, with results from \cite{cang2018integration, cang2018representability, su2018comparative, meng2021persistent, liu2023persistent}. As illustrated in Fig.~\ref{fig.pcc_comparison}, our model outperforms all the other models for the two PDBbind datasets. These results demonstrate the utility and effectiveness of our method in capturing the topological features.

\section{Conclusion}
Although there has had tremendous success of topological data analysis (TDA) \cite{nguyen2020review,nguyen2019mathematical}, particularly, topological deep learning (TDL) on point cloud data \cite{cang2017topologynet,papamarkou2024position}, there are few methods for the topological analysis of data on manifolds or manifold topological analysis \cite{desbrun2006discrete}. To fill this gap,  we presented a new method, persistent Hodge Laplacian (PHL) in the Eulerian representation,  for manifold topological learning (MTL) of real-world data on manifolds. PHL differs from existing state-of-the-art TDA methods on point clouds in the sense that the proposed  PHL is defined on manifolds, for which the traditional TDA methods do not work.  Additionally, 
PHL extends our earlier evolutionary de Rham-Hodge theory constructed on the Lagrangian representation \cite{chen2021evolutionary} to the Eulerian representation, which avoids numerical inconsistency over multiscale manifolds. 
We offer two discrete Hodge stars that mimic the continuous operator and developed both a continuous theory for mapping of normal forms across manifolds in a filtration to enable persistent cohomology analysis and the associated topology-persevering discrete construction on Cartesian grids. A proof-of-principle test on two benchmark datasets validates our  MTL model, highlighting its simplicity and promise for the predictions of data on manifolds. 

The popularity of TDA is facilitated by effective software packages, such as JavaPlex \cite{adams2014javaplex}, Perseus \cite{mischaikow2013morse}, Ripser \cite{bauer2021ripser}, etc. The further development of efficient PHL software is an important task. The computational efficiency has not been studied in this work. Algorithm acceleration and parallel and GPU architecture are to be explored.  Further experimental validations of manifold topological learning are also needed.

\section*{Acknowledgments}
This work was supported in part by NIH grants R01GM126189, R01AI164266, and R35GM148196, National Science Foundation grants DMS2052983, DMS-1761320, and IIS-1900473, NASA  grant 80NSSC21M0023,   Michigan State University Research Foundation, and  Bristol-Myers Squibb  65109. The authors thank Dr. Hongsong Feng for his help on transformer embeddings.

\bibliography{refs}


\begin{thebibliography}{62}
\ifx \bisbn   \undefined \def \bisbn  #1{ISBN #1}\fi
\ifx \binits  \undefined \def \binits#1{#1}\fi
\ifx \bauthor  \undefined \def \bauthor#1{#1}\fi
\ifx \batitle  \undefined \def \batitle#1{#1}\fi
\ifx \bjtitle  \undefined \def \bjtitle#1{#1}\fi
\ifx \bvolume  \undefined \def \bvolume#1{\textbf{#1}}\fi
\ifx \byear  \undefined \def \byear#1{#1}\fi
\ifx \bissue  \undefined \def \bissue#1{#1}\fi
\ifx \bfpage  \undefined \def \bfpage#1{#1}\fi
\ifx \blpage  \undefined \def \blpage #1{#1}\fi
\ifx \burl  \undefined \def \burl#1{\textsf{#1}}\fi
\ifx \doiurl  \undefined \def \doiurl#1{\url{https://doi.org/#1}}\fi
\ifx \betal  \undefined \def \betal{\textit{et al.}}\fi
\ifx \binstitute  \undefined \def \binstitute#1{#1}\fi
\ifx \binstitutionaled  \undefined \def \binstitutionaled#1{#1}\fi
\ifx \bctitle  \undefined \def \bctitle#1{#1}\fi
\ifx \beditor  \undefined \def \beditor#1{#1}\fi
\ifx \bpublisher  \undefined \def \bpublisher#1{#1}\fi
\ifx \bbtitle  \undefined \def \bbtitle#1{#1}\fi
\ifx \bedition  \undefined \def \bedition#1{#1}\fi
\ifx \bseriesno  \undefined \def \bseriesno#1{#1}\fi
\ifx \blocation  \undefined \def \blocation#1{#1}\fi
\ifx \bsertitle  \undefined \def \bsertitle#1{#1}\fi
\ifx \bsnm \undefined \def \bsnm#1{#1}\fi
\ifx \bsuffix \undefined \def \bsuffix#1{#1}\fi
\ifx \bparticle \undefined \def \bparticle#1{#1}\fi
\ifx \barticle \undefined \def \barticle#1{#1}\fi
\bibcommenthead
\ifx \bconfdate \undefined \def \bconfdate #1{#1}\fi
\ifx \botherref \undefined \def \botherref #1{#1}\fi
\ifx \url \undefined \def \url#1{\textsf{#1}}\fi
\ifx \bchapter \undefined \def \bchapter#1{#1}\fi
\ifx \bbook \undefined \def \bbook#1{#1}\fi
\ifx \bcomment \undefined \def \bcomment#1{#1}\fi
\ifx \oauthor \undefined \def \oauthor#1{#1}\fi
\ifx \citeauthoryear \undefined \def \citeauthoryear#1{#1}\fi
\ifx \endbibitem  \undefined \def \endbibitem {}\fi
\ifx \bconflocation  \undefined \def \bconflocation#1{#1}\fi
\ifx \arxivurl  \undefined \def \arxivurl#1{\textsf{#1}}\fi
\csname PreBibitemsHook\endcsname

\bibitem[\protect\citeauthoryear{Wasserman}{2018}]{wasserman2018topological}
\begin{barticle}
\bauthor{\bsnm{Wasserman}, \binits{L.}}:
\batitle{Topological data analysis}.
\bjtitle{Annual Review of Statistics and Its Application}
\bvolume{5}(\bissue{1}),
\bfpage{501}--\blpage{532}
(\byear{2018})
\end{barticle}
\endbibitem

\bibitem[\protect\citeauthoryear{Mischaikow and
  Nanda}{2013}]{mischaikow2013morse}
\begin{barticle}
\bauthor{\bsnm{Mischaikow}, \binits{K.}},
\bauthor{\bsnm{Nanda}, \binits{V.}}:
\batitle{Morse theory for filtrations and efficient computation of persistent
  homology}.
\bjtitle{Discrete \& Computational Geometry}
\bvolume{50}(\bissue{2}),
\bfpage{330}--\blpage{353}
(\byear{2013})
\end{barticle}
\endbibitem

\bibitem[\protect\citeauthoryear{Carlsson}{2009}]{carlsson2009topology}
\begin{barticle}
\bauthor{\bsnm{Carlsson}, \binits{G.}}:
\batitle{Topology and data}.
\bjtitle{Bulletin of the American Mathematical Society}
\bvolume{46}(\bissue{2}),
\bfpage{255}--\blpage{308}
(\byear{2009})
\end{barticle}
\endbibitem

\bibitem[\protect\citeauthoryear{Edelsbrunner
  et~al.}{2008}]{edelsbrunner2008persistent}
\begin{barticle}
\bauthor{\bsnm{Edelsbrunner}, \binits{H.}},
\bauthor{\bsnm{Harer}, \binits{J.}}, \betal:
\batitle{Persistent homology-a survey}.
\bjtitle{Contemporary mathematics}
\bvolume{453}(\bissue{26}),
\bfpage{257}--\blpage{282}
(\byear{2008})
\end{barticle}
\endbibitem

\bibitem[\protect\citeauthoryear{Zomorodian and
  Carlsson}{2005}]{zomorodian2005computing}
\begin{barticle}
\bauthor{\bsnm{Zomorodian}, \binits{A.}},
\bauthor{\bsnm{Carlsson}, \binits{G.}}:
\batitle{Computing persistent homology}.
\bjtitle{Discrete \& Computational Geometry}
\bvolume{33}(\bissue{2}),
\bfpage{249}--\blpage{274}
(\byear{2005})
\end{barticle}
\endbibitem

\bibitem[\protect\citeauthoryear{Ghrist}{2008}]{ghrist2008barcodes}
\begin{barticle}
\bauthor{\bsnm{Ghrist}, \binits{R.}}:
\batitle{Barcodes: the persistent topology of data}.
\bjtitle{Bulletin of the American Mathematical Society}
\bvolume{45}(\bissue{1}),
\bfpage{61}--\blpage{75}
(\byear{2008})
\end{barticle}
\endbibitem

\bibitem[\protect\citeauthoryear{Bubenik et~al.}{2015}]{bubenik2015statistical}
\begin{barticle}
\bauthor{\bsnm{Bubenik}, \binits{P.}}, \betal:
\batitle{Statistical topological data analysis using persistence landscapes.}
\bjtitle{J. Mach. Learn. Res.}
\bvolume{16}(\bissue{1}),
\bfpage{77}--\blpage{102}
(\byear{2015})
\end{barticle}
\endbibitem

\bibitem[\protect\citeauthoryear{Dey et~al.}{2014}]{dey2014computing}
\begin{bchapter}
\bauthor{\bsnm{Dey}, \binits{T.K.}},
\bauthor{\bsnm{Fan}, \binits{F.}},
\bauthor{\bsnm{Wang}, \binits{Y.}}:
\bctitle{Computing topological persistence for simplicial maps}.
In: \bbtitle{Proceedings of the Thirtieth Annual Symposium on Computational
  Geometry},
pp. \bfpage{345}--\blpage{354}
(\byear{2014})
\end{bchapter}
\endbibitem

\bibitem[\protect\citeauthoryear{Adams et~al.}{2017}]{adams2017persistence}
\begin{barticle}
\bauthor{\bsnm{Adams}, \binits{H.}},
\bauthor{\bsnm{Emerson}, \binits{T.}},
\bauthor{\bsnm{Kirby}, \binits{M.}},
\bauthor{\bsnm{Neville}, \binits{R.}},
\bauthor{\bsnm{Peterson}, \binits{C.}},
\bauthor{\bsnm{Shipman}, \binits{P.}},
\bauthor{\bsnm{Chepushtanova}, \binits{S.}},
\bauthor{\bsnm{Hanson}, \binits{E.}},
\bauthor{\bsnm{Motta}, \binits{F.}},
\bauthor{\bsnm{Ziegelmeier}, \binits{L.}}:
\batitle{Persistence images: A stable vector representation of persistent
  homology}.
\bjtitle{Journal of Machine Learning Research}
\bvolume{18}(\bissue{8}),
\bfpage{1}--\blpage{35}
(\byear{2017})
\end{barticle}
\endbibitem

\bibitem[\protect\citeauthoryear{Xia and Wei}{2014}]{xia2014persistent}
\begin{barticle}
\bauthor{\bsnm{Xia}, \binits{K.}},
\bauthor{\bsnm{Wei}, \binits{G.-W.}}:
\batitle{Persistent homology analysis of protein structure, flexibility, and
  folding}.
\bjtitle{International journal for numerical methods in biomedical engineering}
\bvolume{30}(\bissue{8}),
\bfpage{814}--\blpage{844}
(\byear{2014})
\end{barticle}
\endbibitem

\bibitem[\protect\citeauthoryear{Townsend
  et~al.}{2020}]{townsend2020representation}
\begin{barticle}
\bauthor{\bsnm{Townsend}, \binits{J.}},
\bauthor{\bsnm{Micucci}, \binits{C.P.}},
\bauthor{\bsnm{Hymel}, \binits{J.H.}},
\bauthor{\bsnm{Maroulas}, \binits{V.}},
\bauthor{\bsnm{Vogiatzis}, \binits{K.D.}}:
\batitle{Representation of molecular structures with persistent homology for
  machine learning applications in chemistry}.
\bjtitle{Nature communications}
\bvolume{11}(\bissue{1}),
\bfpage{3230}
(\byear{2020})
\end{barticle}
\endbibitem

\bibitem[\protect\citeauthoryear{MacPherson and
  Schweinhart}{2012}]{macpherson2012measuring}
\begin{botherref}
\oauthor{\bsnm{MacPherson}, \binits{R.}},
\oauthor{\bsnm{Schweinhart}, \binits{B.}}:
Measuring shape with topology.
Journal of Mathematical Physics
\textbf{53}(7)
(2012)
\end{botherref}
\endbibitem

\bibitem[\protect\citeauthoryear{Cang et~al.}{2015}]{cang2015topological}
\begin{botherref}
\oauthor{\bsnm{Cang}, \binits{Z.}},
\oauthor{\bsnm{Mu}, \binits{L.}},
\oauthor{\bsnm{Wu}, \binits{K.}},
\oauthor{\bsnm{Opron}, \binits{K.}},
\oauthor{\bsnm{Xia}, \binits{K.}},
\oauthor{\bsnm{Wei}, \binits{G.-W.}}:
A topological approach for protein classification.
Computational and Mathematical Biophysics
\textbf{3}(1)
(2015)
\end{botherref}
\endbibitem

\bibitem[\protect\citeauthoryear{Cang and Wei}{2017}]{cang2017topologynet}
\begin{barticle}
\bauthor{\bsnm{Cang}, \binits{Z.}},
\bauthor{\bsnm{Wei}, \binits{G.-W.}}:
\batitle{Topologynet: Topology based deep convolutional and multi-task neural
  networks for biomolecular property predictions}.
\bjtitle{PLoS computational biology}
\bvolume{13}(\bissue{7}),
\bfpage{1005690}
(\byear{2017})
\end{barticle}
\endbibitem

\bibitem[\protect\citeauthoryear{Papamarkou
  et~al.}{2024}]{papamarkou2024position}
\begin{bchapter}
\bauthor{\bsnm{Papamarkou}, \binits{T.}},
\bauthor{\bsnm{Birdal}, \binits{T.}},
\bauthor{\bsnm{Bronstein}, \binits{M.M.}},
\bauthor{\bsnm{Carlsson}, \binits{G.E.}},
\bauthor{\bsnm{Curry}, \binits{J.}},
\bauthor{\bsnm{Gao}, \binits{Y.}},
\bauthor{\bsnm{Hajij}, \binits{M.}},
\bauthor{\bsnm{Kwitt}, \binits{R.}},
\bauthor{\bsnm{Lio}, \binits{P.}},
\bauthor{\bsnm{Di~Lorenzo}, \binits{P.}}, \betal:
\bctitle{Position: Topological deep learning is the new frontier for relational
  learning}.
In: \bbtitle{Forty-first International Conference on Machine Learning}
(\byear{2024})
\end{bchapter}
\endbibitem

\bibitem[\protect\citeauthoryear{Nguyen et~al.}{2020a}]{nguyen2020review}
\begin{barticle}
\bauthor{\bsnm{Nguyen}, \binits{D.D.}},
\bauthor{\bsnm{Cang}, \binits{Z.}},
\bauthor{\bsnm{Wei}, \binits{G.-W.}}:
\batitle{A review of mathematical representations of biomolecular data}.
\bjtitle{Physical Chemistry Chemical Physics}
\bvolume{22}(\bissue{8}),
\bfpage{4343}--\blpage{4367}
(\byear{2020})
\end{barticle}
\endbibitem

\bibitem[\protect\citeauthoryear{Nguyen et~al.}{2020b}]{nguyen2020mathdl}
\begin{barticle}
\bauthor{\bsnm{Nguyen}, \binits{D.D.}},
\bauthor{\bsnm{Gao}, \binits{K.}},
\bauthor{\bsnm{Wang}, \binits{M.}},
\bauthor{\bsnm{Wei}, \binits{G.-W.}}:
\batitle{Mathdl: mathematical deep learning for d3r grand challenge 4}.
\bjtitle{Journal of computer-aided molecular design}
\bvolume{34}(\bissue{2}),
\bfpage{131}--\blpage{147}
(\byear{2020})
\end{barticle}
\endbibitem

\bibitem[\protect\citeauthoryear{Nguyen et~al.}{2019}]{nguyen2019mathematical}
\begin{barticle}
\bauthor{\bsnm{Nguyen}, \binits{D.D.}},
\bauthor{\bsnm{Cang}, \binits{Z.}},
\bauthor{\bsnm{Wu}, \binits{K.}},
\bauthor{\bsnm{Wang}, \binits{M.}},
\bauthor{\bsnm{Cao}, \binits{Y.}},
\bauthor{\bsnm{Wei}, \binits{G.-W.}}:
\batitle{Mathematical deep learning for pose and binding affinity prediction
  and ranking in d3r grand challenges}.
\bjtitle{Journal of computer-aided molecular design}
\bvolume{33},
\bfpage{71}--\blpage{82}
(\byear{2019})
\end{barticle}
\endbibitem

\bibitem[\protect\citeauthoryear{Chen et~al.}{2020}]{chen2020mutations}
\begin{barticle}
\bauthor{\bsnm{Chen}, \binits{J.}},
\bauthor{\bsnm{Wang}, \binits{R.}},
\bauthor{\bsnm{Wang}, \binits{M.}},
\bauthor{\bsnm{Wei}, \binits{G.-W.}}:
\batitle{Mutations strengthened sars-cov-2 infectivity}.
\bjtitle{Journal of molecular biology}
\bvolume{432}(\bissue{19}),
\bfpage{5212}--\blpage{5226}
(\byear{2020})
\end{barticle}
\endbibitem

\bibitem[\protect\citeauthoryear{Wang et~al.}{2021}]{wang2021mechanisms}
\begin{barticle}
\bauthor{\bsnm{Wang}, \binits{R.}},
\bauthor{\bsnm{Chen}, \binits{J.}},
\bauthor{\bsnm{Wei}, \binits{G.-W.}}:
\batitle{Mechanisms of sars-cov-2 evolution revealing vaccine-resistant
  mutations in europe and america}.
\bjtitle{The journal of physical chemistry letters}
\bvolume{12}(\bissue{49}),
\bfpage{11850}--\blpage{11857}
(\byear{2021})
\end{barticle}
\endbibitem

\bibitem[\protect\citeauthoryear{Chen and Wei}{2022}]{chen2022omicronBA2}
\begin{barticle}
\bauthor{\bsnm{Chen}, \binits{J.}},
\bauthor{\bsnm{Wei}, \binits{G.-W.}}:
\batitle{Omicron ba. 2 (b. 1.1. 529.2): high potential for becoming the next
  dominant variant}.
\bjtitle{The journal of physical chemistry letters}
\bvolume{13}(\bissue{17}),
\bfpage{3840}--\blpage{3849}
(\byear{2022})
\end{barticle}
\endbibitem

\bibitem[\protect\citeauthoryear{Chen et~al.}{2022}]{chen2022persistent}
\begin{barticle}
\bauthor{\bsnm{Chen}, \binits{J.}},
\bauthor{\bsnm{Qiu}, \binits{Y.}},
\bauthor{\bsnm{Wang}, \binits{R.}},
\bauthor{\bsnm{Wei}, \binits{G.-W.}}:
\batitle{Persistent laplacian projected omicron ba. 4 and ba. 5 to become new
  dominating variants}.
\bjtitle{Computers in Biology and Medicine}
\bvolume{151},
\bfpage{106262}
(\byear{2022})
\end{barticle}
\endbibitem

\bibitem[\protect\citeauthoryear{Pun et~al.}{2018}]{pun2018persistent}
\begin{botherref}
\oauthor{\bsnm{Pun}, \binits{C.S.}},
\oauthor{\bsnm{Xia}, \binits{K.}},
\oauthor{\bsnm{Lee}, \binits{S.X.}}:
Persistent-homology-based machine learning and its applications--a survey.
arXiv preprint arXiv:1811.00252
(2018)
\end{botherref}
\endbibitem

\bibitem[\protect\citeauthoryear{Wei and Wei}{2023}]{wei2023persistentb}
\begin{botherref}
\oauthor{\bsnm{Wei}, \binits{X.}},
\oauthor{\bsnm{Wei}, \binits{G.-W.}}:
Persistent topological laplacians--a survey.
arXiv preprint arXiv:2312.07563
(2023)
\end{botherref}
\endbibitem

\bibitem[\protect\citeauthoryear{Wang et~al.}{2020a}]{wang2020persistent}
\begin{barticle}
\bauthor{\bsnm{Wang}, \binits{R.}},
\bauthor{\bsnm{Nguyen}, \binits{D.D.}},
\bauthor{\bsnm{Wei}, \binits{G.-W.}}:
\batitle{Persistent spectral graph}.
\bjtitle{International Journal for Numerical Methods in Biomedical Engineering}
\bvolume{36}(\bissue{9}),
\bfpage{3376}
(\byear{2020})
\end{barticle}
\endbibitem

\bibitem[\protect\citeauthoryear{Wang et~al.}{2020b}]{wang2021hermes}
\begin{barticle}
\bauthor{\bsnm{Wang}, \binits{R.}},
\bauthor{\bsnm{Zhao}, \binits{R.}},
\bauthor{\bsnm{Ribando-Gros}, \binits{E.}},
\bauthor{\bsnm{Chen}, \binits{J.}},
\bauthor{\bsnm{Tong}, \binits{Y.}},
\bauthor{\bsnm{Wei}, \binits{G.-W.}}:
\batitle{Hermes: Persistent spectral graph software}.
\bjtitle{Foundations of Data Science}
\bvolume{3}(\bissue{1}),
\bfpage{67}--\blpage{97}
(\byear{2020})
\end{barticle}
\endbibitem

\bibitem[\protect\citeauthoryear{Dong}{2024}]{dong2024faster}
\begin{botherref}
\oauthor{\bsnm{Dong}, \binits{R.}}:
A faster algorithm of up persistent laplacian over non-branching simplicial
  complexes.
arXiv preprint arXiv:2408.16741
(2024)
\end{botherref}
\endbibitem

\bibitem[\protect\citeauthoryear{Liu et~al.}{2023}]{liu2023algebraic}
\begin{botherref}
\oauthor{\bsnm{Liu}, \binits{J.}},
\oauthor{\bsnm{Li}, \binits{J.}},
\oauthor{\bsnm{Wu}, \binits{J.}}:
The algebraic stability for persistent laplacians.
arXiv preprint arXiv:2302.03902
(2023)
\end{botherref}
\endbibitem

\bibitem[\protect\citeauthoryear{M{\'e}moli
  et~al.}{2022}]{memoli2022persistent}
\begin{barticle}
\bauthor{\bsnm{M{\'e}moli}, \binits{F.}},
\bauthor{\bsnm{Wan}, \binits{Z.}},
\bauthor{\bsnm{Wang}, \binits{Y.}}:
\batitle{Persistent laplacians: Properties, algorithms and implications}.
\bjtitle{SIAM Journal on Mathematics of Data Science}
\bvolume{4}(\bissue{2}),
\bfpage{858}--\blpage{884}
(\byear{2022})
\end{barticle}
\endbibitem

\bibitem[\protect\citeauthoryear{G{\"u}len
  et~al.}{2023}]{gulen2023generalization}
\begin{botherref}
\oauthor{\bsnm{G{\"u}len}, \binits{A.B.}},
\oauthor{\bsnm{M{\'e}moli}, \binits{F.}},
\oauthor{\bsnm{Wan}, \binits{Z.}},
\oauthor{\bsnm{Wang}, \binits{Y.}}:
A generalization of the persistent laplacian to simplicial maps.
arXiv preprint arXiv:2302.03771
(2023)
\end{botherref}
\endbibitem

\bibitem[\protect\citeauthoryear{Wei and Wei}{2024}]{wei2024persistent}
\begin{botherref}
\oauthor{\bsnm{Wei}, \binits{X.}},
\oauthor{\bsnm{Wei}, \binits{G.-W.}}:
Persistent sheaf laplacian.
Foundations of data science (Springfield, Mo.),
10--39342024033
(2024)
\end{botherref}
\endbibitem

\bibitem[\protect\citeauthoryear{Liu et~al.}{2021}]{liu2021persistent}
\begin{barticle}
\bauthor{\bsnm{Liu}, \binits{X.}},
\bauthor{\bsnm{Feng}, \binits{H.}},
\bauthor{\bsnm{Wu}, \binits{J.}},
\bauthor{\bsnm{Xia}, \binits{K.}}:
\batitle{Persistent spectral hypergraph based machine learning (psh-ml) for
  protein-ligand binding affinity prediction}.
\bjtitle{Briefings in Bioinformatics}
\bvolume{22}(\bissue{5}),
\bfpage{127}
(\byear{2021})
\end{barticle}
\endbibitem

\bibitem[\protect\citeauthoryear{Meng and Xia}{2021}]{meng2021persistent}
\begin{barticle}
\bauthor{\bsnm{Meng}, \binits{Z.}},
\bauthor{\bsnm{Xia}, \binits{K.}}:
\batitle{Persistent spectral--based machine learning (perspect ml) for
  protein-ligand binding affinity prediction}.
\bjtitle{Science advances}
\bvolume{7}(\bissue{19}),
\bfpage{5329}
(\byear{2021})
\end{barticle}
\endbibitem

\bibitem[\protect\citeauthoryear{Chen et~al.}{2021}]{chen2021evolutionary}
\begin{barticle}
\bauthor{\bsnm{Chen}, \binits{J.}},
\bauthor{\bsnm{Zhao}, \binits{R.}},
\bauthor{\bsnm{Tong}, \binits{Y.}},
\bauthor{\bsnm{Wei}, \binits{G.-W.}}:
\batitle{Evolutionary de rham-hodge method}.
\bjtitle{Discrete and continuous dynamical systems. Series B}
\bvolume{26}(\bissue{7}),
\bfpage{3785}
(\byear{2021})
\end{barticle}
\endbibitem

\bibitem[\protect\citeauthoryear{Yang and Parr}{1984}]{yang1984electron}
\begin{barticle}
\bauthor{\bsnm{Yang}, \binits{W.}},
\bauthor{\bsnm{Parr}, \binits{R.G.}}:
\batitle{Electron density, kohn--sham frontier orbitals, and fukui functions}.
\bjtitle{The Journal of Chemical Physics}
\bvolume{81}(\bissue{6}),
\bfpage{2862}--\blpage{2863}
(\byear{1984})
\end{barticle}
\endbibitem

\bibitem[\protect\citeauthoryear{Chen et~al.}{2017}]{chen2017low}
\begin{barticle}
\bauthor{\bsnm{Chen}, \binits{H.}},
\bauthor{\bsnm{Zhang}, \binits{Y.}},
\bauthor{\bsnm{Zhang}, \binits{W.}},
\bauthor{\bsnm{Liao}, \binits{P.}},
\bauthor{\bsnm{Li}, \binits{K.}},
\bauthor{\bsnm{Zhou}, \binits{J.}},
\bauthor{\bsnm{Wang}, \binits{G.}}:
\batitle{Low-dose ct via convolutional neural network}.
\bjtitle{Biomedical optics express}
\bvolume{8}(\bissue{2}),
\bfpage{679}--\blpage{694}
(\byear{2017})
\end{barticle}
\endbibitem

\bibitem[\protect\citeauthoryear{Khovanov}{2000}]{khovanov2000categorification}
\begin{barticle}
\bauthor{\bsnm{Khovanov}, \binits{M.}}:
\batitle{A categorification of the jones polynomial}.
\bjtitle{Duke Mathematical Journal}
\bvolume{101}(\bissue{3}),
\bfpage{359}--\blpage{426}
(\byear{2000})
\end{barticle}
\endbibitem

\bibitem[\protect\citeauthoryear{Panagiotou
  et~al.}{2019}]{panagiotou2019topological}
\begin{barticle}
\bauthor{\bsnm{Panagiotou}, \binits{E.}},
\bauthor{\bsnm{Millett}, \binits{K.C.}},
\bauthor{\bsnm{Atzberger}, \binits{P.J.}}:
\batitle{Topological methods for polymeric materials: characterizing the
  relationship between polymer entanglement and viscoelasticity}.
\bjtitle{Polymers}
\bvolume{11}(\bissue{3}),
\bfpage{437}
(\byear{2019})
\end{barticle}
\endbibitem

\bibitem[\protect\citeauthoryear{Shen et~al.}{accepted, 2024}]{shen2024knot}
\begin{botherref}
\oauthor{\bsnm{Shen}, \binits{L.}},
\oauthor{\bsnm{Feng}, \binits{H.}},
\oauthor{\bsnm{Li}, \binits{F.}},
\oauthor{\bsnm{Lei}, \binits{F.}},
\oauthor{\bsnm{Wu}, \binits{J.}},
\oauthor{\bsnm{Wei}, \binits{G.-W.}}:
Knot data analysis using multiscale gauss link integral.
Proceedings of the National Academy of Sciences
(accepted, 2024)
\end{botherref}
\endbibitem

\bibitem[\protect\citeauthoryear{Shen et~al.}{accepted
  2024}]{shen2024evolutionary}
\begin{botherref}
\oauthor{\bsnm{Shen}, \binits{L.}},
\oauthor{\bsnm{Liu}, \binits{J.}},
\oauthor{\bsnm{Wei}, \binits{G.-W.}}:
Evolutionary khovanov homology.
AIMS Mathematics
(accepted 2024)
\end{botherref}
\endbibitem

\bibitem[\protect\citeauthoryear{Ribando-Gros et~al.}{2024}]{ribando2024graph}
\begin{barticle}
\bauthor{\bsnm{Ribando-Gros}, \binits{E.}},
\bauthor{\bsnm{Wang}, \binits{R.}},
\bauthor{\bsnm{Chen}, \binits{J.}},
\bauthor{\bsnm{Tong}, \binits{Y.}},
\bauthor{\bsnm{Wei}, \binits{G.-W.}}:
\batitle{Combinatorial and hodge laplacians: Similarity and difference}.
\bjtitle{SIAM Review}
\bvolume{66}(\bissue{3}),
\bfpage{575}--\blpage{601}
(\byear{2024})
\end{barticle}
\endbibitem

\bibitem[\protect\citeauthoryear{Desbrun et~al.}{2006}]{desbrun2006discrete}
\begin{bchapter}
\bauthor{\bsnm{Desbrun}, \binits{M.}},
\bauthor{\bsnm{Kanso}, \binits{E.}},
\bauthor{\bsnm{Tong}, \binits{Y.}}:
\bctitle{Discrete differential forms for computational modeling}.
In: \bbtitle{ACM SIGGRAPH 2006 Courses},
pp. \bfpage{39}--\blpage{54}
(\byear{2006})
\end{bchapter}
\endbibitem

\bibitem[\protect\citeauthoryear{Dodziuk}{1976}]{dodziuk1976finite}
\begin{barticle}
\bauthor{\bsnm{Dodziuk}, \binits{J.}}:
\batitle{Finite-difference approach to the hodge theory of harmonic forms}.
\bjtitle{American Journal of Mathematics}
\bvolume{98}(\bissue{1}),
\bfpage{79}--\blpage{104}
(\byear{1976})
\end{barticle}
\endbibitem

\bibitem[\protect\citeauthoryear{Arnold et~al.}{2006}]{arnold2006finite}
\begin{barticle}
\bauthor{\bsnm{Arnold}, \binits{D.N.}},
\bauthor{\bsnm{Falk}, \binits{R.S.}},
\bauthor{\bsnm{Winther}, \binits{R.}}:
\batitle{Finite element exterior calculus, homological techniques, and
  applications}.
\bjtitle{Acta numerica}
\bvolume{15},
\bfpage{1}--\blpage{155}
(\byear{2006})
\end{barticle}
\endbibitem

\bibitem[\protect\citeauthoryear{Schwarz}{2006}]{schwarz2006hodge}
\begin{botherref}
\oauthor{\bsnm{Schwarz}, \binits{G.}}:
Hodge decomposition - {A} method for solving boundary value problems.
Springer
(2006)
\end{botherref}
\endbibitem

\bibitem[\protect\citeauthoryear{Morrey}{1956}]{morrey1956variational}
\begin{barticle}
\bauthor{\bsnm{Morrey}, \binits{C.B.}}:
\batitle{A variational method in the theory of harmonic integrals, ii}.
\bjtitle{American Journal of Mathematics}
\bvolume{78}(\bissue{1}),
\bfpage{137}--\blpage{170}
(\byear{1956})
\end{barticle}
\endbibitem

\bibitem[\protect\citeauthoryear{Zhao et~al.}{2019}]{zhao20193d}
\begin{barticle}
\bauthor{\bsnm{Zhao}, \binits{R.}},
\bauthor{\bsnm{Desbrun}, \binits{M.}},
\bauthor{\bsnm{Wei}, \binits{G.-W.}},
\bauthor{\bsnm{Tong}, \binits{Y.}}:
\batitle{3d hodge decompositions of edge-and face-based vector fields}.
\bjtitle{ACM Transactions on Graphics (TOG)}
\bvolume{38}(\bissue{6}),
\bfpage{1}--\blpage{13}
(\byear{2019})
\end{barticle}
\endbibitem

\bibitem[\protect\citeauthoryear{Friedrichs}{1955}]{friedrichs1955differential}
\begin{barticle}
\bauthor{\bsnm{Friedrichs}, \binits{K.O.}}:
\batitle{Differential forms on riemannian manifolds}.
\bjtitle{Communications on Pure and Applied Mathematics}
\bvolume{8}(\bissue{4}),
\bfpage{551}--\blpage{590}
(\byear{1955})
\end{barticle}
\endbibitem

\bibitem[\protect\citeauthoryear{Shonkwiler}{2009}]{shonkwiler2009poincare}
\begin{botherref}
\oauthor{\bsnm{Shonkwiler}, \binits{C.}}:
Poincar{\'e} duality angles on riemannian manifolds with boundary.
PhD thesis,
University of Pennsylvania
(2009)
\end{botherref}
\endbibitem

\bibitem[\protect\citeauthoryear{Su et~al.}{2024}]{su2024hodge}
\begin{barticle}
\bauthor{\bsnm{Su}, \binits{Z.}},
\bauthor{\bsnm{Tong}, \binits{Y.}},
\bauthor{\bsnm{Wei}, \binits{G.-W.}}:
\batitle{Hodge decomposition of single-cell rna velocity}.
\bjtitle{Journal of chemical information and modeling}
\bvolume{64}(\bissue{8}),
\bfpage{3558}--\blpage{3568}
(\byear{2024})
\end{barticle}
\endbibitem

\bibitem[\protect\citeauthoryear{Chen et~al.}{2021}]{chen2019evolutionary}
\begin{barticle}
\bauthor{\bsnm{Chen}, \binits{J.}},
\bauthor{\bsnm{Zhao}, \binits{R.}},
\bauthor{\bsnm{Tong}, \binits{Y.}},
\bauthor{\bsnm{Wei}, \binits{G.-W.}}:
\batitle{Evolutionary de rham-hodge method}.
\bjtitle{iscrete and Continuous Dynamical Systems - B}
\bvolume{26}(\bissue{7}),
\bfpage{3785}--\blpage{3821}
(\byear{2021})
\end{barticle}
\endbibitem

\bibitem[\protect\citeauthoryear{Liu et~al.}{2017}]{liu2017forging}
\begin{barticle}
\bauthor{\bsnm{Liu}, \binits{Z.}},
\bauthor{\bsnm{Su}, \binits{M.}},
\bauthor{\bsnm{Han}, \binits{L.}},
\bauthor{\bsnm{Liu}, \binits{J.}},
\bauthor{\bsnm{Yang}, \binits{Q.}},
\bauthor{\bsnm{Li}, \binits{Y.}},
\bauthor{\bsnm{Wang}, \binits{R.}}:
\batitle{Forging the basis for developing protein--ligand interaction scoring
  functions}.
\bjtitle{Accounts of chemical research}
\bvolume{50}(\bissue{2}),
\bfpage{302}--\blpage{309}
(\byear{2017})
\end{barticle}
\endbibitem

\bibitem[\protect\citeauthoryear{Cang and Wei}{2018}]{cang2018integration}
\begin{barticle}
\bauthor{\bsnm{Cang}, \binits{Z.}},
\bauthor{\bsnm{Wei}, \binits{G.-W.}}:
\batitle{Integration of element specific persistent homology and machine
  learning for protein-ligand binding affinity prediction}.
\bjtitle{International journal for numerical methods in biomedical engineering}
\bvolume{34}(\bissue{2}),
\bfpage{2914}
(\byear{2018})
\end{barticle}
\endbibitem

\bibitem[\protect\citeauthoryear{Cang et~al.}{2018}]{cang2018representability}
\begin{barticle}
\bauthor{\bsnm{Cang}, \binits{Z.}},
\bauthor{\bsnm{Mu}, \binits{L.}},
\bauthor{\bsnm{Wei}, \binits{G.-W.}}:
\batitle{Representability of algebraic topology for biomolecules in machine
  learning based scoring and virtual screening}.
\bjtitle{PLoS computational biology}
\bvolume{14}(\bissue{1}),
\bfpage{1005929}
(\byear{2018})
\end{barticle}
\endbibitem

\bibitem[\protect\citeauthoryear{Su et~al.}{2018}]{su2018comparative}
\begin{barticle}
\bauthor{\bsnm{Su}, \binits{M.}},
\bauthor{\bsnm{Yang}, \binits{Q.}},
\bauthor{\bsnm{Du}, \binits{Y.}},
\bauthor{\bsnm{Feng}, \binits{G.}},
\bauthor{\bsnm{Liu}, \binits{Z.}},
\bauthor{\bsnm{Li}, \binits{Y.}},
\bauthor{\bsnm{Wang}, \binits{R.}}:
\batitle{Comparative assessment of scoring functions: the casf-2016 update}.
\bjtitle{Journal of chemical information and modeling}
\bvolume{59}(\bissue{2}),
\bfpage{895}--\blpage{913}
(\byear{2018})
\end{barticle}
\endbibitem

\bibitem[\protect\citeauthoryear{Francoeur et~al.}{2020}]{francoeur2020three}
\begin{barticle}
\bauthor{\bsnm{Francoeur}, \binits{P.G.}},
\bauthor{\bsnm{Masuda}, \binits{T.}},
\bauthor{\bsnm{Sunseri}, \binits{J.}},
\bauthor{\bsnm{Jia}, \binits{A.}},
\bauthor{\bsnm{Iovanisci}, \binits{R.B.}},
\bauthor{\bsnm{Snyder}, \binits{I.}},
\bauthor{\bsnm{Koes}, \binits{D.R.}}:
\batitle{Three-dimensional convolutional neural networks and a cross-docked
  data set for structure-based drug design}.
\bjtitle{Journal of chemical information and modeling}
\bvolume{60}(\bissue{9}),
\bfpage{4200}--\blpage{4215}
(\byear{2020})
\end{barticle}
\endbibitem

\bibitem[\protect\citeauthoryear{Liu et~al.}{2023}]{liu2023persistent}
\begin{barticle}
\bauthor{\bsnm{Liu}, \binits{R.}},
\bauthor{\bsnm{Liu}, \binits{X.}},
\bauthor{\bsnm{Wu}, \binits{J.}}:
\batitle{Persistent path-spectral ({PPS}) based machine learning for
  protein--ligand binding affinity prediction}.
\bjtitle{Journal of Chemical Information and Modeling}
\bvolume{63}(\bissue{3}),
\bfpage{1066}--\blpage{1075}
(\byear{2023})
\end{barticle}
\endbibitem

\bibitem[\protect\citeauthoryear{Nguyen and Wei}{2019}]{nguyen2019dg}
\begin{barticle}
\bauthor{\bsnm{Nguyen}, \binits{D.D.}},
\bauthor{\bsnm{Wei}, \binits{G.-W.}}:
\batitle{Dg-gl: Differential geometry-based geometric learning of molecular
  datasets}.
\bjtitle{International journal for numerical methods in biomedical engineering}
\bvolume{35}(\bissue{3}),
\bfpage{3179}
(\byear{2019})
\end{barticle}
\endbibitem

\bibitem[\protect\citeauthoryear{Lin et~al.}{2022}]{lin2022language}
\begin{barticle}
\bauthor{\bsnm{Lin}, \binits{Z.}},
\bauthor{\bsnm{Akin}, \binits{H.}},
\bauthor{\bsnm{Rao}, \binits{R.}},
\bauthor{\bsnm{Hie}, \binits{B.}},
\bauthor{\bsnm{Zhu}, \binits{Z.}},
\bauthor{\bsnm{Lu}, \binits{W.}},
\bauthor{\bsnm{Santos~Costa}, \binits{A.}},
\bauthor{\bsnm{Fazel-Zarandi}, \binits{M.}},
\bauthor{\bsnm{Sercu}, \binits{T.}},
\bauthor{\bsnm{Candido}, \binits{S.}}, \betal:
\batitle{Language models of protein sequences at the scale of evolution enable
  accurate structure prediction}.
\bjtitle{BioRxiv}
\bvolume{2022},
\bfpage{500902}
(\byear{2022})
\end{barticle}
\endbibitem

\bibitem[\protect\citeauthoryear{Chen et~al.}{2021}]{chen2021extracting}
\begin{barticle}
\bauthor{\bsnm{Chen}, \binits{D.}},
\bauthor{\bsnm{Zheng}, \binits{J.}},
\bauthor{\bsnm{Wei}, \binits{G.-W.}},
\bauthor{\bsnm{Pan}, \binits{F.}}:
\batitle{Extracting predictive representations from hundreds of millions of
  molecules}.
\bjtitle{The journal of physical chemistry letters}
\bvolume{12}(\bissue{44}),
\bfpage{10793}--\blpage{10801}
(\byear{2021})
\end{barticle}
\endbibitem

\bibitem[\protect\citeauthoryear{Adams et~al.}{2014}]{adams2014javaplex}
\begin{bchapter}
\bauthor{\bsnm{Adams}, \binits{H.}},
\bauthor{\bsnm{Tausz}, \binits{A.}},
\bauthor{\bsnm{Vejdemo-Johansson}, \binits{M.}}:
\bctitle{Javaplex: A research software package for persistent (co) homology}.
In: \bbtitle{Mathematical Software--ICMS 2014: 4th International Congress,
  Seoul, South Korea, August 5-9, 2014. Proceedings 4},
pp. \bfpage{129}--\blpage{136}
(\byear{2014}).
\bcomment{Springer}
\end{bchapter}
\endbibitem

\bibitem[\protect\citeauthoryear{Bauer}{2021}]{bauer2021ripser}
\begin{barticle}
\bauthor{\bsnm{Bauer}, \binits{U.}}:
\batitle{Ripser: efficient computation of vietoris--rips persistence barcodes}.
\bjtitle{Journal of Applied and Computational Topology}
\bvolume{5}(\bissue{3}),
\bfpage{391}--\blpage{423}
(\byear{2021})
\end{barticle}
\endbibitem

\end{thebibliography}

\end{document}